\documentclass[final]{article}

\usepackage{graphicx}
\usepackage{amsfonts}
\usepackage{amsmath}
\usepackage{amssymb}
\usepackage{fancyhdr}
\usepackage{titlesec}
\usepackage{indentfirst}
\usepackage{booktabs}
\usepackage{verbatim}
\usepackage{color}
\usepackage{latexsym}
\usepackage{amscd}
\usepackage{esvect}
\usepackage{graphicx}
\usepackage{upgreek}
\usepackage{subcaption}
\usepackage{caption}
\usepackage[page,header]{appendix}
\usepackage{titletoc}
\usepackage{mathrsfs}
\usepackage[numbers]{natbib}
\usepackage{float}
\usepackage{lipsum}
\usepackage[export]{adjustbox}

\numberwithin{equation}{section}
\newtheorem{theorem}{Theorem}[section]

\newtheorem{proposition}[theorem]{Proposition}

\newtheorem{remark}[theorem]{Remark}

\allowdisplaybreaks

\def\bz{{\bf{z}}}

\def\bk{{\bf{k}}}
\def\bl{{\bf{l}}}

\begin{document}

\title{Nonlinear Geometric Optics Based Multiscale Stochastic Galerkin Methods for Highly Oscillatory Transport Equations with Random Inputs
\footnote{This work was partially supported by NSF grants
 DMS-1522184 and DMS-1107291: RNMS KI-Net,
 and by the Office of the Vice Chancellor for Research and
Graduate Education at the University of Wisconsin-Madison with funding from the Wisconsin
Alumni Research Foundation.}}
\date{\today}

\author{Nicolas Crouseilles\footnote{Inria (IPSO team), IRMAR, University of Rennes 1, Rennes, France},
 Shi Jin\footnote{Department of Mathematics, University of Wisconsin-Madison, Madison, WI 53706, USA (sjin@wisc.edu).},
Mohammed Lemou\footnote{CNRS, IRMAR, Inria (IPSO team), University of Rennes 1, Rennes, France},
Liu Liu \footnote{Department of Mathematics, University of Wisconsin-Madison, Madison, WI 53706, USA (lliu@math.wisc.edu).}}

\maketitle
\abstract{
We develop generalized polynomial chaos (gPC) based stochastic Galerkin
(SG) methods for a class of highly oscillatory transport equations that arise
in semiclassical modeling of non-adiabatic quantum dynamics. These models
contain uncertainties, particularly in coefficients that correspond to the
potentials of the molecular system.  We  first focus on a highly oscillatory scalar model with random uncertainty. Our method is built upon
the nonlinear geometrical optics (NGO) based method, developed in
\cite{NGO} for numerical approximations of deterministic equations, which can
obtain accurate pointwise solution even without numerically resolving spatially and temporally the oscillations. With the random uncertainty, we show that
such a method has oscillatory higher order derivatives in the random space, thus requires a frequency
dependent discretization in the random space. We modify this method by
introducing a new ``time'' variable based on the phase, which is shown
to be non-oscillatory in the random space, based on which we develop
a gPC-SG method that can capture oscillations with the frequency-independent
time step, mesh size as well as the degree of polynomial chaos.
A similar approach is then extended to a semiclassical surface hopping
model system with a similar numerical conclusion.
Various numerical examples attest that these methods indeed capture accurately the solution statistics {\em pointwisely} even though none of the numerical parameters resolve the high frequencies of the solution. }


\section{Introduction}
\noindent

Computational high frequency waves is challenging since one needs to
numerically resolve the small wave length which is often prohibitively
expensive \cite{ER2}.
Recently, a nonlinear geometric optics (NGO) based numerical method was
introduced in \cite{NGO} for a class of highly oscillatory transport
equations which allows the use of mesh size and time step {\it independent} of
the wave length. The transport equations that were solved by this method
 are deterministic and are relevant
to semiclassical approximations to quantum dynamics with band-crossings
\cite{CJL, CQJ, Morandi2, Morandi},
a general non-adiabatic quantum mechanics phenomenon that can be found
in a variety of important physics and chemistry problems, such as chemical
reaction \cite{Tully1, Tully2}, Bose-Einstein condensation \cite{WuNiu}, and graphene
\cite{Geim, FW}.

Since these semiclassical models, and many other quantum models, use
potential matrices that are obtained in empirical or {\it ad hoc} ways,
they inevitably contain {\it uncertainties}. Uncertainty could arise
in the potential, boundary or initial data, and forcing terms. In this paper,
we are interested in the random uncertainty of the potential energy surface, and
the goal is to develop efficient computational methods to compute the
propagation of the uncertainty. To deal with uncertainty in general, the so-called gPC-SG methods (combination of generalized polynomial chaos (gPC) approximation with stochastic Galerkin (SG) projections) are known to be
efficient for a wide range of partial differential equations with random uncertainties (see for example
\cite{GS, GWZ, LMK, Xiubook, XK}). However, direct application of such methods to highly-oscillatory problems with uncertainty becomes computationally expensive
if one wants to correctly capture the effect of these oscillations, since one
needs to resolve numerically the oscillations.
In this paper, we are interested in extending
the nonlinear geometric optics (NGO)-based method, which was developed in \cite{NGO} for a class of highly-oscillatory deterministic problems,
to highly-oscillatory problems with  uncertainty, in the framework of
gPC-SG methods. In particular, we will develop here
 gPC-SG methods that not only allow the mesh size and time step, {\it but also
the order of the gPC approximation, to be independent of the small wave length}.

We develop the gPC-SG methods for two model equations: a scalar model transport
equation and a semiclassical surface hopping model developed in \cite{CQJ}.
Both models were studied in \cite{NGO} in the deterministic case. Here we
assume random coefficients in these models that correspond to the band-gap in the non-adiabatic
quantum dynamics, and the band gap could become small to model the so-called
avoided crossing in which the quantum transition between bands is significant.
With the random uncertainty, we first show that
such a method is oscillatory for higher order derivatives in the random space, thus requires a frequency
dependent discretization in the random space, which becomes prohibitively
expensive to compute. We then modify this method by
introducing a new ``time'' variable based on the phase, which is shown
to be non-oscillatory in the random space, based on which we develop
a gPC-SG method that can capture the pointwise solution with {\it frequency-independent}
time step, mesh size as well as the degree of polynomial chaos.
This method is then extended to the semiclassical model of surface hopping
with random band gap, with the same numerical property.

The paper is organized as follows. In section 2, using a scalar equation with uncertain random coefficient,
we develop the gPC-SG method by either solving the equation directly,
or using the NGO approach. We prove that in both cases the method
needs to use wave frequency dependent gPC order (although for the latter method
such a dependence is one order more milder). Based on a theoretical result of
\cite{NGO}, we introduce a new time variable using the phase, which will
be shown to be non-oscillatory also in the random space, allowing us
the develop a new NGO-based method which are capable of obtaining accurate
pointwise solution without resolving the oscillations by any of the numerical
parameters. Numerical examples will verify the theoretical property as well
as the aforementioned numerical properties.
This method is then extended to a semiclassical
surface hopping model in section 3, with the same numerical properties which
are demonstrated numerically in section 4. We conclude the paper in section 5.

\section{The one dimensional scalar equation with random inputs}
\label{scalar}
In this section we consider a one dimensional problem with random inputs. Here the unknown is $u(t,x,\bz)\in\mathbb C$, where
$x\in\Omega\subset \mathbb R$, $t\geq 0$, and $\bz\in I_{\bz}\subset \mathbb R ^n, n\geq 1,$ is the random variable with a prescribed probability density function $\pi(\bz)$. We consider the following scalar model for $u$:
 \begin{equation}
\partial_t u+c(x)\partial_x u+r(u)=\frac{i a(x,\bz)}{\varepsilon}u, \qquad u(0,x,\bz)=u_{in}(x,\bz),
\label{model}
\end{equation}
where $a(x,\bz)$, $c(x)$, $r(u)$ and $u_{in}(x,\bz)$ are all given functions, with $a\geq 0$. We assume that $\Omega$ is a bounded interval and that periodic boundary
conditions are considered for the space variable $x$.

\subsection{A gPC-SG framework}
We briefly describe the SG method. Let $\mathbb P_P^n$ be the space of the $n$-variate polynomials of degree less than or equal to $P$,
$P\geq 1$, and recall that
$$ \text{dim}(\mathbb P_P^n)= \mbox{card}\{\bk \in \mathbb N ^n, |\bk|\leq P\}= \left(\begin{array}{c} n+P \\ P \end{array}  \right):=K,$$
where we have denoted $\bk= (k_1,\dots, k_n)$ and $|\bk|=k_1+\dots+k_n$.
We consider the following inner product
$$  \langle f, g\rangle_\pi =  \int_{I_{\bz}} f(\bz)g(\bz)\pi(\bz)d\bz, \quad \forall f, g\in L^2(\pi(\bz)d\bz),$$
where $L^2(\pi(\bz)d\bz)$ is the usual weighted Lebesgue space, and its associated norm is
$$ \|f\|_{L^2(\pi(\bz)d\bz)}^2 = \int_{I_{\bz}} |f(\bz)|^2\pi(\bz)d\bz.$$
Consider a corresponding orthonormal basis $\{\psi_{\bk}(\bz)\}_{\bk\in \mathbb N ^n, |\bk|\leq P}$ of  the space $\mathbb P_P^n$,  where the degree of $\psi_\bk$ is $\mbox{deg}(\psi_\bk)= |\bk|$. In particular
\begin{equation*}
\langle \psi_\bk, \psi_\bl\rangle_\pi = \int_{I_{\bz}} \psi_\bk(\bz)\psi_\bl(\bz)\pi(\bz)d\bz=\delta_{\bk\bl}, \qquad |\bk|, |\bl|\leq P,
\end{equation*}
where $\delta_{\bk\bl}$ is the Kronecker symbol. The commonly used pairs of $\{\psi_\bk(\bz)\}$ and $\pi(\bz)$ include Hermite-Gaussian, Legendre-uniform, Laguerre-Gamma, etc.  Since  the family $\{\psi_\bk(\bz)\}_{|\bk|\leq P}$ has $K$ elements we introduce its renumbering  family $( \tilde \psi_1,\tilde \psi_2, \cdots \tilde \psi_K)$, that is
$$ \{\psi_\bk(\bz)\}_{\bk\in \mathbb N^n, |\bk|\leq P}= ( \tilde \psi_1,\tilde \psi_2, \cdots \tilde \psi_K),\quad \mbox{deg}(\tilde \psi_j)\leq \mbox{deg}(\tilde \psi_{j+1}).$$
The SG method consists in seeking the solution to \eqref{model} as a projection onto the space $\mathbb P_P^n$, that is
\begin{equation}u(t,x,\bz)\approx \sum_{|\bk|\leq P}u_\bk(t,x)\psi_\bk(\bz)
=\sum_{j=1}^K \tilde u_j(t,x) \tilde \psi_j(\bz), \label{soln u}\end{equation}
with
$$u_\bk = \langle u, \psi_\bk\rangle_\pi, \qquad  \tilde u_j = \langle  u, \tilde \psi_j\rangle_\pi, \qquad |\bk|\leq P,\quad  \  j=1,2, ..., K.$$
From this approximation one can easily compute statistical moments, such as
the mean and standard deviation, as
\begin{equation}
\mathbb{E}(u)\approx u_{(0,\cdots, 0)}=\tilde u_1, \quad
\text{SD}(u)\approx \left( \sum_{1\leq |\bk|\leq P}| u_{\bk}|^2\right)^{1/2}=  \left( \sum_{j=2}^K |\tilde u_{j}|^2\right)^{1/2}.
\label{moments}
\end{equation}

\subsection{The direct gPC-SG method}
\label{gPC_Dir_1D}
We first introduce the gPC-SG method for the direct method (hereafter
called gPC-SG-D)--solving (\ref{model}) directly.
For simplicity of illustration, assume $c(x)$ does not depend on $\bz$.
The case when $c$ depends on $\bz$ can be easily incorporated into the gPC-SG
framework \cite{GX}. The gPC-SG solution for $u$ is computed through the projection formula \eqref{soln u} as follows.

Denote the gPC coefficients $$\displaystyle \vv u(t,x)=(\tilde{u}_1(t,x) \cdots, \tilde{u}_K(t,x))^{T}.$$
Inserting the approximated solution (\ref{soln u}) to (\ref{model}) and conducting the standard Galerkin projection, one gets
\begin{equation}\partial_t \vv u+c(x)\partial_x\vv u+\vv{\boldsymbol{\gamma}}(u)=\frac{i}{\varepsilon}A(x)\vv u, \qquad
\vv u(0,x)=\int_{I_{\bz}}u_{in}(x,\bz)\tilde{\boldsymbol{\psi}}(\bz)\pi(\bz)d\bz\,,
\label{V_Dir}
\end{equation}
with
\begin{equation}
\tilde{\boldsymbol{\psi}}(\bz)= (\tilde \psi_1(z),\tilde \psi_2(z), \cdots,\tilde \psi_K(z))^T, \quad u(t,x,\bz)= \sum_{j=1}^K \tilde u_j(t,x)\tilde \psi_j(\bz). \label{psi}
\end{equation}
The $j$-th component of the vector $\vv{\boldsymbol{\gamma}}= (\tilde \gamma_1, \tilde \gamma_2,\cdots, \tilde \gamma_K)$ is given by
\begin{equation}
{\tilde \gamma}_{j}(u)=\int_{I_{\bz}}r(u)\tilde \psi_{j}(\bz)\pi(\bz)d\bz\,, \qquad\text{for}\,\ j=1,\cdots, K,
\label{gamma0}
\end{equation}
and the symmetric, non-negative definite matrix $A\in\mathbb R_{K\times K}$  (recall that $a\geq 0$) is defined by
\begin{equation}
A_{sj}(x)=\int_{I_{\bz}}a(x,\bz)\tilde\psi_{s}(\bz)\tilde \psi_{j}(\bz)\pi(\bz)d\bz\,, \qquad\text{for}\,\ s, j=1, \cdots, K.
\label{matrix A}
\end{equation}

To solve numerically \eqref{V_Dir}, we use a simple time splitting.
To this aim, we fix a time step $\Delta t >0$ and set $t_n=n\Delta t $ for $n\in \mathbb N$.
As usual, we denote by $\vv u^{n}(x)= \left( \tilde u^{n}_1(x),\tilde u^{n}_2(x),\cdots,\tilde  u^{n}_K(x)\right)$ an approximation of $\vv u(t_n,x)$. We split system (\ref{V_Dir}) in three parts which can be solved as follows. \\[2pt]
\begin{itemize}
\item {\bf Oscillatory part}
\end{itemize}
The oscillatory part $$\partial_t\vv u=\frac{i}{\varepsilon}A(x)\vv u, $$
 is solved exactly in time, for each $x$, to get
$$\vv u^{n+1}(x)=\exp\left(\frac{i}{\varepsilon}A(x)\Delta t\right)\vv u^n(x)\,. $$

\newpage
\begin{itemize}
\item {\bf Nonlinear part}
\end{itemize}
To solve the  nonlinear part
$$\partial_t\vv u+\vv{\boldsymbol{\gamma}}(u)=0, $$
 we choose to use the forward Euler method, that is
\begin{equation}\vv u^{n+1}=\vv u^n -\Delta t \vv{\boldsymbol{\gamma}}(u^n)\,.
\label{source}
\end{equation}
To compute $\vv{\boldsymbol{\gamma}}(u^n)$, we first compute $\displaystyle u^n(x, \bz)=\sum_{s=1}^{K}\tilde u_{s}^n(x)\tilde \psi_{s}(\bz)$ using
the Gauss-quadrature rule (where the quadrature points are chosen as the roots of the orthogonal polynomials determined by
the distribution $\pi$ of the random variables, see \cite{XK})
to get $\vv{\boldsymbol{\gamma}}(u)$ given by (\ref{gamma0}) as follows. Define the Gauss-quadrature points as $z^{(l)}$ and the corresponding weights $\omega_l,$ $l=1, \cdots, N_g$.
The value of $r(u)$ evaluated at $x$, $t$ and at  the Gauss-quadrature point $z^{(l)}$ is given by
$\displaystyle r(u(z^{(l)}))=r\left(\sum_{j=1}^{K}\tilde u_{j}(t,x)\tilde \psi_{j}(z^{(l)})\right)$.
By (\ref{gamma0}), one has
$$\tilde \gamma_{j}(u)=\sum_{l=1}^{N_g} r(u(z^{(l)}))\tilde \psi_{j}(z^{(l)})\omega_l.$$
which is used in (\ref{source}) to get $\vv u^{n+1}$.   \\[2pt]

\begin{itemize}
\item {\bf Transport part}
\end{itemize}
To solve the transport part $$\partial_t\vv u+c(x)\partial_x\vv u=0, $$
we  use a pseudo-spectral method in space and a three-stage Runge-Kutta method in time which was introduced in \cite{Choi-Liu}.
Denote $\mathcal T (u^n)=-c(x)\mathcal F^{-1}(i\xi \mathcal F(u^n))$, where $\mathcal F$ and $\mathcal F^{-1}$ are the (discrete) Fourier and inverse Fourier
Transforms respectively, and $\xi$ is the Fourier space variable. Then
\begin{align}
\displaystyle
&\nonumber\vv u^{n,(1)}=\vv u^n + \frac{1}{2}\Delta t \mathcal T(\vv u^n), \\[4pt]
&\label{RK1} \vv u^{n,(2)}=\vv u^n + \frac{1}{2}\Delta t \mathcal T(\vv u^{n,(1)}), \\[4pt]
&\nonumber\vv u^{n+1}=\vv u^n+\Delta t\mathcal T(\vv u^{n,(2)}).
\end{align}
Notice that the spectral approximation is central type finite difference
approximation with purely imaginary spectrum, thus one needs to use ODE solvers
that have a stability region that contains part of the purely {\it imaginary}
axis \cite{Choi-Liu}. The scheme (\ref{RK1}) is such an ODE solver where the stability region
takes the largest part of the imaginary axis among three-stage
ODE solvers, and is of second order accuracy in time.

\subsection{The NGO-based gPC-SG method}
We first review the NGO based method introduced in \cite{NGO} for the deterministic one-dimensional scalar equation,
\begin{equation}
\partial_t u+c(x)\partial_x u+r(u)=\frac{i a(x)}{\varepsilon}u, \qquad u(0,x)=u_{in}(x),
\label{model2}
\end{equation}
where the functions $u_{in}$, $a$, $c$ and $r$ are given and $a(x)\geq 0$.
Periodic boundary conditions in space are considered.  For the sake of simplicity, we assume that the initial data $u_{in}$ does not depend on $\varepsilon$,
namely, it is non-oscillatory.


We focus on the case when $r$ is nonlinear. We recall how $\it{nonlinear\, geometric\, optics}$ (NGO) is utilized to solve problem (\ref{model2}) in \cite{NGO}. Introduce a profile function $U(t,x,\tau)$ which depends on the $2\pi\hspace{-0.1cm}-\hspace{-0.1cm}\it{periodic}$ variable $\tau$, and a phase $S(t,x)$ such that
\begin{equation}
\label{USu}
U(t,x,S(t,x)/\varepsilon)=u(t,x),
\end{equation}
with $u$ solving (\ref{model2}). Inserting this ansatz into (\ref{model2}) gives
$$\partial_t U+c(x)\partial_x U+\frac{1}{\varepsilon}\left[\partial_t S+c(x)\partial_x S\right]\partial_{\tau}U+r(U)=\frac{ia(x)}{\varepsilon}U.$$
Due to the periodicity constraint on $U$, the following equation on the phase $S$ should be imposed (see \cite{NGO} for details),
\begin{equation}\partial_t S+c(x)\partial_x S=a(x), \qquad S(0,x)=0,
\label{Seqn}
\end{equation}
then the equation for $U$ is given by
\begin{equation}\partial_t U+c(x)\partial_x U+r(U)=-\frac{a(x)}{\varepsilon}(\partial_{\tau}U-iU), \qquad U(0,x,0)=u_{in}(x).
\label{U}
\end{equation}
For convenience, we write  \eqref{U} in terms of  $V=e^{-i\tau}U$ and get
\begin{equation}
\partial_t V+c(x)\partial_x V+e^{-i\tau}r(e^{i\tau}V)=-\frac{a(x)}{\varepsilon}\partial_{\tau}V, \qquad V(0,x,0)=u_{in}(x).
\label{Veqn}
\end{equation}


One needs initial data $V(0,x,\tau)$ for all $\tau$ to solve equation (\ref{Veqn}). The only requirement we have to ensure is $V(0,x,0)=u_{in}(x)$.
One critical idea is to use initial data for $V$ so that the solution to (\ref{Veqn}) remains bounded uniformly in $\varepsilon$ up to certain order of derivatives.

Introduce the operators $\mathcal L$ and $\Pi$,
$$\displaystyle\mathcal Lg=\partial_{\tau}g, \qquad \Pi g=\frac{1}{2\pi}\int_0^{2\pi}g(\tau)d\tau,$$
and $$\displaystyle\mathcal L^{-1}g=(\mathcal I-\Pi)\int_0^{\tau}g(\sigma)d\sigma=\int_0^{\tau}g(\sigma)d\sigma+\frac{1}{2\pi}\int_0^{2\pi}\sigma g(\sigma)d\sigma.$$
Following the work of \cite{PCL, FCL, NGO}, a Chapman-Enskog expansion is used to give a suitable initial condition given by
\begin{align}
& \label{IC_V} V(0,x,\tau)=u_{in}(x)+\frac{\varepsilon}{a(x)}\left[ {\cal G}(0,u_{in})-{\cal G}(\tau,u_{in})\right]\,, \\[4pt]
& \quad\text{with}\,\ {\cal G}(\tau,u_{in})=\mathcal L^{-1}(\mathcal I-\Pi)\left[ e^{-i\tau}r(e^{i\tau}u_{in})\right]\,.
\end{align}
It was proved in \cite{NGO} that the solution $V$ to (\ref{Veqn}) and (\ref{IC_V}) has bounded ({\em uniformly in $\varepsilon$}) derivatives in both $x$ and $t$ up to second order. Thus it allows the construction of a scheme
with a uniform accuracy with respect to $\varepsilon$. We recall in the sequel the scheme introduced in \cite{NGO}.
\\[2pt]

Consider a uniform partition in time $t_n=n\Delta t$ ($\Delta t>0$ the time step) of a time interval $[0,T]$, $n=0,1,\cdots, N$, $N\Delta t=t_N$ and in space $x_j=j\Delta x$,
$j=0,1,\cdots, N_x$, $\Delta x=1/N_x$ of the spatial interval $[0,1]$. A uniform mesh is assumed in the direction
$\tau\in\mathbb T=[0, 2\pi]$,
$\tau_l=l\Delta\tau$, $l=0, \cdots, N_{\tau}$, $\Delta\tau=2\pi/N_{\tau}$. Denote by $V_j^n(\tau)\approx V(t_n, x_j,\tau)$ and $S_j^n\approx S(t_n, x_j)$ the discrete unknowns at time $t_n$,
evaluated at $x_j$.

To solve the equation (\ref{Seqn}) for $S$, the pseudo-spectral method in space and $4$-th order Runge-Kutta method in time are used. In other words, the
following ODE system on $S_j(t) \approx S(t, j\Delta x)$:
\begin{align}
\label{eqn S} \partial_t S_j + c(x_j)\mathcal F^{-1}(i\xi \mathcal F(S^n))_j = a(x_j), \quad S_j(0)=0,
\end{align}
is solved by the $4$-th order Runge-Kutta method. To recover the original solution $u$ at the final time $t_N$,
we use \eqref{USu} at $t=t_N$ and $x=x_j$.
Since $S_j^N/\varepsilon$ does not necessarily coincide with a grid point $\tau_l$ at the final time $t_N=N\Delta t$,
thus one can perform a trigonometric interpolation
\begin{equation}
u(t_N, x_j)\approx V_j^N(\tau=S_j^N/\varepsilon).
\label{reconstruct}
\end{equation}
Note that a higher order method is necessary here to solve $S$ since we need to construct the quantity $S(t,x)/\varepsilon$, and the error of $S$ is divided by $\varepsilon$.

 To solve (\ref{Veqn}) for $V$, we start with $V_j^0$ given by (\ref{IC_V}) at $x=x_j$.
 Then, the scheme  for $V$ reads (assume $c(x)>0$)
\begin{align}
&\label{eqn V}\frac{V_j^{n+1}- V_j^{n}}{\Delta t}+c(x_j) \frac{V_j^n- V_{j-1}^n}{\Delta x}+e^{-i\tau}r(e^{i\tau}V_j^n)=-\frac{a(x_j)}{\varepsilon}\partial_{\tau}V_j^{n+1}.
\end{align}
Here the right hand side is discretized implicitly due to its numerical stiffness and the $\tau$ variable
is discretized using the Fourier transform.
\subsubsection{The gPC-SG-N1 method}
\label{gPC1D_N}
Our aim now is to extend the NGO method presented in the previous section to solve a highly oscillatory problem with uncertainty, namely \eqref{model}.
By \cite{NGO}, and similarly to (\ref{Seqn}) and (\ref{Veqn}), the NGO-based method solves $S$ from
\begin{equation}\partial_t S+c(x)\partial_x S=a(x,\bz), \qquad S(0,x,\bz)=0,
\label{Seqn1}
\end{equation}
and solves $V$ from
\begin{equation}
\partial_t V+c(x)\partial_x V+e^{-i\tau}r(e^{i\tau}V)=-\frac{a(x,\bz)}{\varepsilon}\partial_{\tau}V,
\label{Veqn1}
\end{equation}
with initial data
\begin{align}
\displaystyle
&\label{IC_1D} V(0,x,\tau,\bz)=u_{in}(x,\bz)+\frac{\varepsilon}{a(x,\bz)}\left[{\cal G}(0,u_{in},\bz)-{\cal G}(\tau,u_{in},\bz)\right]\\[4pt]
&\label{IC_1D_G} \,\ \text{with} \,\ {\cal G}(\tau,u_{in},\bz)=\mathcal L^{-1}(\mathcal I-\Pi)\left[e^{-i\tau}r(e^{i\tau}u_{in})\right],
\end{align}
where $u_{in}(x,\bz)$ and $a(x,\bz)$ now depend on $\bz$. Below we detail the gPC approach for this system.

\paragraph{{The gPC formulation for $S$\\}}
One inserts the gPC-SG ansatz
$$S(t,x,\bz)\approx\sum_{ \ |\bk|\leq P}S_{\bk}(t,x)\psi_{\bk}(\bz)= \sum_{j=1}^{K}\tilde S_{j}(t,x)\tilde \psi_{j}(\bz)$$
into (\ref{Seqn1}) and conducts the Galerkin projection, to get
\begin{align}
\partial_t \vv S+c(x)\partial_x\vv S=\vv{\mathcal R}(x), \qquad \vv S(0,x)=0,
\label{S_gPC}
\end{align}
where the gPC coefficients $\vv S$ is defined by $$\displaystyle \vv S(t,x)=(\tilde{S}_1(t,x), \cdots, \tilde{S}_K(t,x))^{T},$$
and
\begin{equation}
\label{RR}
\vv{\mathcal R}(x)=\int_{I_{\bz}}a(x,\bz) \tilde{\boldsymbol{\psi}}(\bz)\pi(\bz)d\bz.
\end{equation}
We recall the notation:  $ \tilde{\boldsymbol{\psi}}= (\tilde\psi_1, \tilde\psi_2, \cdots, \tilde \psi_K)$.

\paragraph{The gPC formulation for $V$\\}
We insert the gPC-SG ansatz
$$V(t,x,\tau,\bz)\approx\sum_{ |\bk|\leq P}V_{\bk}(t,x,\tau)\psi_{\bk}(\bz)=\sum_{j=1}^{K}\tilde V_{j}(t,x,\tau)\tilde \psi_{j}(\bz)$$
into (\ref{Veqn1}) and conduct the Galerkin projection to get
\begin{equation}
\displaystyle\partial_t\vv V+c(x)\partial_x\vv V+e^{-i\tau}\vv{\boldsymbol{\gamma}}(e^{i\tau}V)=-\frac{1}{\varepsilon}A(x)\partial_{\tau}\vv V,
\end{equation}
where $A$ is the matrix given by \eqref{matrix A} and  $$\displaystyle \vv V(t,x,\tau)=(\tilde{V}_1(t,x,\tau), \cdots, \tilde{V}_K(t,x,\tau))^{T}.$$
The $j$-th component $\tilde \gamma_j$ of the vector $\vv{\boldsymbol{\gamma}}$ is
\begin{equation}
\tilde \gamma_j(e^{i\tau}V)=\int_{I_{\bz}}r(e^{i\tau}V)\tilde \psi_{j}(\bz)\pi(\bz)d\bz, \qquad\text{for}\,\ j=1,\cdots, K,
\label{gamma1}
\end{equation}
which is computed by the Gauss-quadrature formula. \\[2pt]

\paragraph{The initial data\\}
Here we consider that $a(x,\bz)$ depends on the random variable $\bz$.
The gPC approximation of the initial condition (\ref{IC_1D}) becomes
\begin{align}
\displaystyle\vv V(0,x,\tau)=\boldsymbol T+\varepsilon\boldsymbol Y,
\label{IC}
\end{align}
where the vectors $\boldsymbol T\in\mathbb R^{K}$ and $\boldsymbol Y\in\mathbb R^{K}$ are given by
\begin{align*}
&\displaystyle\boldsymbol T(x)=\int_{I_{\bz}}u_{in}(x,\bz)\tilde{\boldsymbol{\psi}}(\bz)\pi(\bz)d\bz, \\[2pt]
&\displaystyle\boldsymbol Y(x, \tau)=\int_{I_{\bz}}\frac{1}{a(x,\bz)}\left[{\cal G}(0,u_{in},\bz)-{\cal G}(\tau,u_{in},\bz)\right]\tilde{\boldsymbol{\psi}}(\bz)\pi(\bz)d\bz,
\end{align*}
with $\tilde{\boldsymbol{\psi}}(\bz)$ given in (\ref{psi}).

\subsubsection{The fully discrete gPC-SG-N1 method}
\label{gPC-SG-N}
Now we discretize (\ref{Seqn1}), (\ref{Veqn1}) in time as in the deterministic case. The phase vector $\vv S$ is obtained by applying a pseudo-spectral method in space, which leads to
the following semi-discretized  equation
\begin{align}
\displaystyle
\label{S_gPC_N1} \partial_t \vv S_j + c(x_j)\mathcal F^{-1}(i\xi \mathcal F(\vv S))_j = \vv{\mathcal R}(x_j), \quad \vv S_j(0)=0\,, \ \ \forall
j=0,\cdots N_x,
\end{align}
where $\vv{\mathcal R}(x)$ is given by (\ref{RR}).  Here, $\vv S_j(t)$ is an approximation of $\vv S (t,x_j)$ and the same notation
is used for other quantities. Then a $4$-th order Runge-Kutta method is applied to march in time.

To solve the equation on $\vv V$, knowing
$\vv S$, we start by the initial value $\vv V_j^0$ given by (\ref{IC}) at grid points $x=x_j$. Then the system is advanced in time by
a simple time-splitting algorithm as explained for the deterministic case in section \ref{gPC_Dir_1D}. The nonlinear part (the $\vv \gamma$ term) and the transport part are treated in the same way as in section \ref{gPC_Dir_1D}. For the oscillatory part, we use the backward Euler method in time:
\begin{equation}
\displaystyle \frac{\vv{V_j}^{n+1}-\vv{V_j}^n}{\Delta t}=-\frac{1}{\varepsilon}A(x_j)\partial_{\tau}\vv{V_j}^{n+1},
\label{tau}
\end{equation}
where $\vv{V_j}^n(\tau)$ is an approximation of $\vv{V}(n\Delta t, x_j, \tau)$ and $A$ is given by \eqref{matrix A}.
Let
$\displaystyle\vv{\hat V}_j^n(\zeta)$ be the Fourier transform of $\vv{V_j}^n(\tau)$ in the periodic variable $\tau$, where $\zeta$ is the Fourier variable. The spectral method is used to discretize the $\tau$-derivative, then \eqref{tau} becomes (removing the dependency in $\zeta$)
\begin{equation}
\displaystyle\vv{\hat V_j}^{n+1}=\left[I+i\zeta \frac{\Delta t}{\varepsilon}A(x_j)\right]^{-1}\vv{\hat V_j}^n.
\end{equation}
The updated values $\vv V_j^{n+1}$ are then obtained by the inverse Fourier transform. Note that since $A$ is non-negative definite,
all its eigenvalues are real and non-negative. Thus the matrix $\displaystyle I+i\zeta \frac{\Delta t}{\varepsilon}A$ has non-zero eigenvalues, thus is invertible.
In the sequel, we refer this method as the gPC-SG-N1 method.
\subsection{The relation between the gPC order $K$ and $\varepsilon$}
\label{K_eps}

The main advantage of the NGO-based method over the direct method is that the
former allows one to use $\Delta t$ and $\Delta x$ {\it independent} of
$\varepsilon$, thus when $\varepsilon$ is small, one can still use relatively
larger $\Delta t$ and $\Delta x$ to get accurate solution at the mesh
points. In this section, we prove that this remarkable property is not
true for the random approximation, namely {\it one can not use the gPC order $K$
independent of $\varepsilon$} either for the gPC-SG-D or gPC-SG-N1.

\begin{itemize}
\item {\bf{The gPC-SG-D}}
\end{itemize}

To show the gPC-SG-D method does not satisfy the uniform accuracy property, we consider the simple case where
the functions $c$ and $a$ are independent of $x$, and where $\bz\in \mathbb R$:
$$\displaystyle\partial_t u+c\partial_x u= \frac{ia(\bz)}{\varepsilon}u\,.$$
By the method of characteristics, the analytic solution is given by
 $$\displaystyle u(t,x,\bz)=e^{\frac{ia(\bz)}{\varepsilon}t}u_{in}(x-ct,\bz), $$
where $u_{in}$ is the initial data.
Taking the $q$-th order derivative in $\bz$, one sees
that  $\displaystyle \partial_\bz^q u=O(1/\varepsilon^q)$.

Assume $u_{in}\in H_{\bz}^q(\pi(\bz)d\bz)$, where $H_{\bz}^q$ is the Sobolev space in $\bz$ defined as
$$||u||_{H_{\bz}^q}=\left(\sum_{l=0}^q ||\partial_{\bz}^l u||_{L^2(\pi(\bz)d\bz)}\right)^{\frac{1}{2}}.$$
By the standard approximation theory of
orthogonal polynomials (see for example \cite{STW}),
one has
$$\displaystyle ||u-\sum_{j=1}^K \tilde u_j \tilde \psi_j ||_{L^2(\pi(\bz)d\bz)}\leq\frac{C}{\varepsilon^{q}K^q}, $$
where $C$ is a constant depending on $t$, the derivatives of $a(\bz)$ and $u_{in}$ with respect to $\bz$.
This shows that one needs to take $K$ inversely proportional to
$\varepsilon$ in order to ensure an accurate approximation in $\bz$.
\\[2pt]
\begin{itemize}
\item {\bf{The gPC-SG-N1}}
\end{itemize}

We consider the gPC-SG-N1 method and show that the degree of the polynomial approximation needs to be chosen dependent on $\varepsilon$ to ensure a spectral accuracy.
To this aim, we analyze the projection error on the function $V$, solution to (\ref{Veqn1}), 
\begin{equation}
\label{V-eq}
\partial_t V + c(x)\partial_x V + e^{-i\tau}r(e^{i\tau}V) = -\frac{a(x,\bz)}{\varepsilon}\partial_{\tau}V \,.
\end{equation}

From (\ref{V-eq}), clearly, the $l$-th derivative in $z$ of $V$ satisfies basically the {\it same} equation
as the $l$-th derivative in $x$ of $V$ (except one involves the derivatives of
$a$ in $z$ while the other involves the derivatives of $a$ in  $x$). Therefore,
as long as $a$ and $u_{in}$ have the same regularity in $x$ and $z$, $V$ will have the {\it same} regularity in $x$ and $z$.  Since with the choice of initial
data (\ref{IC_1D}), we only gaurantee bounded (in terms of $\varepsilon$)
$x$ derivatives of $V$ up to
the second order, therefore, we only have up to second bounded derivatives in
$z$ for $V$.  This is clearly not enough for SG methods (nor for the stochastic
collocation method, to be introduced later) which can achieve high order--up to spectral--accuracy if there is sufficient regularity in the solution.

\subsection{A new NGO-based gPC-SG method (gPC-SG-N2)}
In the following, we show how to construct a numerical method for which the derivatives of $V$ with respect to $\bz$ and $x$
are bounded uniformly in $\varepsilon$ at {\it any} order.
To this aim, let us consider the following change of function
\begin{equation}
\label{relation_V_W}
V(t, x, \bz, \tau) = W(S(t, x, \bz), x, \bz, \tau),
\end{equation}
where $S$ solves \eqref{Seqn1}.
Then, $W$ is the solution to the following problem:
\begin{eqnarray}
&&\partial_s W+\frac{c}{a(x, \bz)}\partial_x W + \frac{1}{a(x, \bz)} e^{-i\tau} r(e^{i\tau}W) = -\frac{1}{\varepsilon}\partial_\tau W, \nonumber\\
&&W(0, x, \bz, \tau)= V(0,x, \bz,\tau)= u_{in}(x, \bz)+ \frac{\varepsilon}{a(x, \bz)}  \Big[ {\cal G}(0, u_{in},\bz) -{\cal G}(\tau, u_{in},\bz)\Big]\,, \nonumber\\
 &&\mbox{with } {\cal G}(\tau, u_{in},\bz) = {\cal L}^{-1} ({\cal I}-\Pi)[e^{-i\tau}r(e^{i\tau}u_{in}(x, \bz))]\,.
\label{eq_W}
\end{eqnarray}
We then have the following result for $W$.
\begin{proposition}
\label{theoremW}
Let $W$ be the solution of \eqref{eq_W} on $[0, \bar T]$, $\bar T>0$, with periodic boundary condition in $x$ and $\tau$.
Then, up to the second order derivative in $s$, and arbitrary order in $x$ and $\bz$ derivatives of $W$
are bounded uniformly in $\varepsilon\in [0,1]$, that is, $\exists C>0$ independent of $\varepsilon$ such that, $\forall s\in [0, \bar T]$
$$
\|\partial_s^p\partial_x^q\partial_\bz^r W(s)\|_{L^\infty_{\tau, x}(L^2(\pi(\bz)d\bz))} \leq C, \;\; \mbox{ for } \;\; p=0, 1, 2, \;\; \mbox{ and } \;\; q, r\in \mathbb{N}.
$$
\end{proposition}
This proposition can be deduced from Lemma 2.6 of \cite{NGO} in which we replace $x$ by $(x,\bz)$.
In fact, the $W$-equation (\ref{eq_W}) was introduced in \cite{NGO} to prove the uniform in $\varepsilon$ regularity and convergence of the
NGO-based method.
Note in particular that the high order derivatives in $x$ and $\bz$ are uniformly bounded at \textit{any} order with respect
to $\varepsilon$ although the initial data for $W$ is only prepared at the first order in $\varepsilon$.
This guarantees the uniform accuracy in $\bz$ and $x$, that is in particular
 $$
 \displaystyle ||W-\sum_{j=1}^K \tilde W_j \tilde \psi_j ||_{L^2(\pi(\bz)d\bz)}\leq\frac{C}{K^q}, \;\; \forall 
 q \text{  positive integer},
 $$
where $C$ is a positive constant which is independent of $\varepsilon$, and $q$ depends only on the regularity
in $\bz$ of the initial data $u_{in}(x,\bz)$. Thus one can choose $K$ \textit{independent of} $\varepsilon$.

To summarize, to obtain a spectral accuracy in $\bz$ at any order uniformly in $\varepsilon$, one has to
solve equation \eqref{eq_W} satisfied by $W$ (instead of the equation on $V$) and then recover $V$
by formula \eqref{relation_V_W}. Note that the mapping $t\mapsto S(t, x, \bz)$ is a increasing function in $t$
for any fixed $(x, \bz)$. Of course, this needs an interpolation in time to get back to the $t$ variable for $V$.

\subsubsection{The gPC-SG-N2 method}
The above discussion motivates us to design the following new scheme, denoted by gPC-SG-N2 in the sequel.

The gPC system for solving the equation for $S$ is the same as in section \ref{gPC1D_N}.
For the equation on $W$, we insert the gPC-SG ansatz
$$W(s,x,\tau,\bz)\approx\sum_{ |\bk|\leq P}W_{\bk}(s,x,\tau)\psi_{\bk}(\bz)=\sum_{j=1}^{K}\tilde W_{j}(s,x,\tau)\tilde \psi_{j}(\bz)$$
into (\ref{eq_W}) and conduct the Galerkin projection to get
\begin{equation}
\displaystyle\partial_s\vv W+c(x)A^{\star}(x)\partial_x\vv W+e^{-i\tau}\vv{\boldsymbol{\gamma^{\star}}}(e^{i\tau}W)=-\frac{1}{\varepsilon}\partial_{\tau}\vv W,
\end{equation}
where $A^{\star}$ is the symmetric and positive-definite matrix given by
$$A^{\star}_{ij}(x)=\int_{I_{\bz}}\frac{1}{a(x,\bz)}\psi_i(\bz)\psi_j(\bz)\pi(\bz)d\bz,$$
and $$\displaystyle \vv W(s,x,\tau)=(\tilde{W}_1(s,x,\tau), \cdots, \tilde{W}_K(s,x,\tau))^{T}.$$
The $j$-th component $\tilde\gamma^{\star}_j$ of the vector $\vv{\boldsymbol{\gamma^{\star}}}$ is ,
\begin{equation}
\tilde\gamma^{\star}_j(e^{i\tau}W)=\int_{I_{\bz}}\frac{1}{a(x,\bz)}r(e^{i\tau}W)\tilde\psi_{j}(\bz)\pi(\bz)d\bz, \qquad\text{for}\,\ j=1,\cdots, K,
\label{gamma1_N2}
\end{equation}
which is computed by the Gauss-quadrature formula.
The initial data for $W$ is the same as for $V$: $W(0,x,\tau,\bz)=V(0,x,\tau,\bz)$, which is shown in section \ref{gPC1D_N}.

\subsubsection{The fully discrete gPC-SG-N2 method}
\label{N2_1D}
Let the final time be $T$. Define the grid points in space $x_j$, $j=0,1,\cdots, N_x$. As discussed in (\ref{reconstruct}),
to recover the solution $u$, one needs a trigonometric interpolation on $\displaystyle\tau=\frac{S(T,x_j, \bz^{(l)})}{\varepsilon}$, for each quadrature points
$\bz^{(l)}$, $l=1, \cdots, N_s$,
$$u(T,x_j,\bz^{(l)})\approx V(T,x_j,\tau=\frac{S(T,x_j,\bz^{(l)})}{\varepsilon},\bz^{(l)}),$$
then one can compute the mean and standard deviation of $u$ at $(T,x_j)$ in the random space,
$$
\mathbb{E}(u)(T,x_j)=\sum_{l=1}^{N_s}u(T,x_j,\bz^{(l)})\omega_l\,, \quad
\text{SD}(u)(T,x_j)=\left(\sum_{l=1}^{N_s}|u(T,x_j,\bz^{(l)})|^2 \omega_l-\left(\mathbb{E}(u)(T,x_j)\right)^2\right)^{1/2}\,,
$$
where $\omega_l$ are the corresponding quadrature weights.
After obtaining the values of $S$ at $(T, x_j, \bz^{(l)})$ for all $j$, $l$, let $$S^{\star}=\max_{j,l}\left\{S(T, x_j, \bz^{(l)})\right\}. $$ We choose
$s\in[0, S^{\star}]$ to be the time domain for computing $W$.
Define the time discretization of $W$ as $s_l$, for $l=0,1, \cdots, L$, with time step $\Delta s$ and the total number of time steps $\displaystyle L=\left[\frac{S^{\star}}{\Delta s}\right]+1$.

To solve the equation on $\vv W$, knowing
$\vv S$, we start by the initial value $\vv W(0,x,\tau)=\vv V(0,x,\tau)$ given by (\ref{IC}) at $x=x_j$. Then the system is advanced in time by
a simple time-splitting algorithm as explained for the deterministic case in section \ref{gPC_Dir_1D}. The nonlinear part (the $\vv \gamma^{\star}$ term) and the transport part are treated in a similar way as in section \ref{gPC_Dir_1D}. For the oscillatory part, we use the backward Euler method in time,
\begin{equation}\frac{\vv W_j^{n+1}-\vv W_j^n}{\Delta s} = -\frac{1}{\varepsilon}\partial_{\tau}\vv W_j^{n+1}, \label{tau_part} \end{equation}
where $\vv W_j^n$ is an approximation of $\vv W(n\Delta s, x_j)$.
Let
$\displaystyle\vv{\hat W}_j^n(\zeta)$ be the Fourier transform of $\vv{W_j}^n$ in the periodic variable $\tau$, where $\zeta$ is the Fourier variable. The spectral method is used to discretize the $\tau$-derivative, then \eqref{tau_part} becomes
$$\vv{\hat W_j}^{n+1}=\left(I+i\zeta\frac{\Delta s}{\varepsilon}I\right)^{-1}\vv{\hat W_j}^n\,.$$
Note that the matrix $\displaystyle I+i\zeta\frac{\Delta s}{\varepsilon}I$ has non-zero eigenvalues, thus is invertible. \\[2pt]

{\bf The interpolation step}

To find the values of $V$ at time $T$, since $\displaystyle V(T, x_j, \tau_k, \bz^{(l)})=W(S(T, x_j, \bz^{(l)}), x_j, \tau_k, \bz^{(l)})$,
one uses linear interpolation to find
$W$ at $(S(T, x_j, \bz^{(l)}), x_j, \tau_k, \bz^{(l)})$, for $l=1,\cdots, N_s$.

We search for the interval $[s_{n_l}, s_{n_l+1}]$ that contains $S(T, x_j, \bz^{(l)})$. Denote $S(T, x_j, \bz^{(l)})=s^{\dagger}$, $s_{n_l}=s^{(1)}$, $s_{n_l+1}
=s^{(2)}$, $W(s_{n_l})=W^{(1)}$, $W(s_{n_l+1})=W^{(2)}$, then by the Lagrange interpolation formula,
$$
W(s^{\dagger})=W^{(1)}\frac{s^{\dagger}-s^{(2)}}{s^{(1)}-s^{(2)}}
+W^{(2)}\frac{s^{\dagger}-s^{(1)}}{s^{(2)}-s^{(1)}},
$$
which gives $V(T,x_j,\tau_k, \bz^{(l)})=W(s^{\dagger}, x_j, \tau_k, \bz^{(l)})$.

\begin{remark}
Compare to gPC-SG-N1, the transport step is more costly since one has to multiply by the matrix $A^{\star}$.  If $a=O(1)$, then the CFL condition for the transport equation requires just $\Delta t=O(\Delta x)$ thus independent of $\varepsilon$. However,
 for avoided crossing, typically $\displaystyle a(x,\bz)=\mathcal{O}(\sqrt{\varepsilon})$ for some region of $x$,
thus the eigenvalue of $A^{\star}$ are of $\displaystyle\mathcal{O}(1/\sqrt{\varepsilon})$, thus an explicit time discretization of the transport step
will require $\displaystyle\Delta t=\mathcal{O}(\sqrt{\varepsilon}\Delta x)$, further increasing the computational cost.
One may use implicit time discretization in the transport step. This will be studied in our future research.
\label{rmk}
\end{remark}

\subsection{Numerical tests for the 1D scalar equation}
For numerical comparison, we also use the stochastic collocation (SC) method
\cite{GWZ, XH}.
Let $\{\bz^{(j)}\}_{j=1}^{N_c}\subset I_{\bz}$ be the set of collocation nodes, $N_c$ the number of samples.
For each fixed individual sample $\bz^{(j)}$, $j=1,\ldots,N_c$, one applies the deterministic solver to the deterministic equations as in \cite{NGO}, obtains the solution ensemble $u_j(t,x)=u(t,x,\bz^{(j)})$, then adopts an interpolation approach to construct a gPC approximation, such as $$u(t,x,\bz)=\sum_{j=1}^{N_c}u_j(t,x)l_j(\bz),$$ where $l_j(\bz)$ depends on the construction method. The Lagrange interpolation method is used here by choosing $l_j(\bz^{(i)})=\delta_{ij}$. Depending on whether the direct or NGO-based method is used for the
deterministic method, we will have the gPC-SC-D and gPC-SG-N1, gPC-SG-N2 methods.
\\[20pt]
{\bf Example $2.1$}

\noindent We consider the numerical example in \cite{NGO}, with
$a$ involving a 1D random variable $z$ that follows a uniform distribution,
\begin{align*}
r(u)=u^2/(u^2+2|u|^2), \quad c(x)=\cos^2(x), \quad a(x,z)=(3/2+\cos(2x))(1+0.5z)>0,
\end{align*}
and the non-oscillatory initial data given by
$$u(0,x)=1+\frac{1}{2}\cos(2x)+i\left[1+\frac{1}{2}\sin(2x)\right], \quad x\in [-\pi/2,\pi/2].$$

In the following tests, for all the reference solutions obtained from gPC-SC-D method, $N_c=64$, $N_g=32$ quadrature points are used.
In gPC-SG-N1 and gPC-SG-N2 methods, $N_s$ stands for the number of quadrature points used in the final reconstruction step
(to get $u$ from $V$ or $W$).

\begin{figure} [H]
\begin{subfigure}{0.5\textwidth}
 \includegraphics[trim={0 8cm 0 0},width=1.1\textwidth, height=0.8\textwidth]{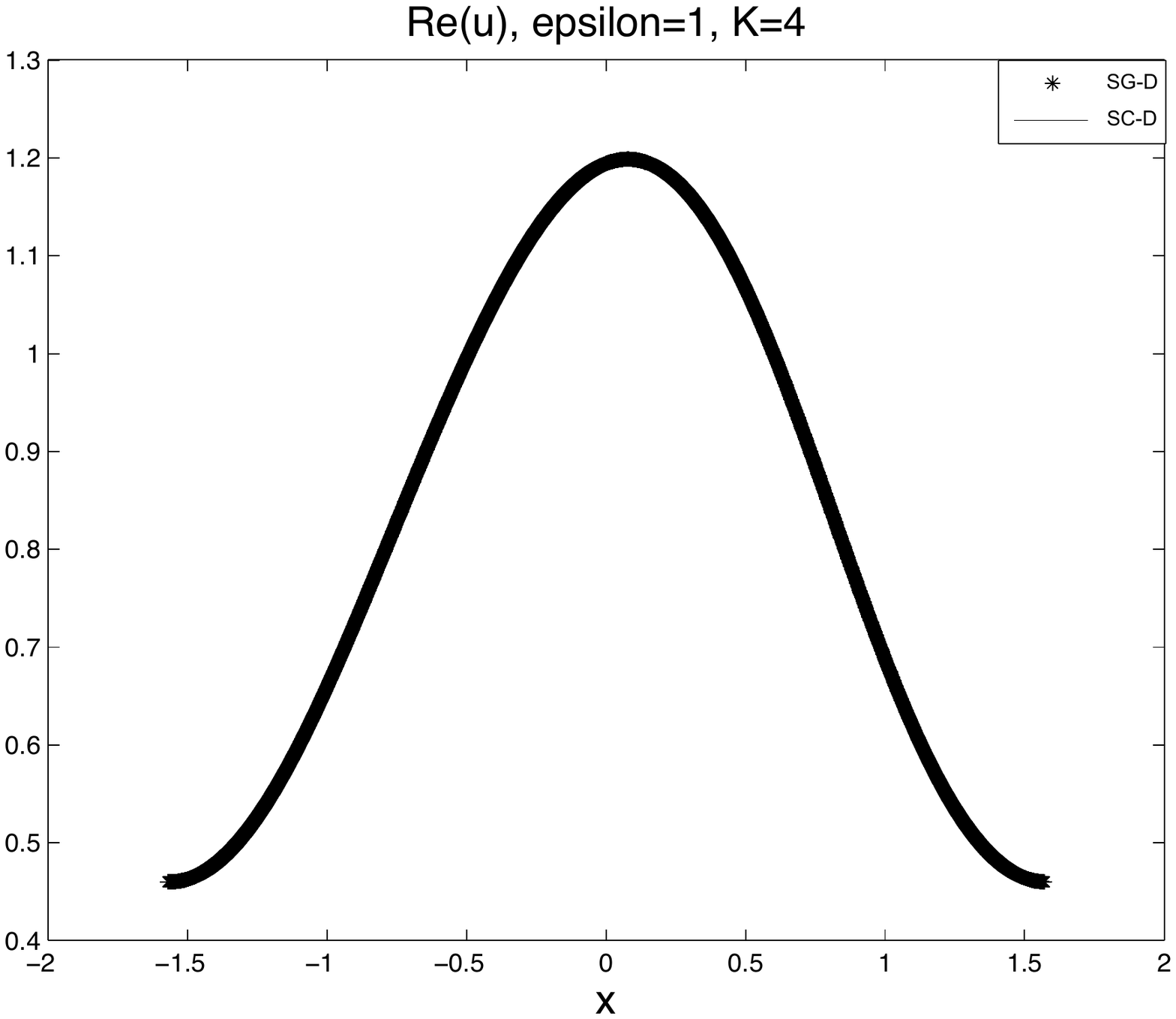}
  \end{subfigure}
  \begin{subfigure}{0.5\textwidth}
 \includegraphics[trim={0 8cm 0 0},width=1.1\textwidth, height=0.8\textwidth]{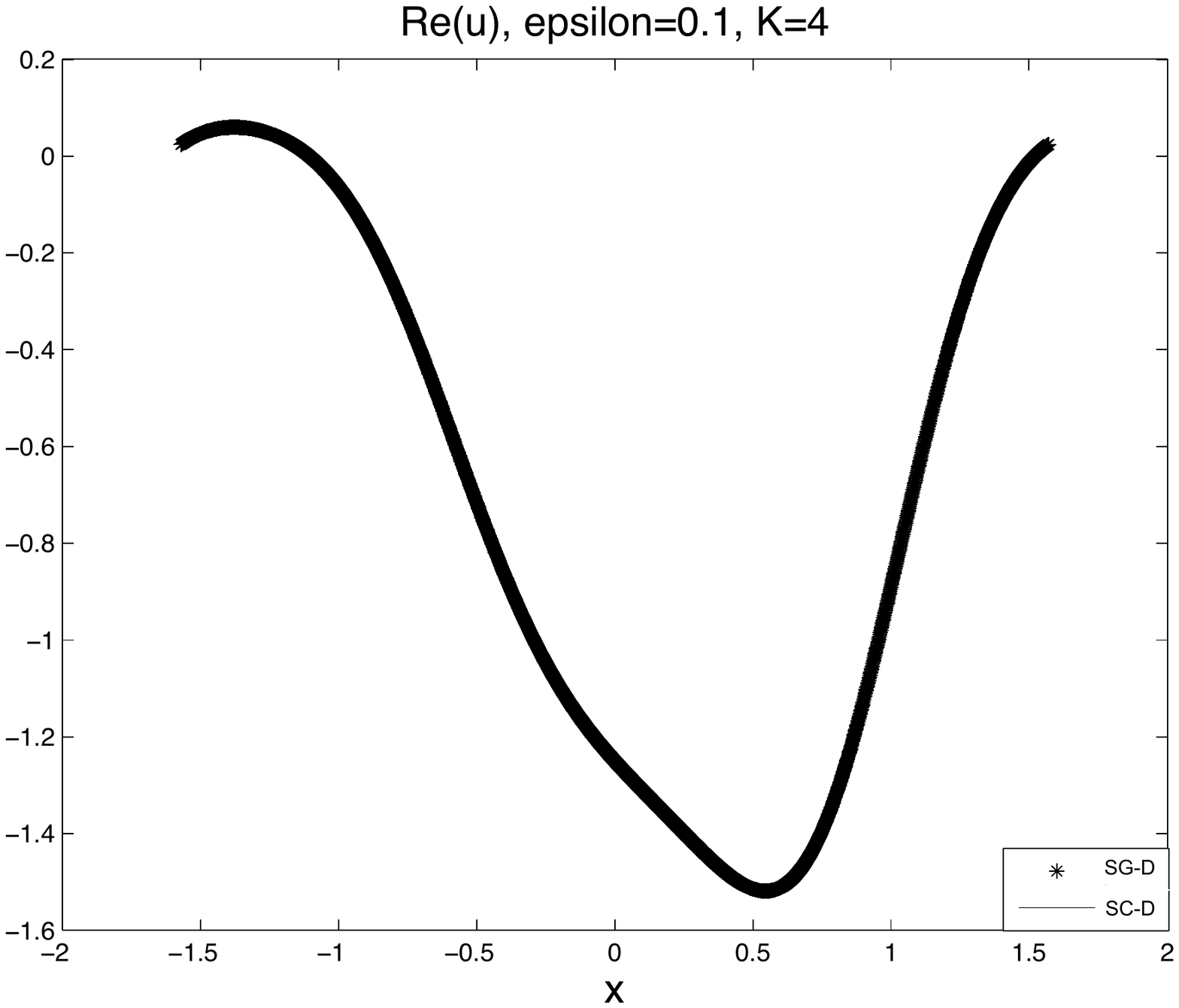}
  \end{subfigure}
 \begin{subfigure}{0.5\textwidth}
   \centering
 \includegraphics[trim={0 8cm 0 0},width=1.1\textwidth, height=0.8\textwidth]{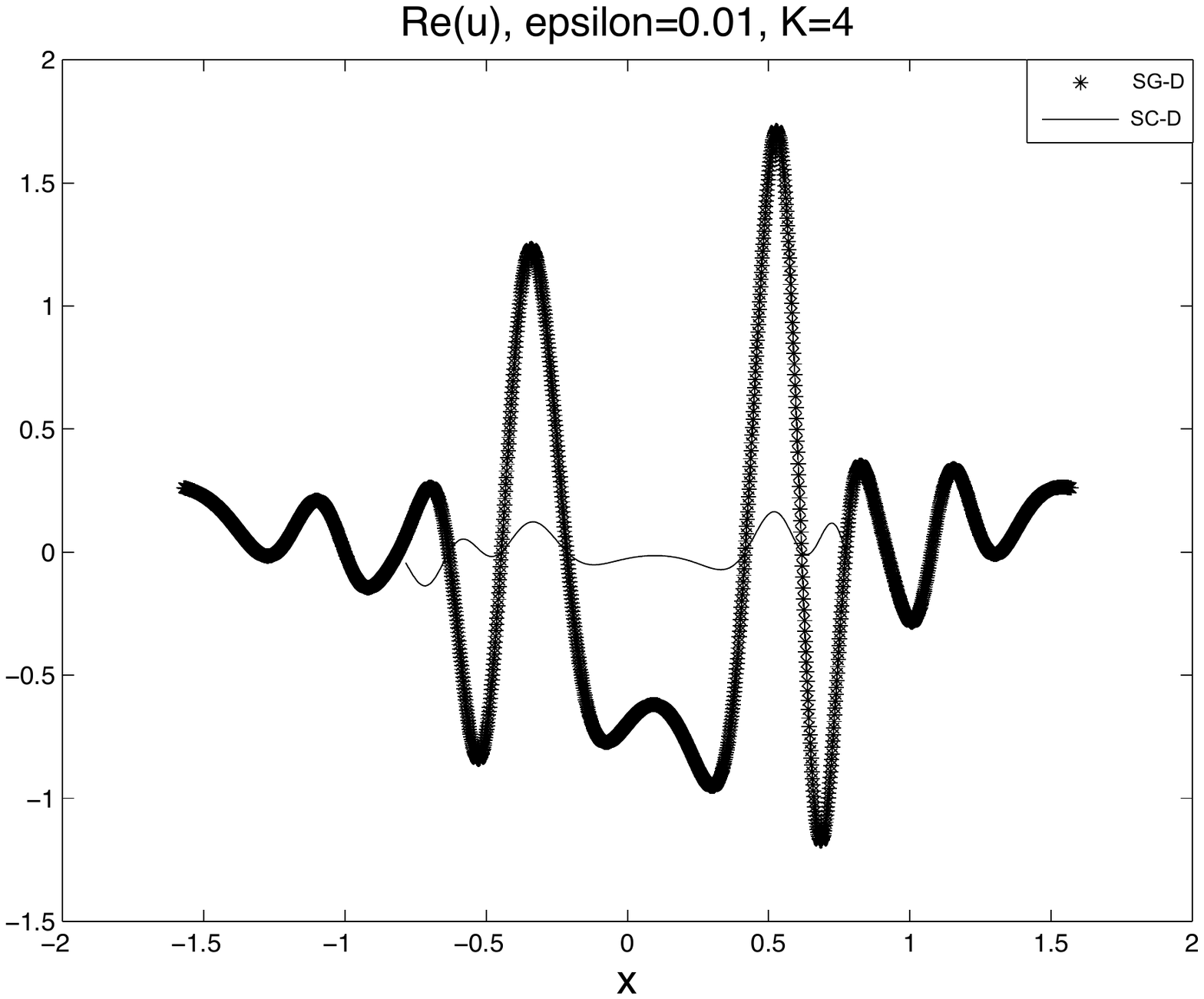}
  \end{subfigure}
  \begin{subfigure}{0.5\textwidth}
   \centering
 \includegraphics[trim={0 8cm 0 0},width=1.1\textwidth, height=0.8\textwidth]{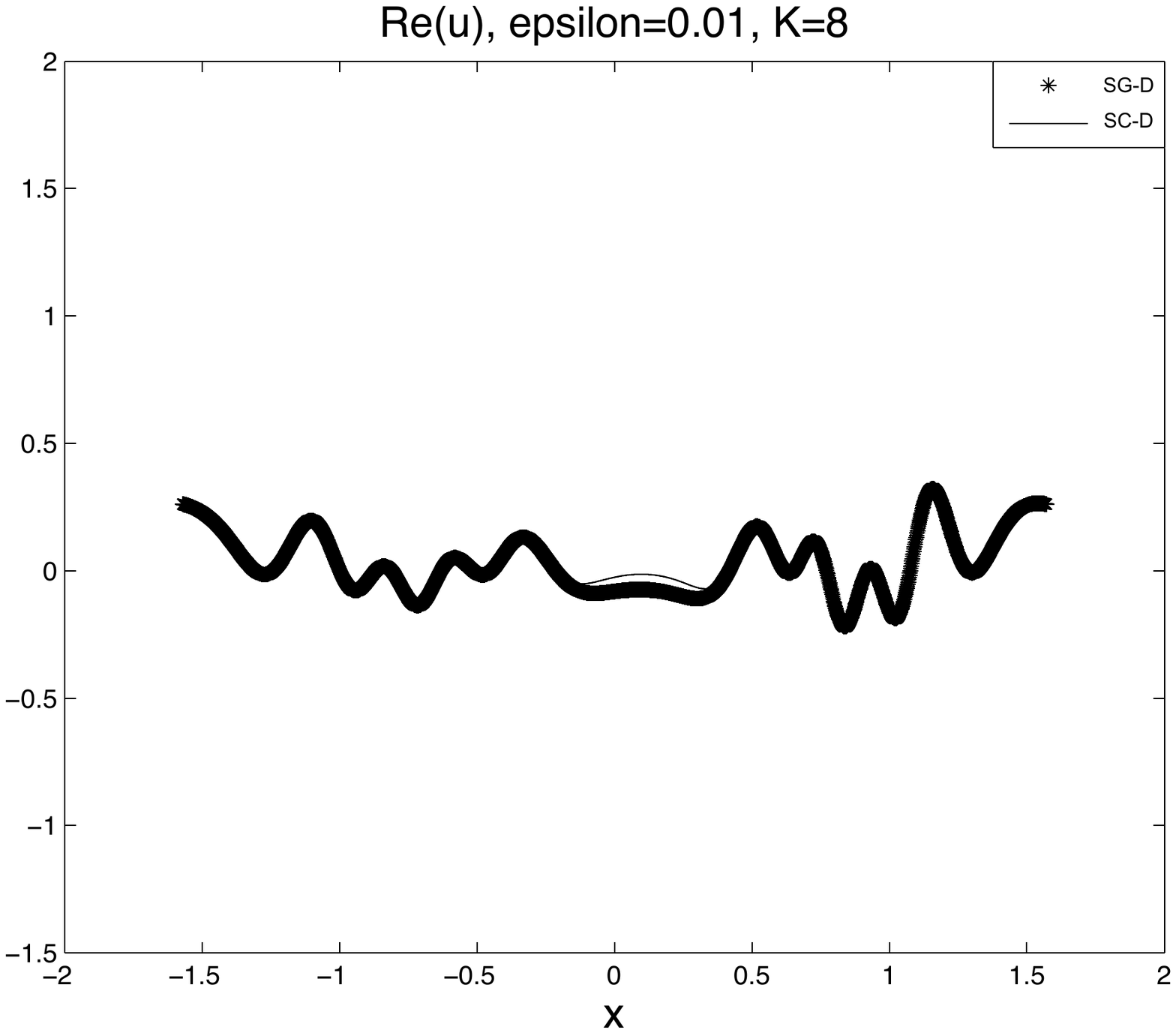}
  \end{subfigure}
  \begin{subfigure}{0.5\textwidth}
   \centering
 \includegraphics[trim={0 8cm 0 0},width=1.1\textwidth, height=0.8\textwidth]{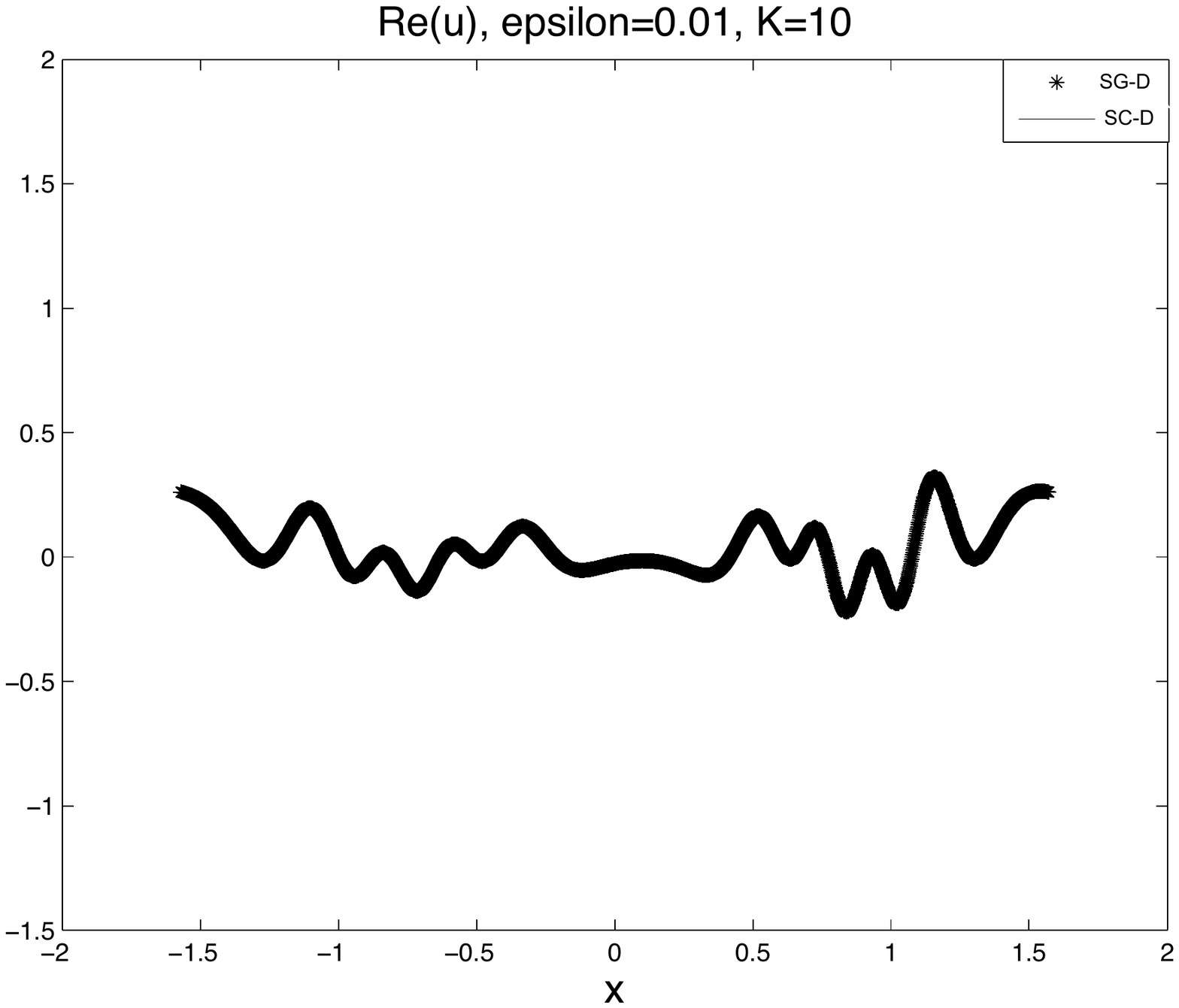}
  \end{subfigure}
  \begin{subfigure}{0.5\textwidth}
   \centering
 \includegraphics[trim={0 8cm 0 0},width=1.1\textwidth, height=0.8\textwidth]{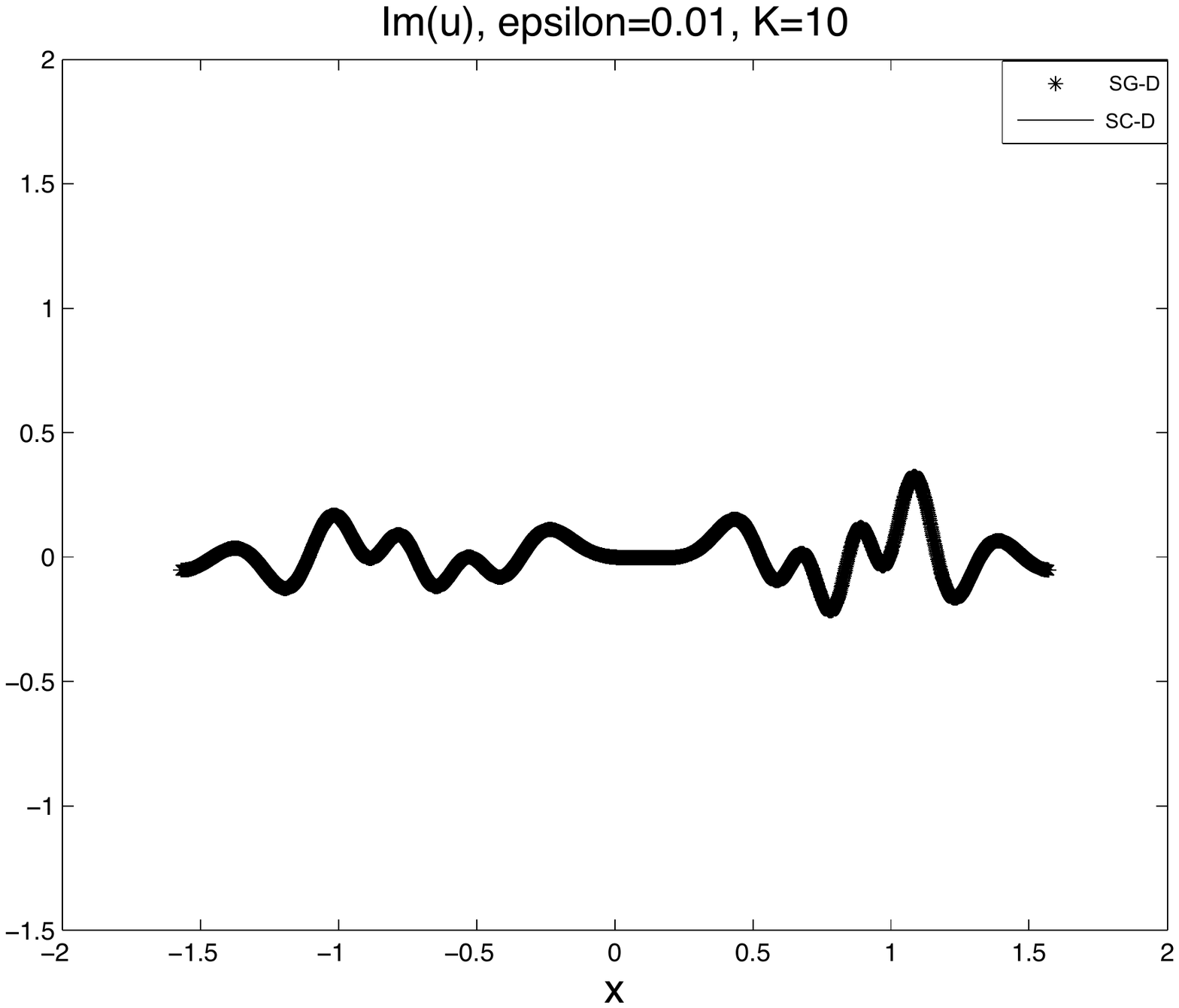}
  \end{subfigure}
  \caption{Example $2.1$. Mean of real parts of $u$ at $t=0.1$, $\varepsilon=1$, $0.1$ and $0.01$.
  $N_x=2000$, $\Delta t=10^{-5}$. Stars: gPC-SG-D with $K=4$.
  Solid lines: reference solution by gPC-SC-D. }
  \label{Dir1}
  \end{figure}

\begin{figure} [H]
\begin{subfigure}{0.55\textwidth}
   \centering
 \includegraphics[width=1.1\textwidth, height=1.0\textwidth]{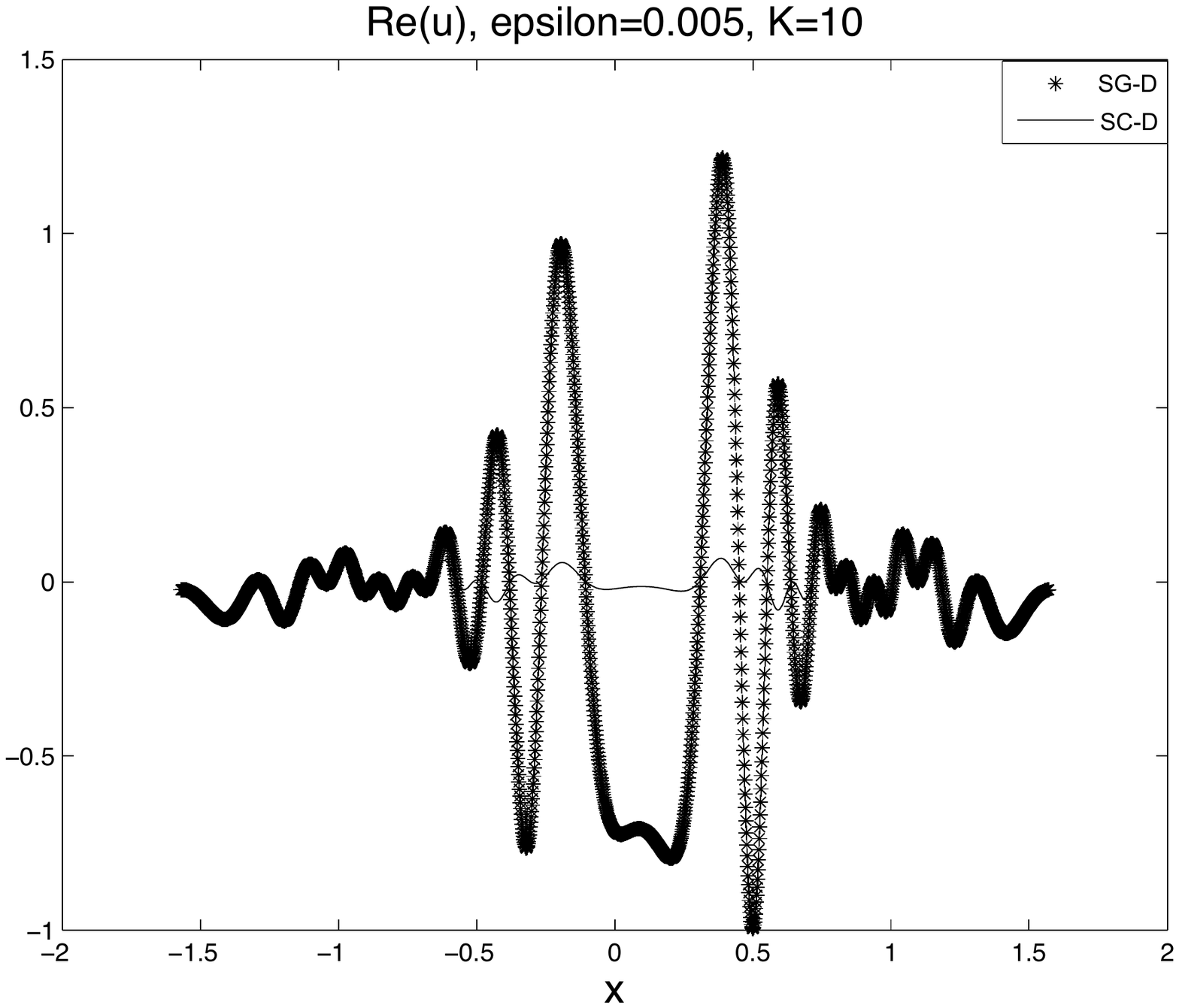}
  \end{subfigure}
    \begin{subfigure}{0.55\textwidth}
   \centering
 \includegraphics[width=1.1\textwidth, height=1.0\textwidth]{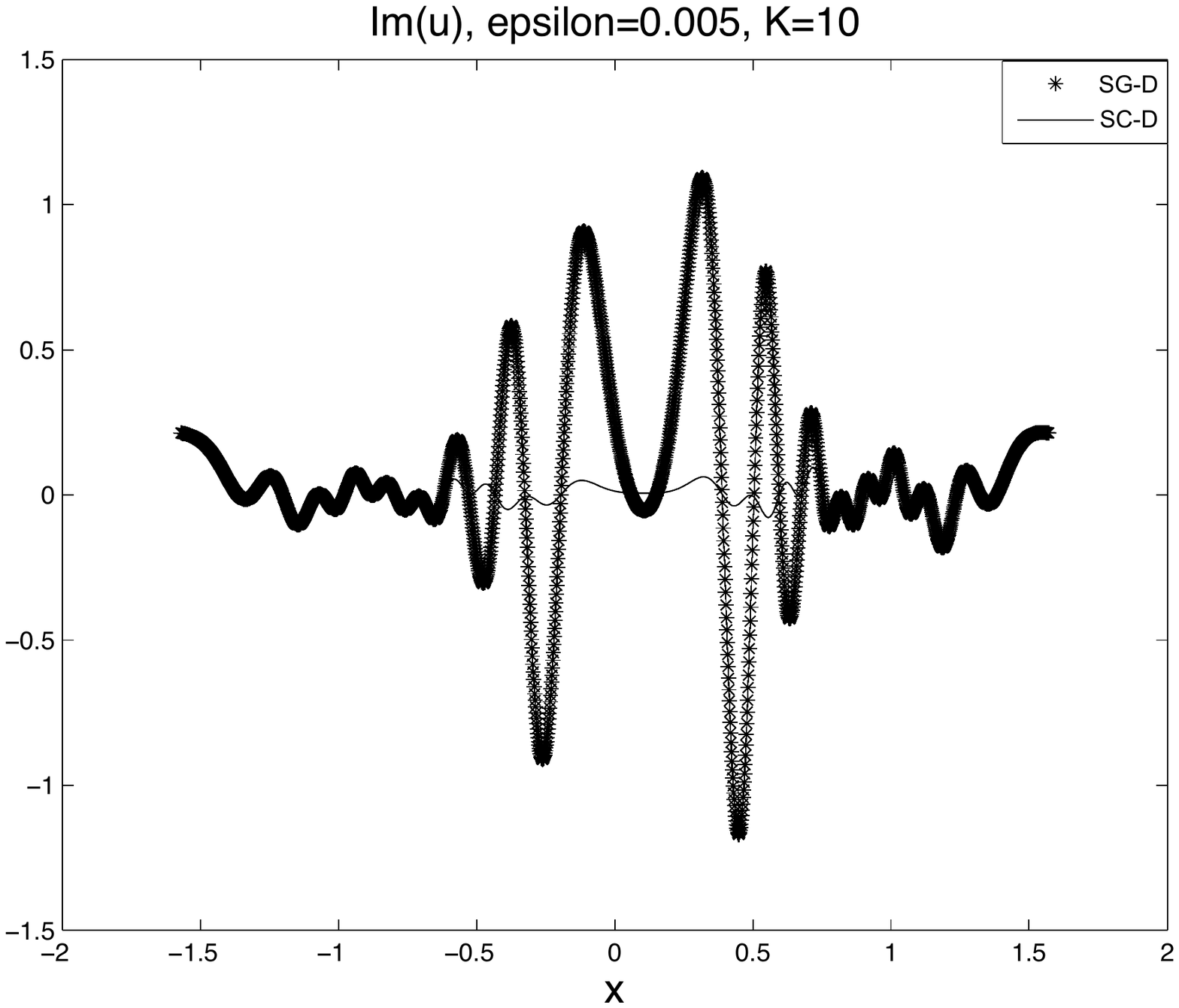}
  \end{subfigure}
  \begin{subfigure}{0.55\textwidth}
   \centering
 \includegraphics[width=1.1\textwidth, height=1.0\textwidth]{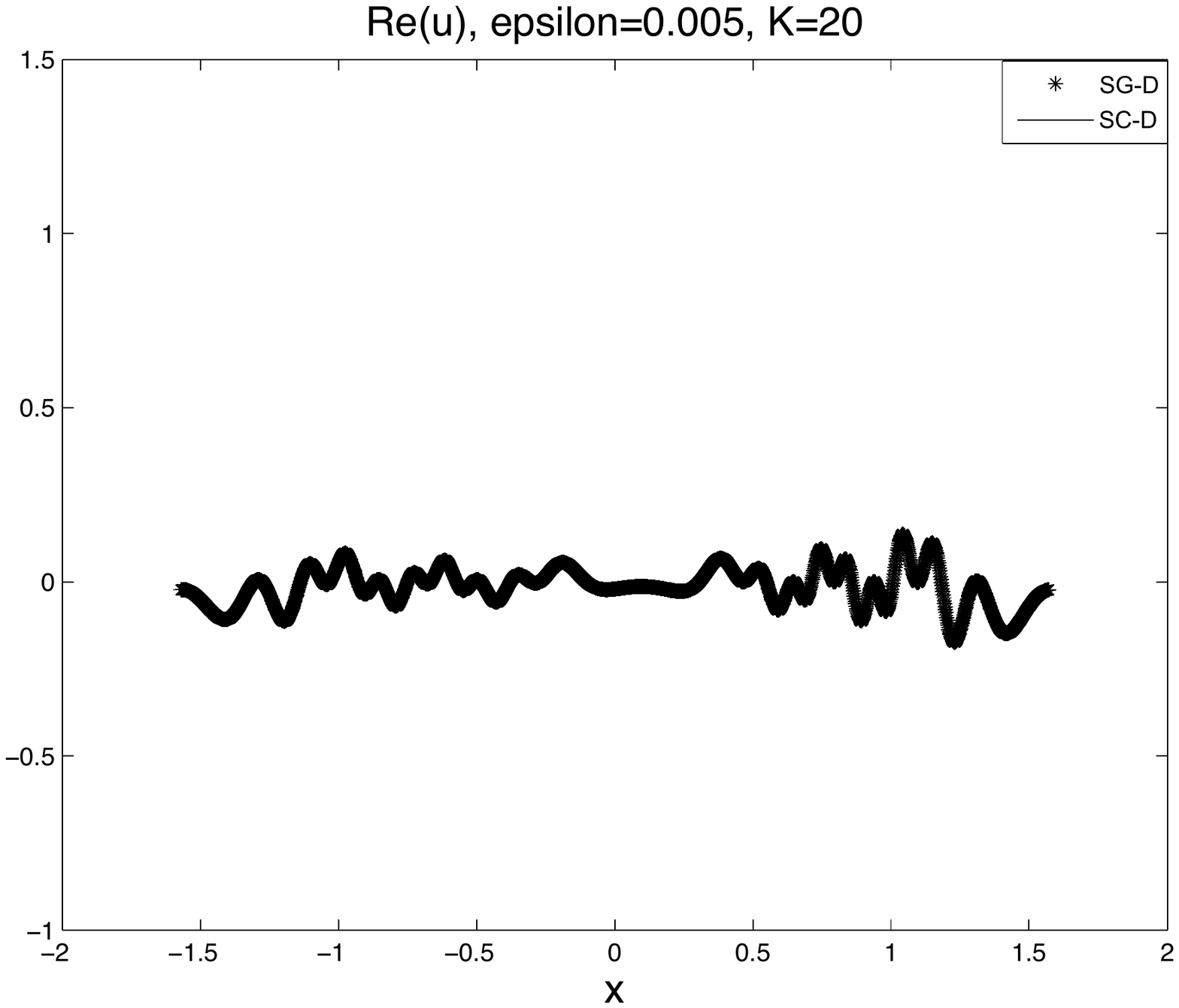}
  \end{subfigure}
  \begin{subfigure}{0.55\textwidth}
   \centering
 \includegraphics[width=1.1\textwidth, height=1.0\textwidth]{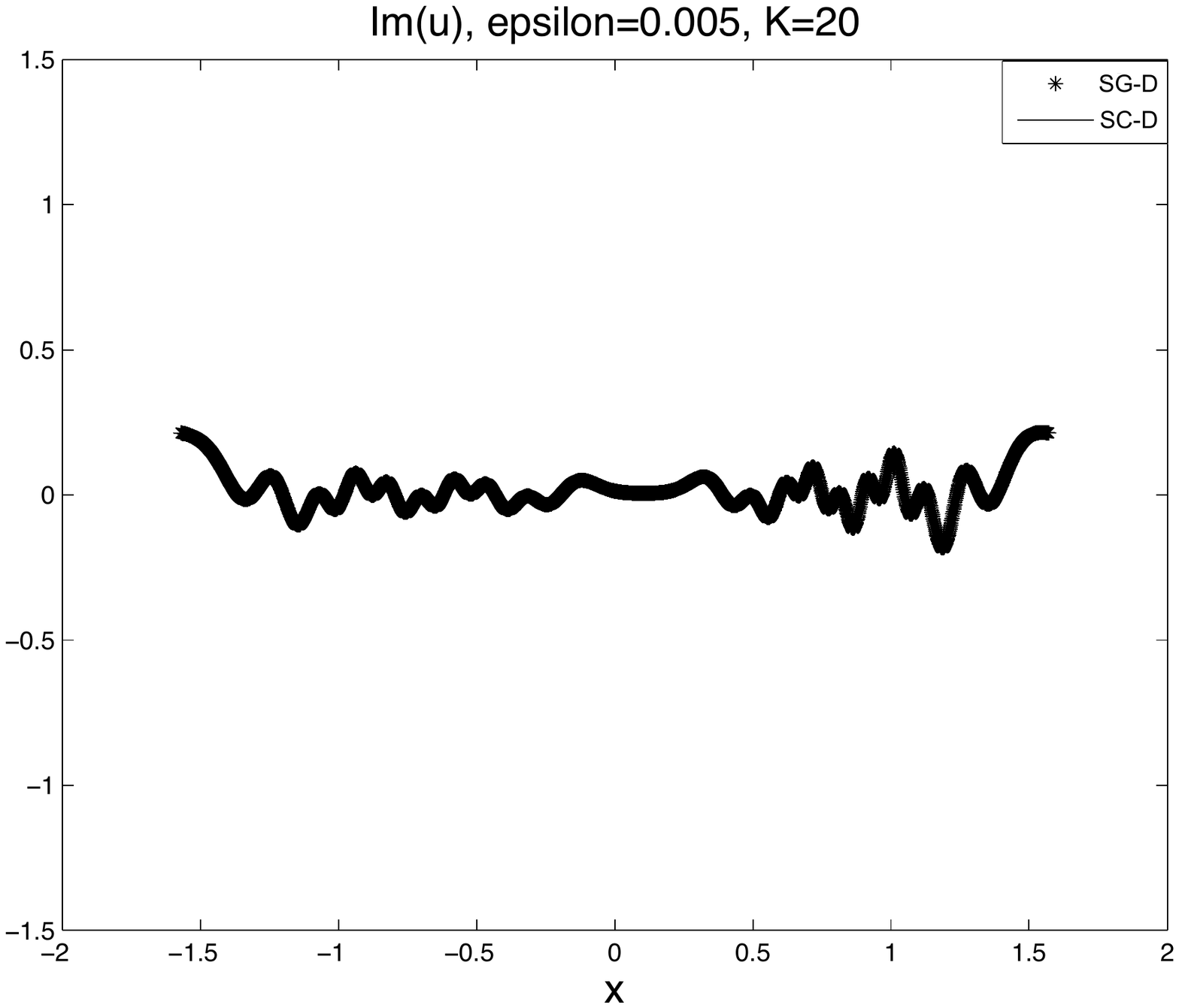}
  \end{subfigure}
  \caption{Example $2.1$. Mean of real parts of $u$ at $t=0.1$, $\varepsilon=5\times 10^{-3}$.
  $N_x=2000$, $\Delta t=10^{-5}$. Stars:  gPC-SG-D. Solid lines: reference solution by gPC-SC-D. }
\label{Dir2}
  \end{figure}

Figures \ref{Dir1} and \ref{Dir2} demonstrate how
$\varepsilon$ affects the choices of $K$ for gPC-SG-D.
Here the gPC-SC-D, computed with very small $\Delta t$ and $\Delta x$, is
used to obtain the reference solutions. For gPC-SG-D, we also use
very small $\Delta x$ and $\Delta t$, in order to concentrate on effect of
$\varepsilon$ on $K$.
We output the solution at $t=0.1$. In Figure \ref{Dir1}, for $\varepsilon=1$ and $0.1$, one can use
$K=4$. However, for $\varepsilon=0.01$, neither $K=4$ nor $K=8$
gives the correct solutions, especially around the center of the domain.
One starts to see satisfactory numerical solutions when $K=10$.
For $\varepsilon=0.005$ one needs to use $K=20$, see Figure \ref{Dir2}.

\begin{figure} [H]
 \begin{subfigure}{0.5\textwidth}
 \includegraphics[width=1.1\textwidth, height=1.1\textwidth]{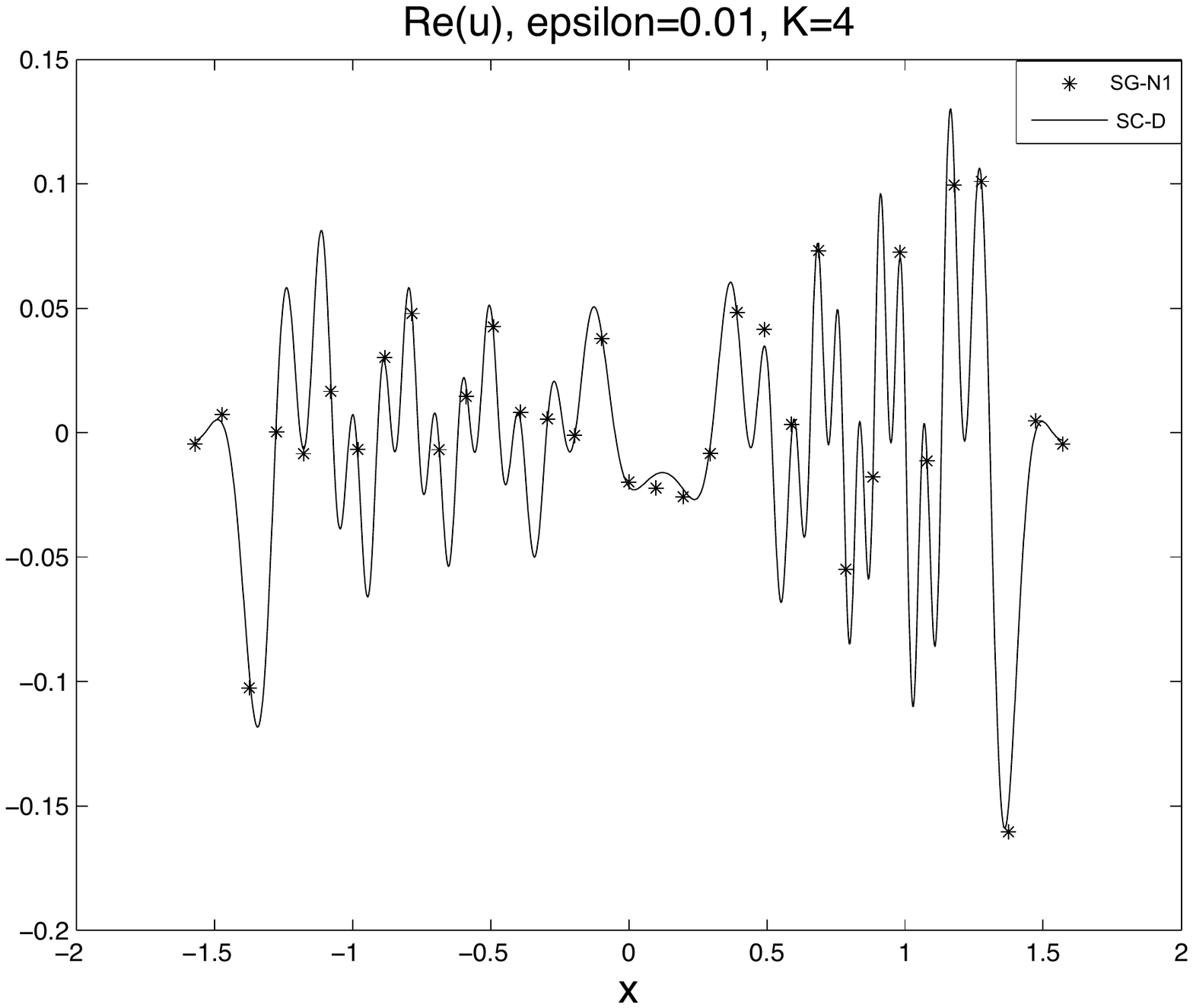}
 \end{subfigure}
 \begin{subfigure}{0.5\textwidth}
 \includegraphics[width=1.1\textwidth, height=1.1\textwidth]{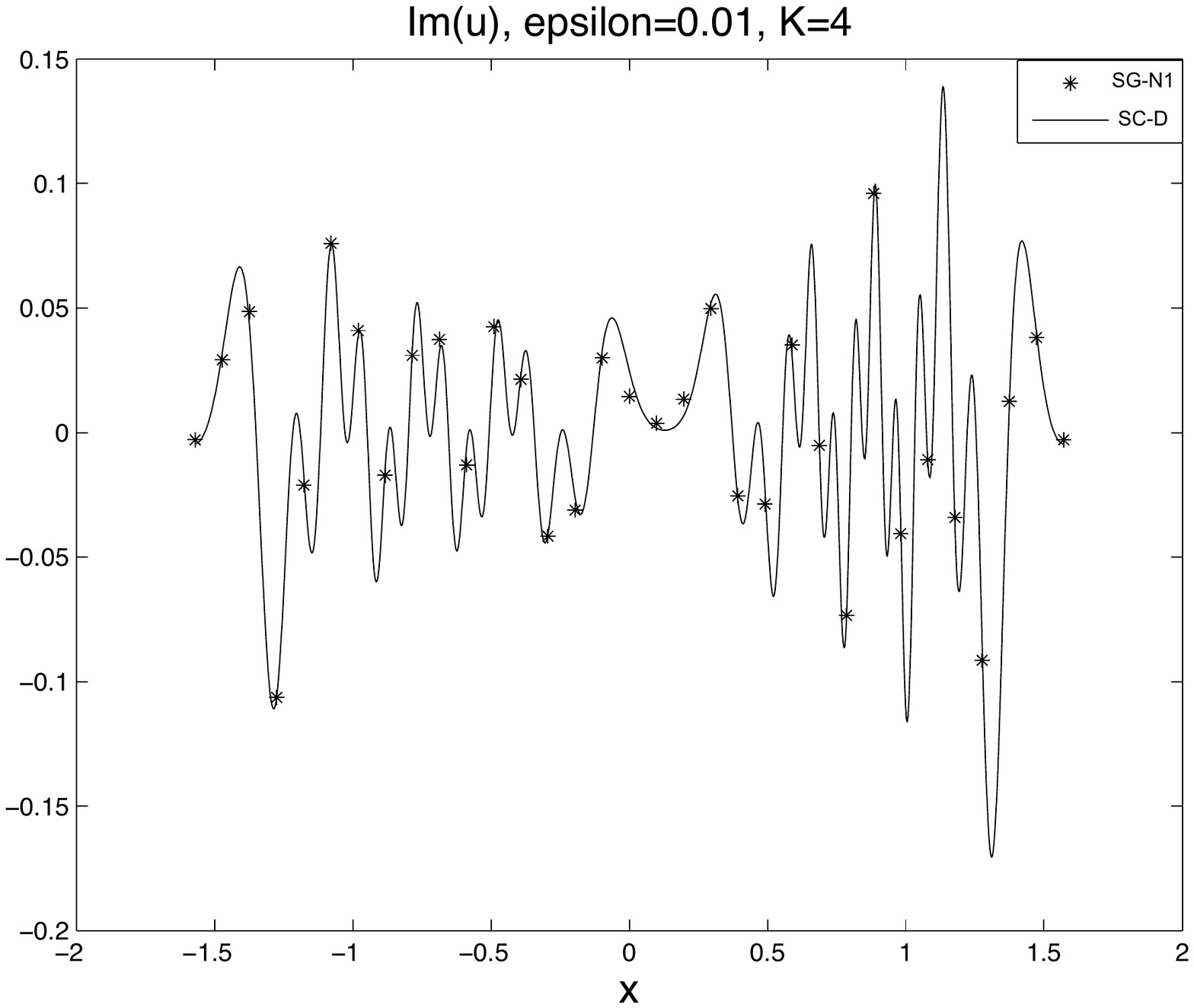}
  \end{subfigure}
 \begin{subfigure}{0.5\textwidth}
  \includegraphics[width=1.1\textwidth, height=1.1\textwidth]{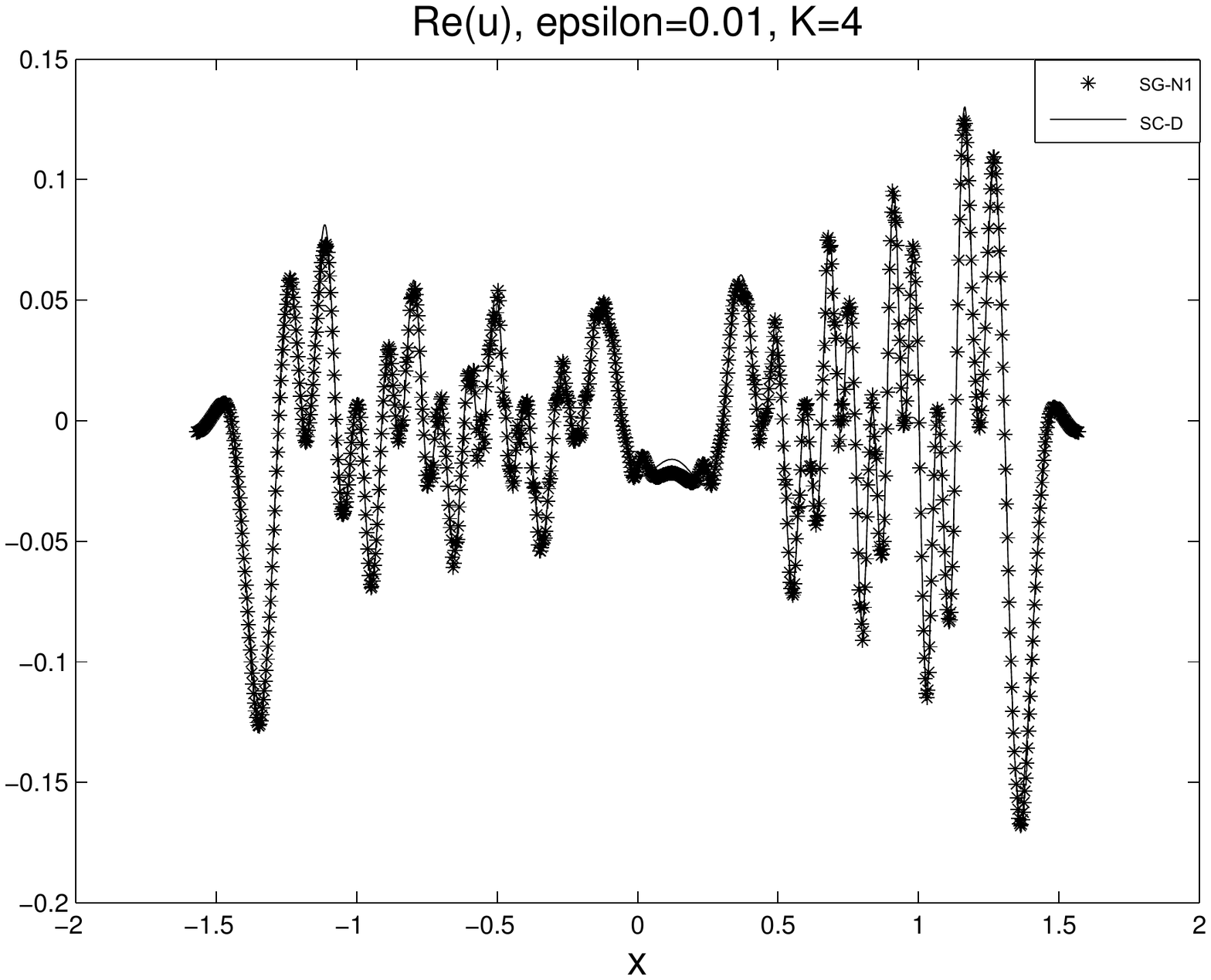}
 \end{subfigure}
  \begin{subfigure}{0.5\textwidth}
 \includegraphics[width=1.1\textwidth, height=1.1\textwidth]{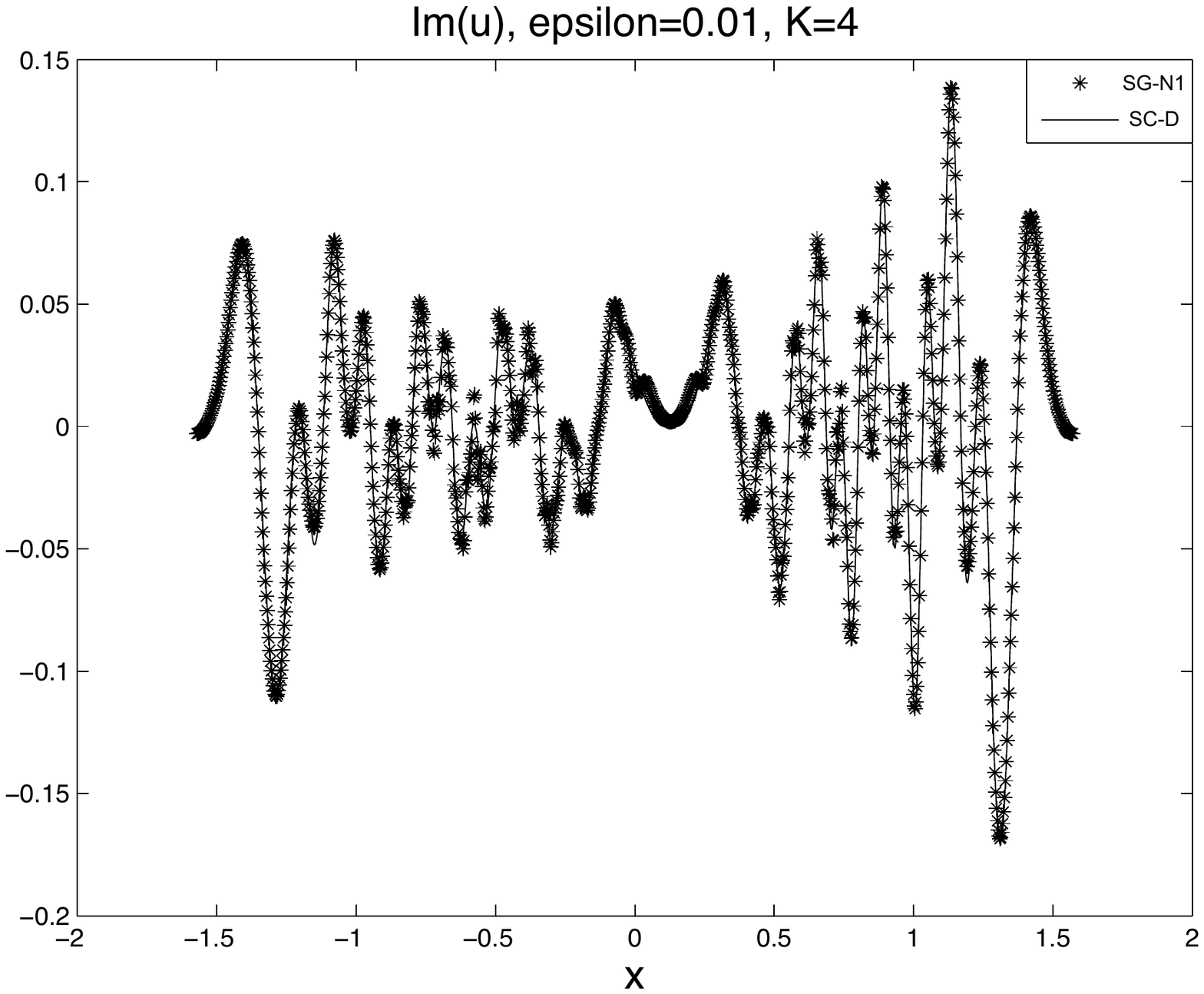}
 \end{subfigure}
  \caption{Example $2.1$. Mean of real and imaginary parts of $u$ at $t=0.25$, $\varepsilon=0.01$.
$\Delta x=\pi/32$, $\Delta t=0.01$ for the first row (gPC-SG-N1), and $\Delta x=\pi/1000$, $\Delta t=10^{-4}$ for the second row (gPC-SG-N1).
$\Delta x=\pi/1000$, $\Delta t=5\times 10^{-5}$ (gPC-SC-D).
Stars: gPC-SG-N1 with $K=4$. Solid lines: reference solution by gPC-SC-D. }
\label{New_Dir1_1D}
 \end{figure}

\begin{figure} [H]
 \begin{subfigure}{0.5\textwidth}
   \centering
 \includegraphics[trim={0 8cm 0 0},width=1.0\textwidth, height=0.8\textwidth]{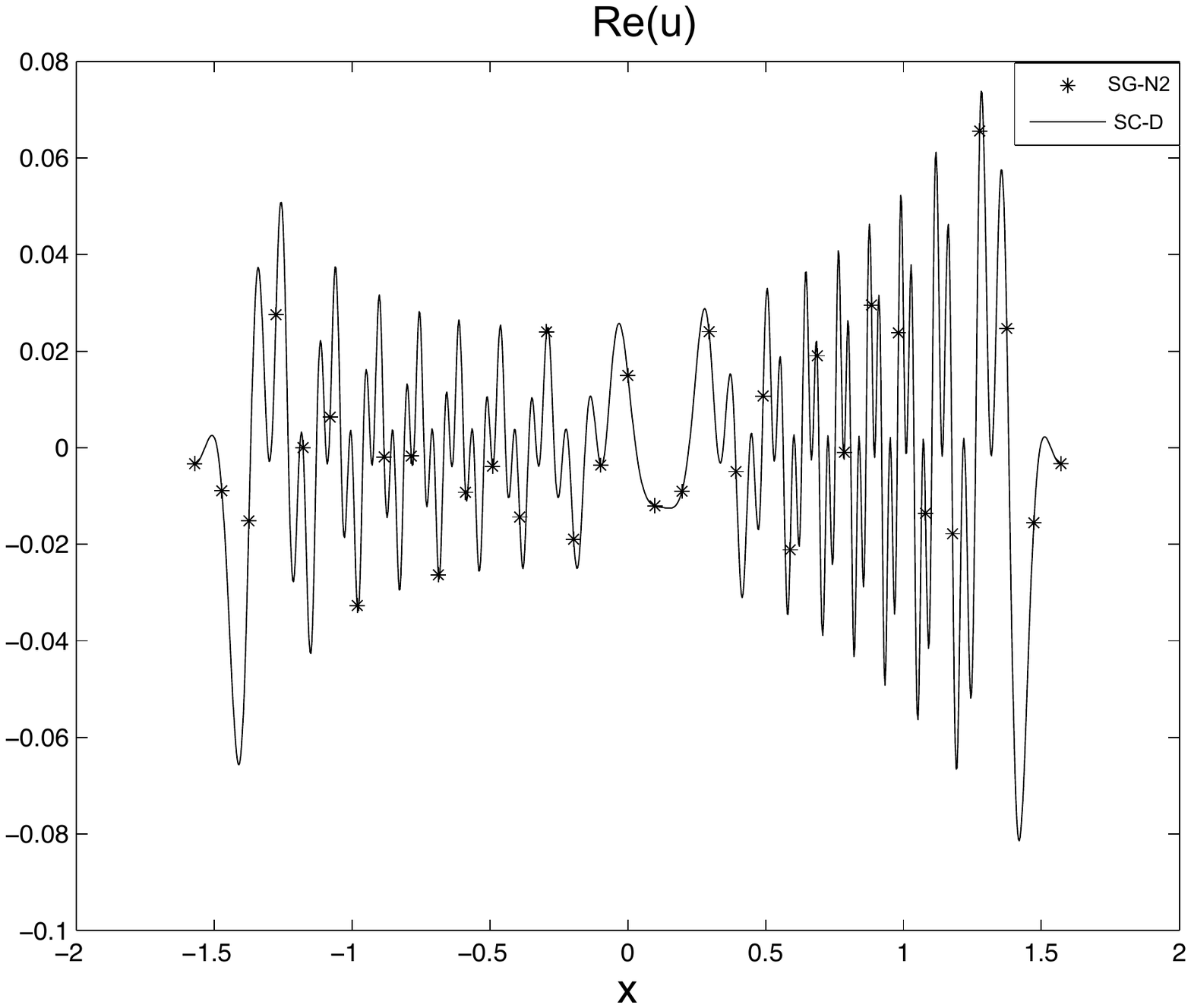}
   \end{subfigure}
  \begin{subfigure}{0.5\textwidth}
 \centering
 \includegraphics[trim={0 8cm 0 0},width=1.0\textwidth, height=0.8\textwidth]{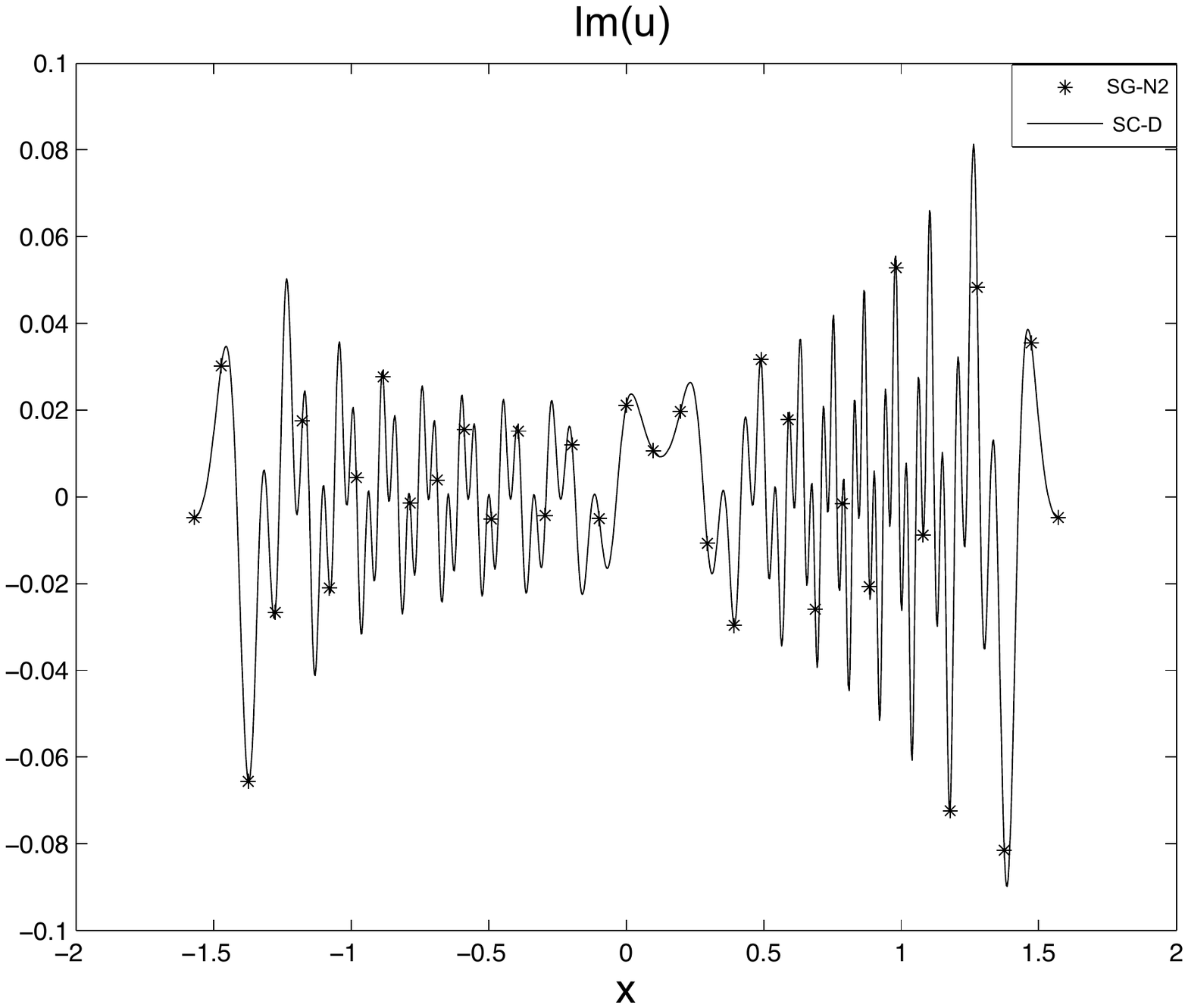}
   \end{subfigure}
  \begin{subfigure}{0.5\textwidth}
   \centering
 \includegraphics[trim={0 8cm 0 0},width=1.0\textwidth, height=0.8\textwidth]{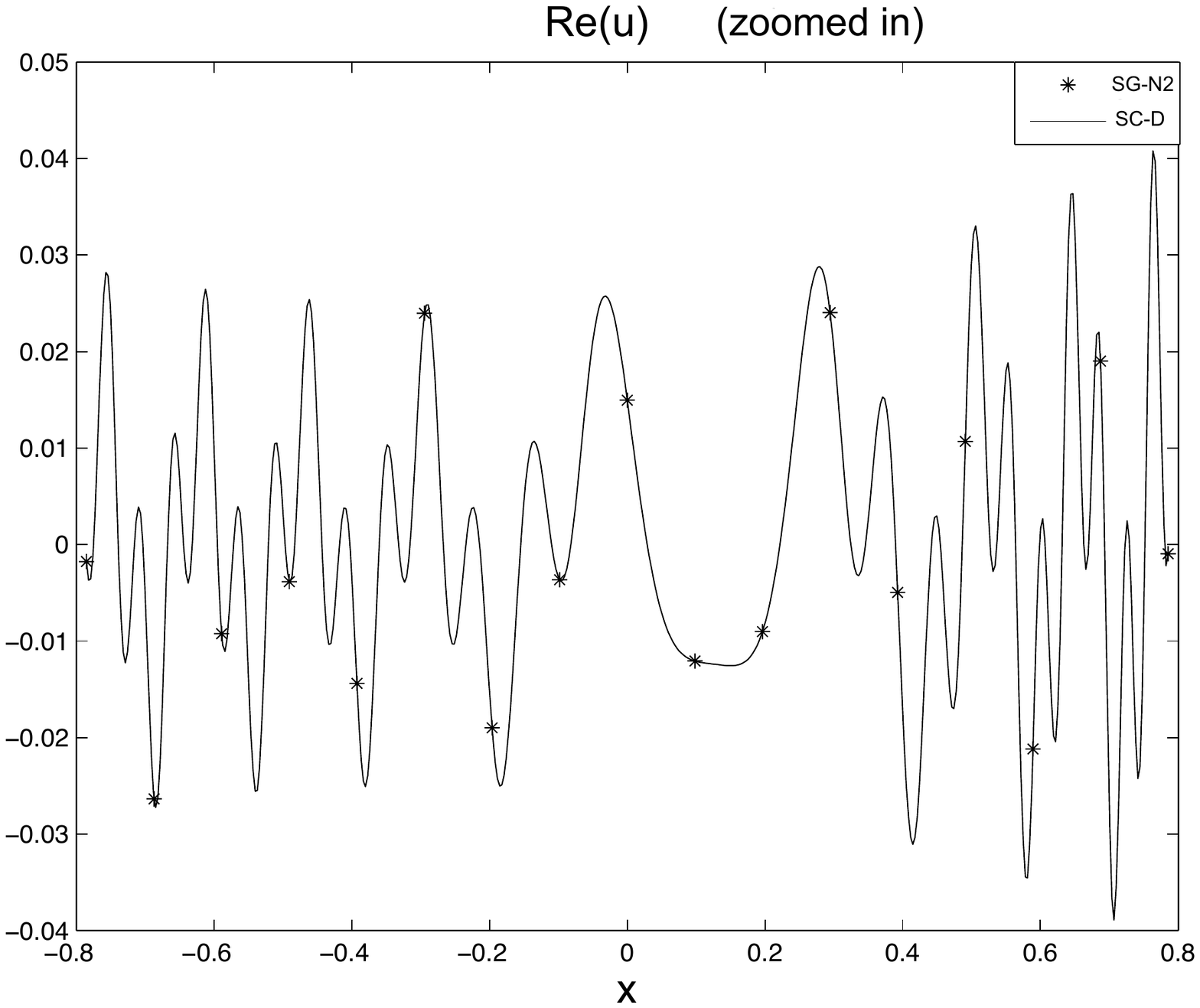}
  \end{subfigure}
  \begin{subfigure}{0.5\textwidth}
 \centering
 \includegraphics[trim={0 8cm 0 0},width=1.0\textwidth, height=0.8\textwidth]{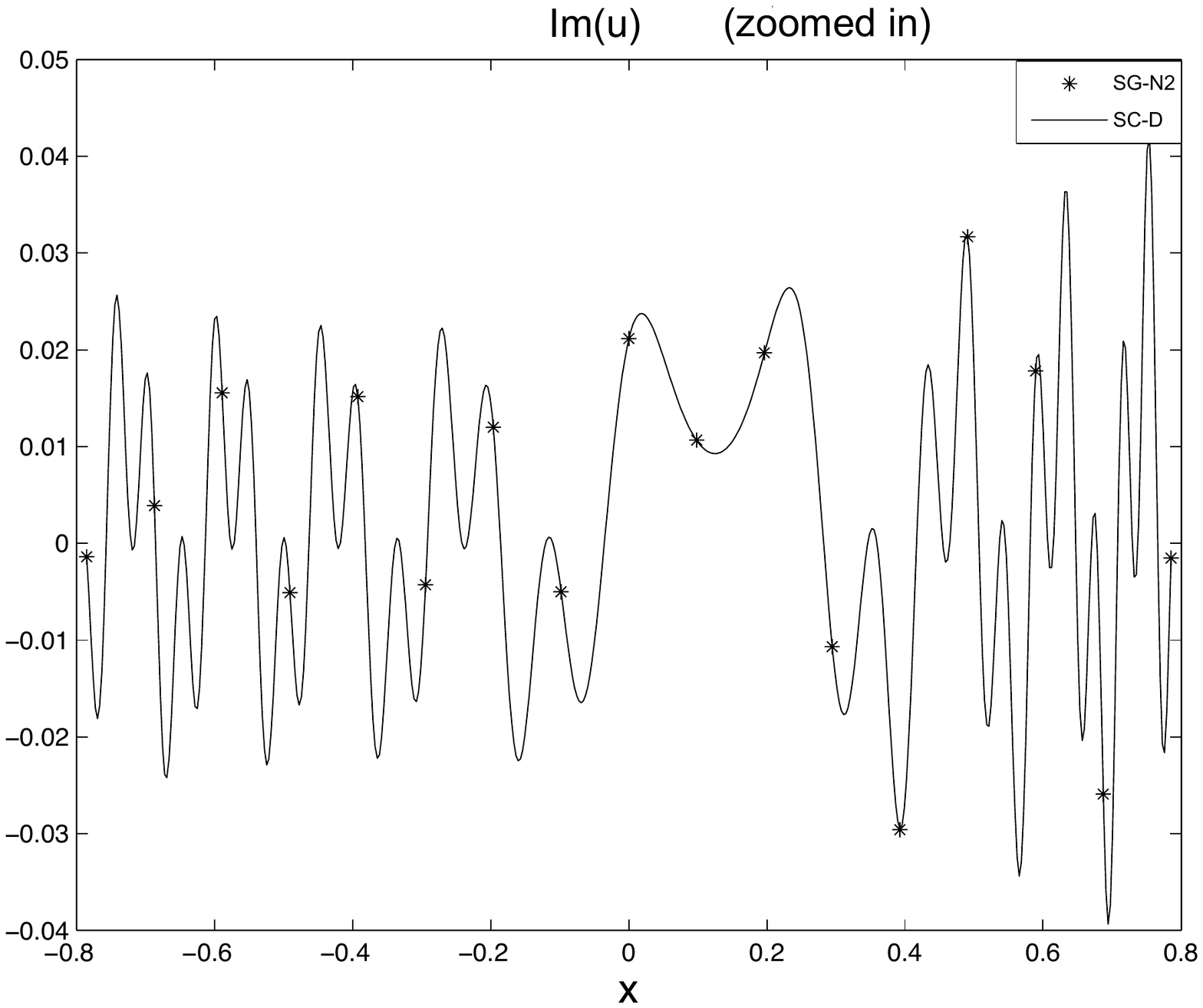}
   \end{subfigure}
   \begin{subfigure}{0.5\textwidth}
   \centering
 \includegraphics[trim={0 8cm 0 0},width=1.0\textwidth, height=0.8\textwidth]{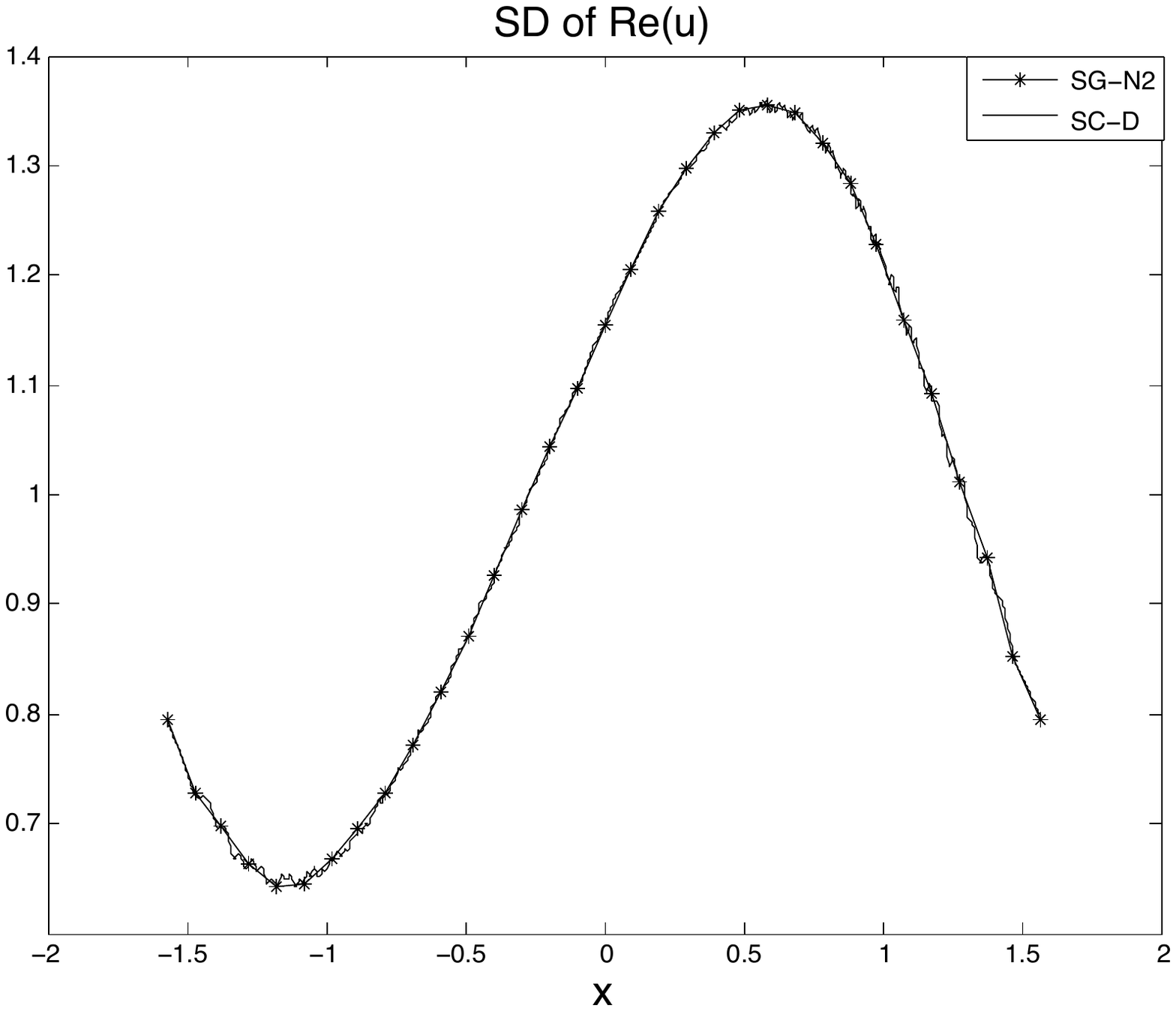}         
   \end{subfigure}
  \begin{subfigure}{0.5\textwidth}
 \centering
 \includegraphics[trim={0 8cm 0 0},width=1.0\textwidth, height=0.8\textwidth]{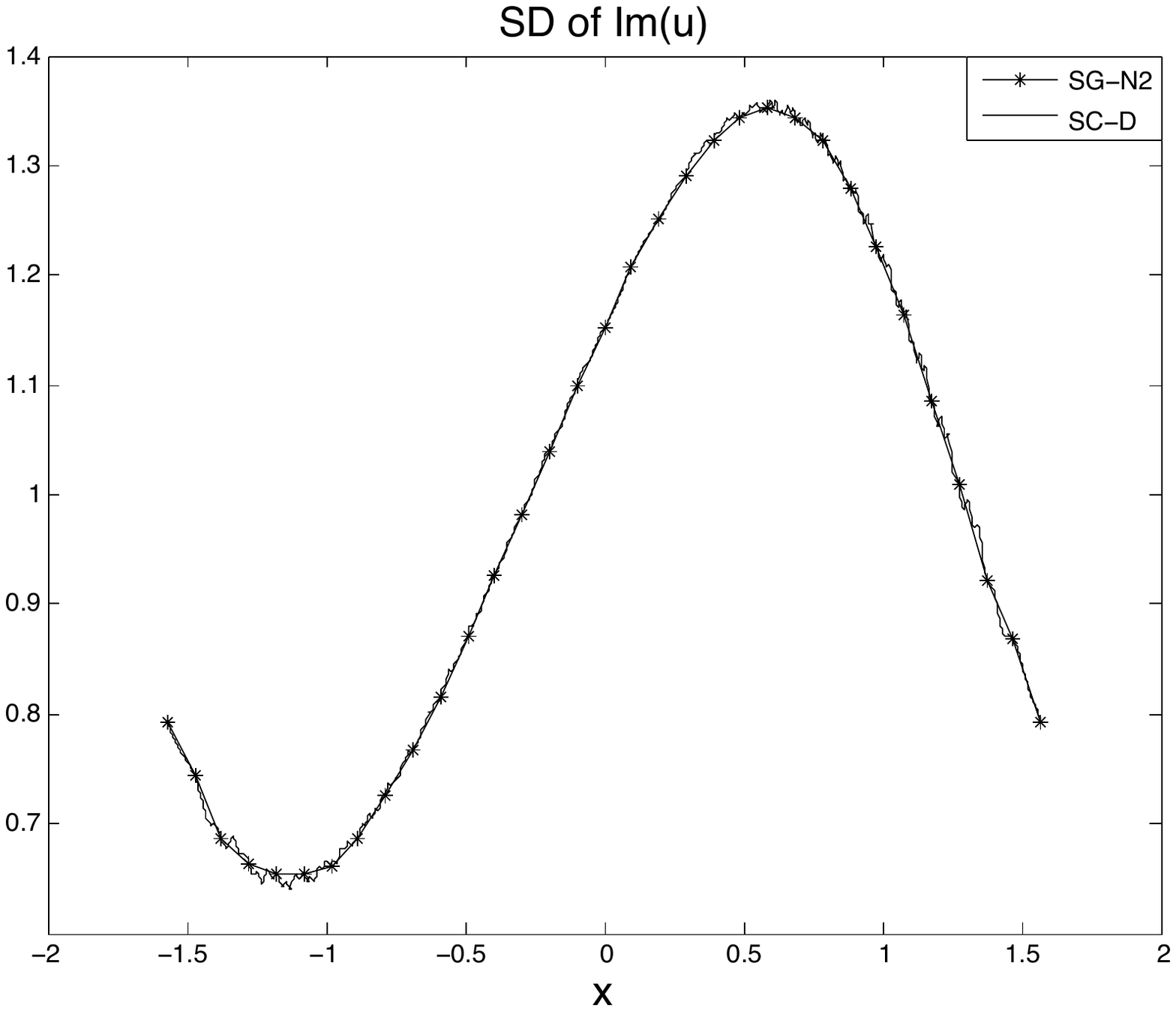}           
   \end{subfigure}
 \caption{Example $2.1$. Mean and standard deviation of real and imaginary parts of $u$ and their zoomed in solutions at $t=0.25$, $\varepsilon=5\times 10^{-3}$.
$\Delta x=\pi/32$, $\Delta t=0.01$, $K=4$ (gPC-SG-N2), and $\Delta x=\pi/1000$, $\Delta t=5\times10^{-5}$ (gPC-SC-D).
Stars: gPC-SG-N2. Solid lines: reference solution by gPC-SC-D.}
\label{New_Dir2_1D}
 \end{figure}

Figures \ref{New_Dir1_1D} and \ref{New_Dir2_1D} compare mean of the solutions of gPC-SG-N1 and gPC-SG-N2 for $\varepsilon=0.01$, and
$\varepsilon=0.005$ respectively. We take small $\varepsilon$ and use gPC-SC-D to fully resolve the oscillations (to serve as a reference solution).
One can observe from Figure \ref{New_Dir2_1D} that gPC-SG-N2 is able to capture correctly the solution (its mean and standard deviation)
at mesh points with $\Delta x, \Delta t$ much larger than
$\varepsilon$ and $K=4$. For the same grids, gPC-SG-N1 is a bit off near $x=0$.
If one uses a much finer spatial size $\Delta x=\pi/1000$, gPC-SG-N1 becomes more accurate but still is a bit off near $x=0$,
as shown in the second row of Figure \ref{New_Dir1_1D}.

\begin{figure} [H]
 \begin{subfigure}{0.5\textwidth}
 \includegraphics[trim={0 10cm 0 0cm},width=1.0\textwidth, height=0.7\textwidth]{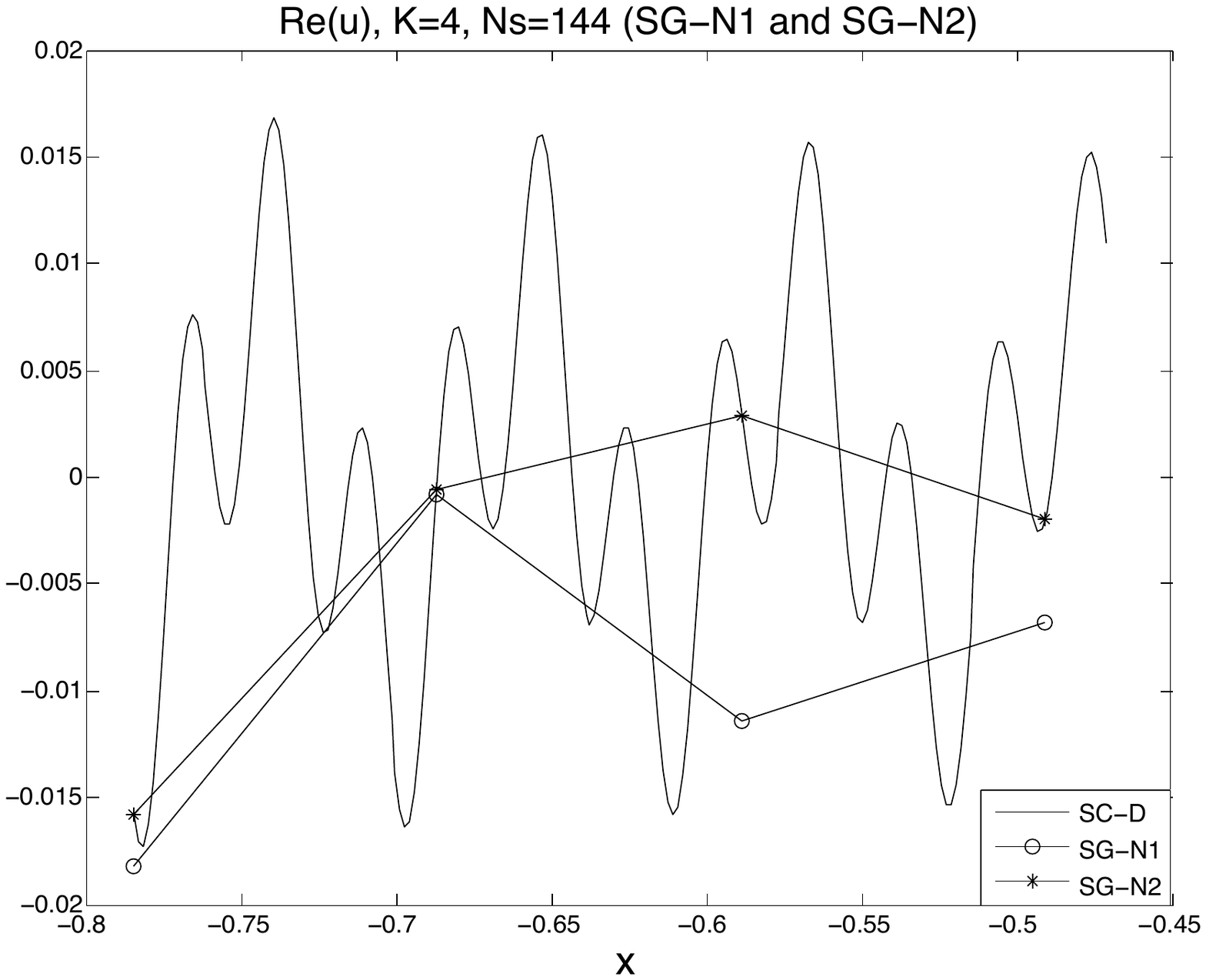}
 \end{subfigure}
 \begin{subfigure}{0.5\textwidth}
 \includegraphics[trim={0 10cm 0 0cm},width=1.0\textwidth, height=0.7\textwidth]{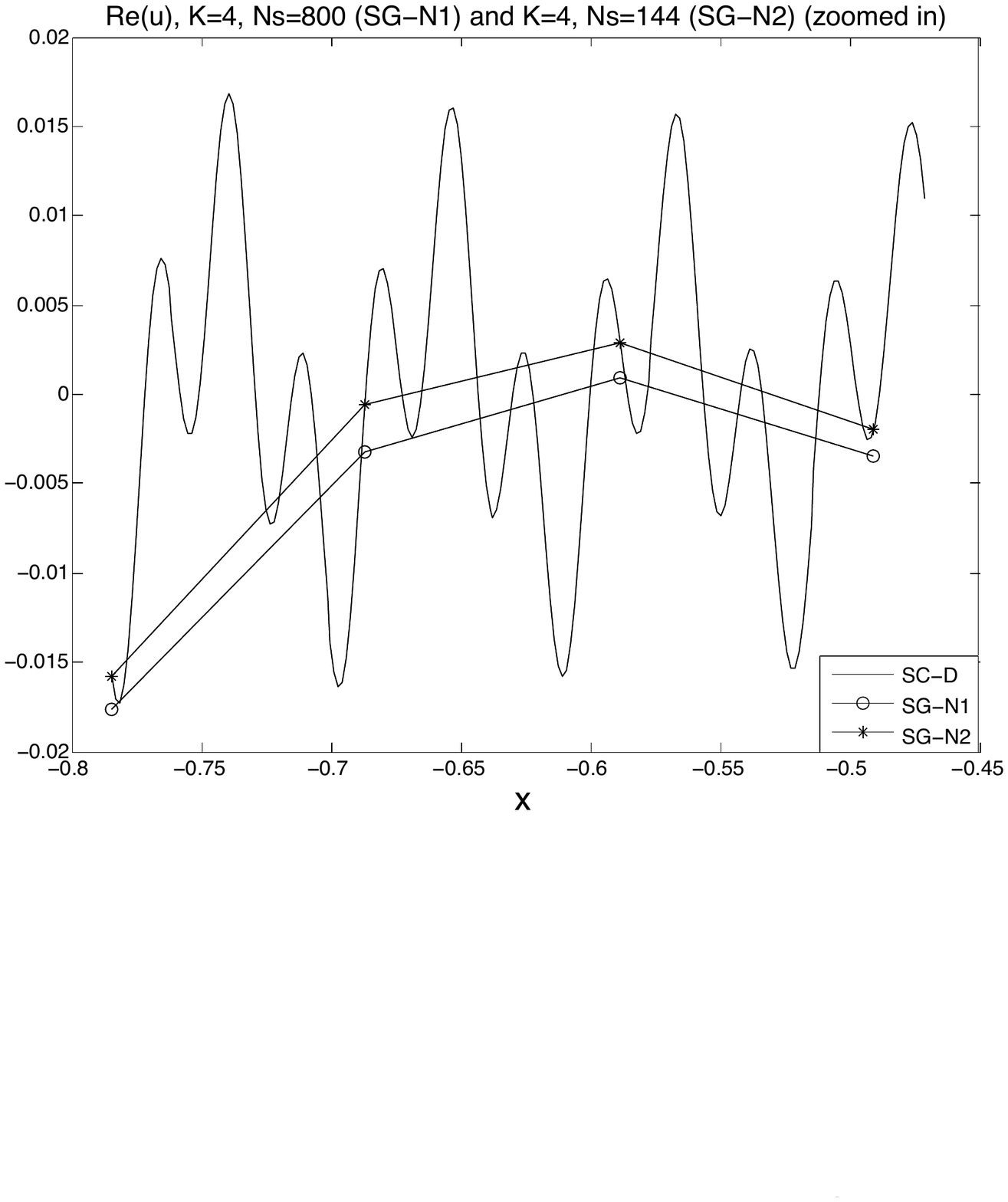}
\end{subfigure}
 \begin{subfigure}{0.5\textwidth}
  \includegraphics[trim={0 10cm 0 0},width=1.0\textwidth, height=0.7\textwidth]{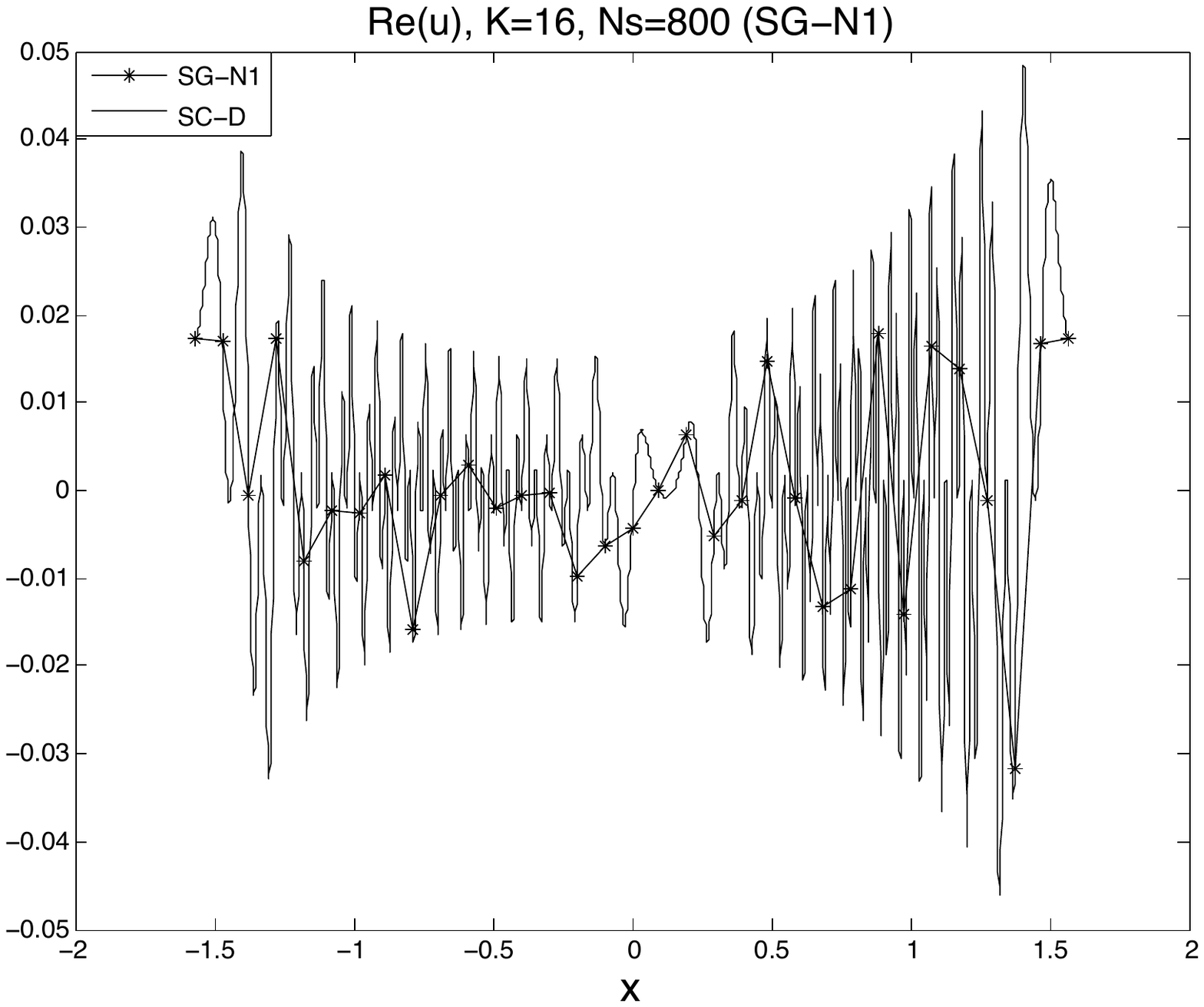}
\end{subfigure}
  \begin{subfigure}{0.5\textwidth}
 \includegraphics[trim={0 10cm 0 0},width=1.0\textwidth, height=0.7\textwidth]{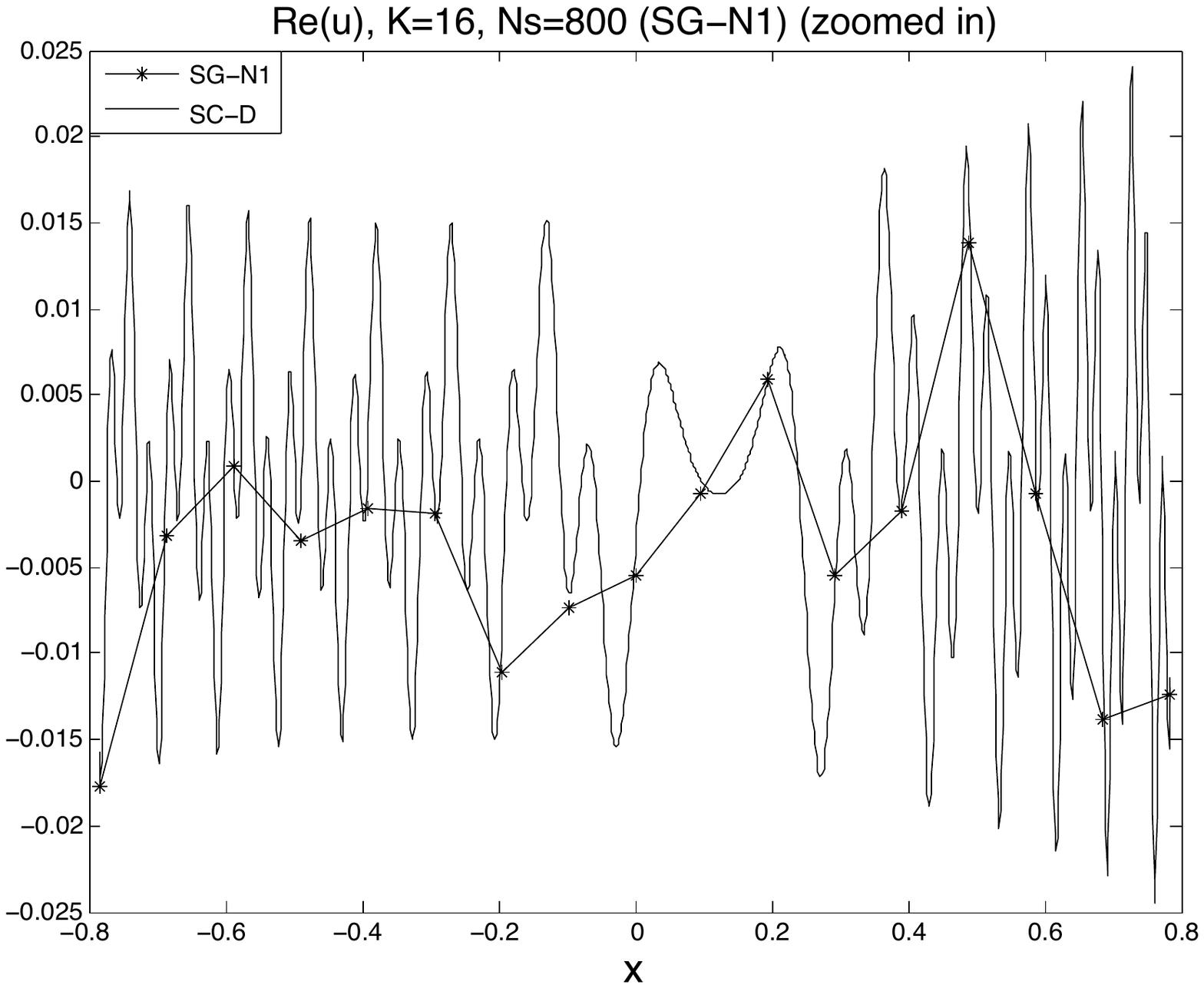}
 \end{subfigure}
  \begin{subfigure}{0.5\textwidth}
  \includegraphics[trim={0 10cm 0 0},width=1.0\textwidth, height=0.7\textwidth]{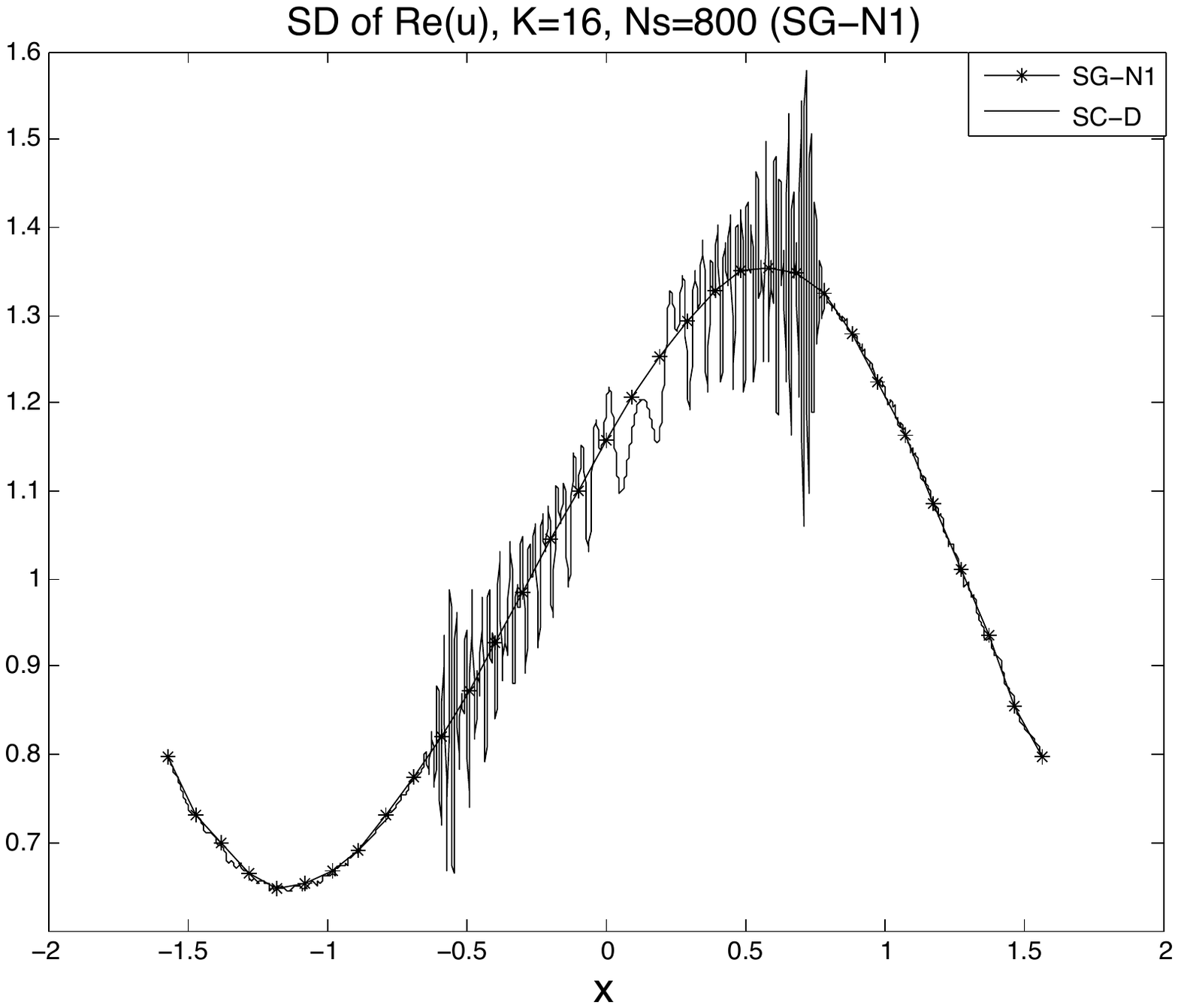}     
\end{subfigure}
  \begin{subfigure}{0.5\textwidth}
 \includegraphics[trim={0 10cm 0 0},width=1.0\textwidth, height=0.7\textwidth]{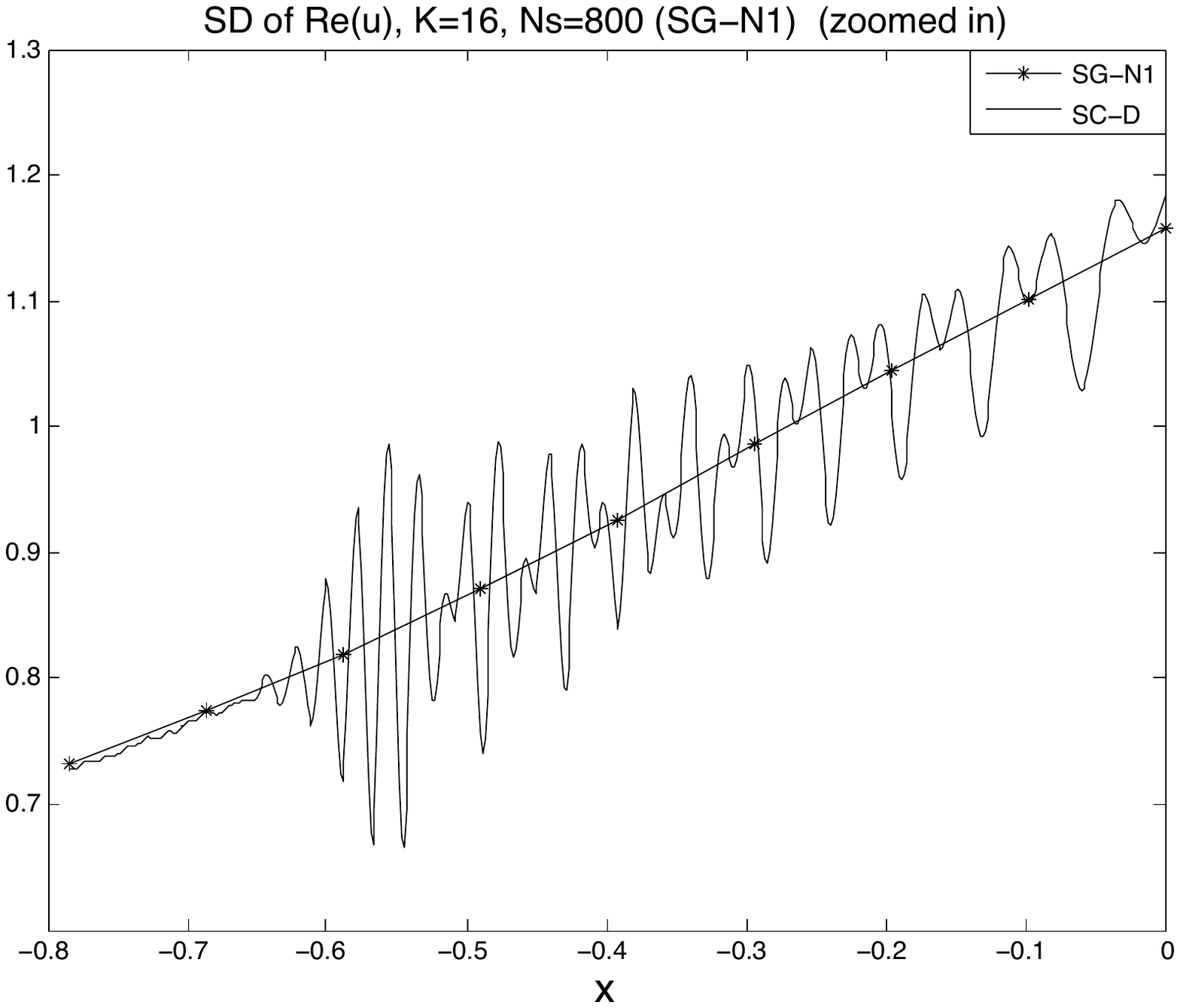}
 \end{subfigure}
  \begin{subfigure}{0.5\textwidth}
  \includegraphics[trim={0 10cm 0 0},width=1.0\textwidth, height=0.7\textwidth]{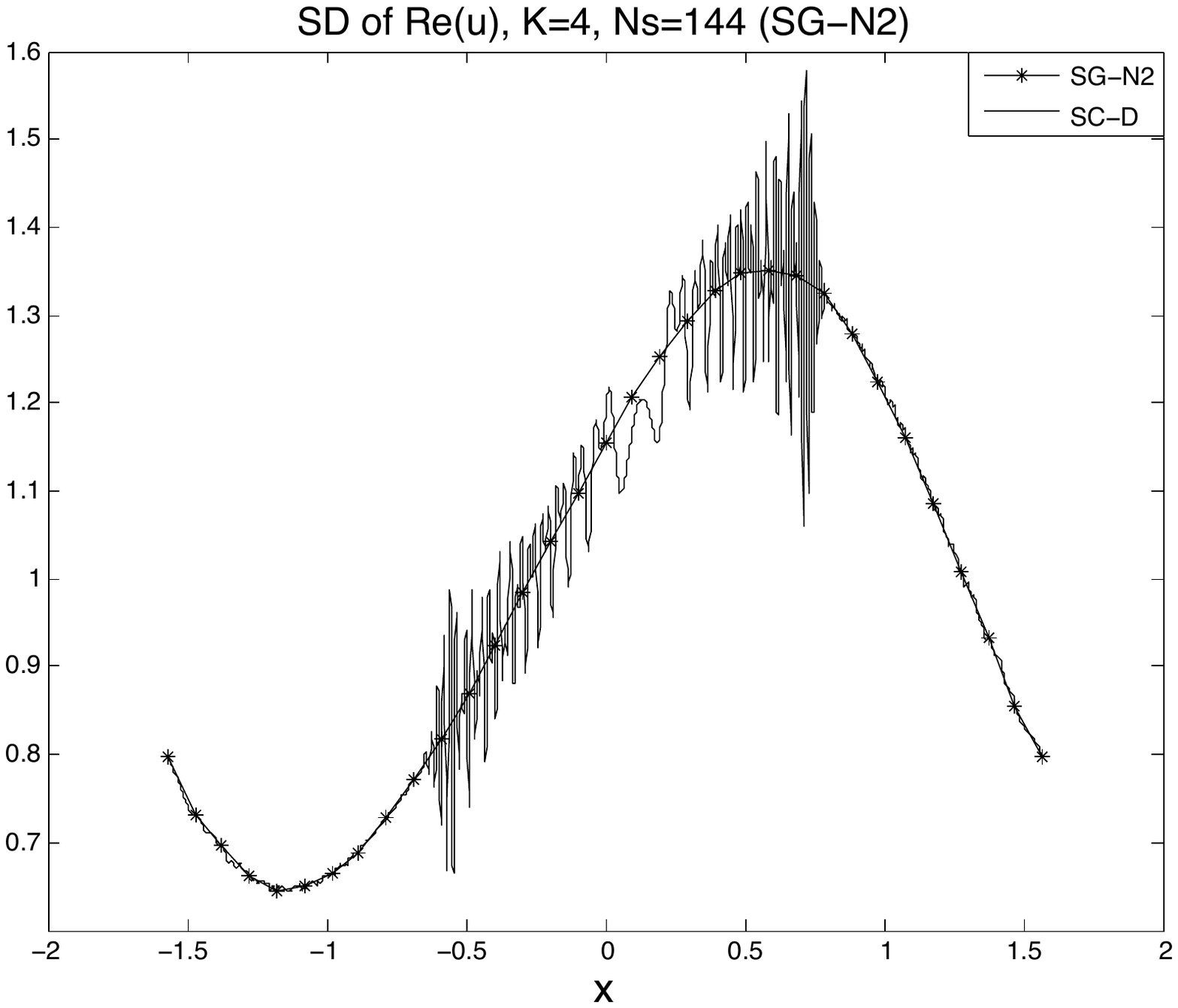}    
\end{subfigure}
  \begin{subfigure}{0.5\textwidth}
 \includegraphics[trim={0 10cm 0 0},width=1.0\textwidth, height=0.7\textwidth]{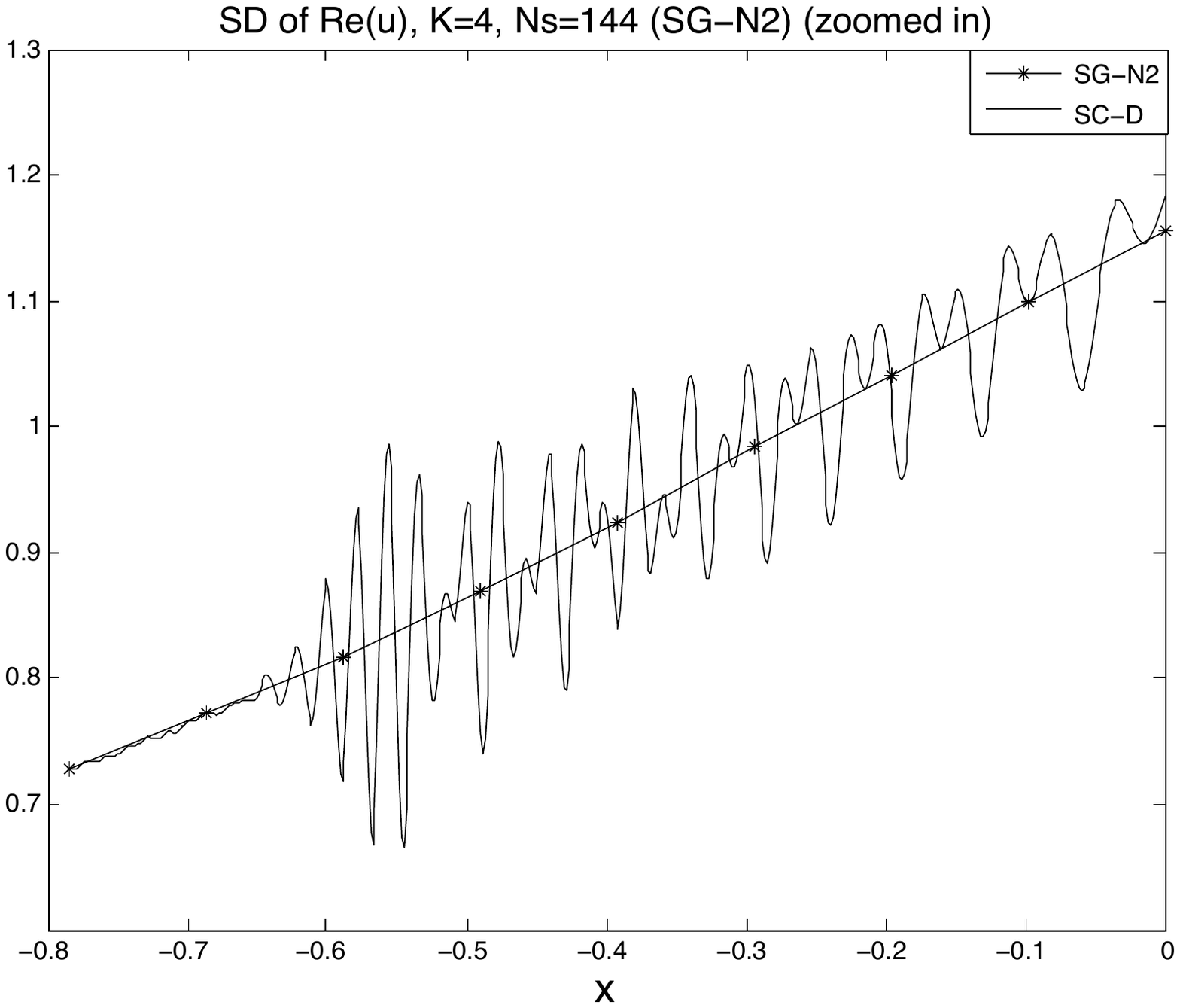}
 \end{subfigure}
\caption{Example $2.1$. Mean (first two rows) and standard deviation (last row) of real parts of $u$ (and their zoomed in solutions) at $t=0.25$, $\varepsilon=3\times 10^{-3}$. $\Delta x=\pi/32$, $\Delta t=0.01$ (gPC-SG-N1 and gPC-SG-N2), $\Delta x=\pi/2000$, $\Delta t=5\times10^{-5}$
 (gPC-SC-D). 
 Solid lines: reference solution by gPC-SC-D.}
 \label{Compare_1D}
 \end{figure}

Some comparisons between the gPC-SG-N1 and gPC-SG-N2 methods are shown in Figure \ref{Compare_1D} in the case of $\varepsilon=3\times 10^{-3}$.
We first plot the mean of real part of $u$ (zoomed in solutions) with $K=4$, $N_s=144$ in gPC-SG-N1 and gPC-SG-N2
and observes that gPC-SG-N2 can capture the correct mean and standard deviation, while gPC-SG-N1 can not.
Using $K=4$, $N_s=800$ in gPC-SG-N1 still does not give the good result, and further increasing the gPC order ($K=16$) enables
the mean and standard deviation of $\text{Re}(u)$ to capture the oscillations.
Thus gPC-SG-N2 does not require larger $K$ or $N_s$ when $\varepsilon$ is small, whereas larger $K$ and $N_s$ are needed in gPC-SG-N1.




\section{A semiclassical surface hopping model with random inputs}
\noindent
We consider a semiclassical surface hopping model in one dimension in space and momentum,
with random inputs,
\begin{align}
\left\{
\begin{array}{l}
\displaystyle\partial_t f^{+}+p \partial_{x}f^{+}-\partial_{x}(U+E) \partial_p f^{+}=\bar b^{i}f^{i}+b^{i}\bar f^{i}, \\[4pt]
\displaystyle\partial_t f^{-}+p \partial_{x}f^{-}-\partial_{x}(U-E) \partial_p f^{-}=-\bar b^{i}f^{i}-b^{i}\bar f^{i}, \\[4pt]
\displaystyle\partial_t f^{i}+p\partial_{x}f^{i}+\partial_{x}U\partial_p f^{i}=-i\frac{2E}{\varepsilon}f^{i}+b^{i}(f^{-}-f^{+})+
(b^{+}-b^{-})f^{i},
\label{SH old}
\end{array}\right.
\end{align}
where $(f^{+}(t,x,p,\bz), f^{-}(t,x,p,\bz), f^i(t,x,p,\bz))\in\mathbb R_+\times\mathbb R_+\times\mathbb C$, $(t,x,p)\in\mathbb R_{+}\times
\mathbb R\times\mathbb R$, and $b^{\pm}\in\mathbb C$, $b^i\in\mathbb C$, $E=E(x,\bz)\in\mathbb R$
are given functions depending on $(x,p,\bz)$.

This system is equipped with initial conditions
$f^\pm(0, x, p)=f^\pm_{in}(x,p)$ and $f^i(0, x, p)=f^i_{in}(x,p)$.

We also define the zeroth moments of $f^{\pm, i}$ as
\[
\rho^{\pm, i}(t, x, \bz) = \int_{\mathbb R} f^{\pm, i}(t, x, p, \bz) dp\,.
\]
Here $\rho^\pm $ denote the densities of particles in the two bands.

This model, introduced in \cite{CQJ}, was a semiclassical approximation to the
nucleonic Schr\"odinger system that arises from the Born-Oppenheimer
approximation with {\it non-adiabatic} phenomenon, in which particles
can ``hop'' from one potential energy surface to the other ones.
Here $f^\pm$ stand for the particle density distributions in the two
energy surfaces respectively, while $f^i$ is the off-diagonal entry of the
semiclassical Wigner matrix.
The right hand side of the system describes the interband transition between different potential energy surfaces with the $\it{energy \, gap}\, \Delta E=2E$, which can be random in reality.
In particular, the so-called $\it{avoided\, crossing}$ corresponds to where the minimum energy gap of order $\sqrt{\varepsilon}$.


\subsection{The direct method}
In this part and the following sections, we assume that the function $E(x,\bz)$ is given and random,
whereas no randomness is introduced to the initial function.

Introduce the following matrix
\begin{equation}\mathscr A=\text{diag}(-\partial_x (U+E), \partial_x (U-E), \partial_x U, \partial_x U), \end{equation}
The specific form of $b^i$ depends on the potential matrix in the
nucleonic Schr\"odinger equation. Here we consider $b^i=b^i(x,p)\in \mathbb R$, $b^{\pm}=0$,
which correspond to two specific potential matrices
considered in \cite{CQJ}. Let
\begin{align*}
\displaystyle
B=\begin{pmatrix}
0 & 0 & 2b^i & 0 \\
0 & 0 & -2b^i & 0\\
-b^i & b^i & 0 & 2E/\varepsilon \\
0 & 0 & -2E/\varepsilon & 0
\end{pmatrix},
\end{align*}
and write $f=(f^{+}, f^{-}, g:=\text{Re}(f^{i}), h:=\text{Im}(f^{i}))\in\mathbb R^4$,
then (\ref{SH old}) can be written as
\begin{equation}
\partial_t f+p\partial_x f+\mathscr A\partial_p f=Bf,  \label{direct}
\end{equation}
which will be solved using simple operator splitting method:

\noindent\textbf{1}. Solve $\partial_t f+p\partial_x f=0$ using spectral method in space and exact integration in time,  \\[4pt]
\textbf{2}. Solve $\partial_t f+\mathscr A\partial_p f=0$ using spectral method in velocity and exact integration in time,  \\[4pt]
\textbf{3}. Solve $\partial_t f=Bf$ to be specified below.
\\[2pt]

The first two steps are the same for both the direct and the
NGO-based method (to be introduced in the next section \ref{NGO_3D}).
For step 3, since the matrix $B$ now involves randomness, instead of solving the whole system exactly in time,
we use a Crank-Nicolson method in time coupled with a substitution method to save computations. More precisely, we have
\begin{align}
&\displaystyle\frac{f^{+,n+1}-f^{+,n}}{\Delta t}=\frac{1}{2}2b^i\left(g^n+g^{n+1}\right), \label{eqn1} \\[6pt]
&\displaystyle\frac{f^{-,n+1}-f^{-,n}}{\Delta t}=-\frac{1}{2}2b^i\left(g^n+g^{n+1}\right),  \label{eqn2}\\[6pt]
&\displaystyle\frac{g^{n+1}-g^n}{\Delta t}=\frac{1}{2}\left[-b^{i}f^{+,n}+b^{i}f^{-,n}+\frac{2E}{\varepsilon}h^{n}
 -b^{i}f^{+,n+1}+b^{i}f^{-,n+1}+\frac{2E}{\varepsilon}h^{n+1}\right], \label{eqn3} \\[6pt]
&\displaystyle\frac{h^{n+1}-h^n}{\Delta t}=-\frac{1}{2}\frac{2E}{\varepsilon}\left(g^{n}+g^{n+1}\right). \label{eqn4}
\end{align}
We used classical notations for the time approximation, for example $f^{+, n}\approx f^+(t_n)$.
From (\ref{eqn1}),
one solves $f^{+,n+1}$ in terms of $g^{n+1}$, and then substitute it into (\ref{eqn3}) to get $g^{n+1}$,
\begin{align}
\displaystyle
&\quad\left[1+(b^i)^2 (\Delta t)^2+\frac{E^2}{\varepsilon^2}(\Delta t)^2\right] g^{n+1} \nonumber\\[4pt]
&=g^{n}+\Delta t \left[-b^{i}f^{-,n}+b^{i}f^{-,n}+
\frac{2E}{\varepsilon}h^n-(b^i)^2 \Delta t g^n -\frac{E^2}{\varepsilon^2}\Delta t g^n \right]  \quad\,.
\label{F_3}
\end{align}
Then the other unknowns $f^{\pm, n+1}$ and $h^{n+1}$ can be obtained from (\ref{eqn1}), (\ref{eqn3}) and (\ref{eqn4}).

\subsection{The NGO based method}
\label{NGO_3D}
We first review the NGO-based method introduced in \cite{NGO}, wherein deterministic $E=E(x)$ and deterministic initial data is considered.

We introduce the augmented unknowns $(F^{\pm}, G, H)(t,x,p,\tau)$ satisfying
\begin{equation} f^{\pm}(t,x,p)=F^{\pm}(t,x,p,S(t,x,p)/\varepsilon), \quad f^i(t,x,p)=e^{iS(t, x, p)/\varepsilon}(G+iH)(t,x,p,S(t,x,p)/\varepsilon).
\end{equation}
The phase function $S(t,x,p)$ is designed to follow the main oscillations in the model. Assume periodicity in $S/\varepsilon$, $S$ solves
\begin{equation}
\partial_t S+p\partial_x S+\partial_x U\partial_p S=2E, \qquad S(0,x,p)=0.
\label{S_eqn}\end{equation}
Then $F^{\pm}, G$ and $H$ solve
\begin{align}
\left\{
\begin{array}{l}
\displaystyle
\partial_t F^{+}+p\partial_{x}F^{+}-\partial_{x}(U+E)\partial_{p}F^{+}=-\frac{\mathscr E^{+}}{\varepsilon}\partial_{\tau}F^{+}
+2b^i(G\cos\tau + H\sin\tau), \\[4pt]
\displaystyle
\partial_t F^{-}+p\partial_{x}F^{-}-\partial_{x}(U-E)\partial_{p}F^{-}=-\frac{\mathscr E^{-}}{\varepsilon}\partial_{\tau}F^{-}
-2b^i(G\cos\tau + H\sin\tau), \\[4pt]
\displaystyle
\partial_t G+p\partial_{x}G+\partial_{x}U\partial_{p}G=-\frac{2E}{\varepsilon}\partial_{\tau}G+b^{i}(F^{-}-F^{+})\cos\tau
, \\[4pt]
\displaystyle
\partial_t H+p\partial_{x}H+\partial_{x}U\partial_{p}H=-\frac{2E}{\varepsilon}\partial_{\tau}H+b^{i}(F^{-}-F^{+})\sin\tau,
\label{SH new}
\end{array}\right.
\end{align}
where we denoted
\begin{equation}
\mathscr E^{\pm}(t,x,p)=2E-\partial_{x}(2U\pm E)\partial_p S \,.
\end{equation}
The well-prepared initial conditions are given by (see \cite{NGO})
\begin{align}
\displaystyle
&\label{IC1}F^{+}(0,x,p,\tau)=f_{{in}}^{+}-i\frac{\varepsilon}{\mathscr E^{+}}\left(\bar b^{i}f_{{in}}^i (1-e^{-i\tau})-b^i\bar f_{{in}}^i
(1-e^{i\tau})\right), \\[4pt]
&\label{IC2} F^{-}(0,x,p,\tau)=f_{{in}}^{-}+i\frac{\varepsilon}{\mathscr E^{-}}\left(\bar b^{i}f_{{in}}^i (1-e^{-i\tau})-b^i\bar f_{{in}}^i
(1-e^{i\tau})\right), \\[4pt]
&\label{IC3}G(0,x,p,\tau)=\mbox{Re}(f_{{in}}^i)-\frac{\varepsilon}{2E}b^{i}(f_{{in}}^{+}-f_{ in}^{-})\sin\tau,  \\[4pt]
&\label{IC4}H(0,x,p,\tau)=\mbox{Im}(f_{{in}}^i)+\frac{\varepsilon}{2E}b^{i}(f_{{in}}^{+}-f_{ in}^{-})(\cos\tau-1),
\end{align}
where $f^\pm_{in}, f^i_{in}$ are the initial data of \eqref{SH old}.

One can also develop the gPC-SG-D and gPC-SG-N1 methods as in section \ref{scalar}. Since our focus is on gPC-SG-N2, we will just put gPC-SG-N1 and gPC-SG-D
in the appendix for future reference.

\subsubsection {The gPC-SG-N2 for the system with random inputs}
\label{N2_3D}
Similar to the discussion in section \ref{K_eps}, which allows the gPC order independent of $\varepsilon$, we implement the gPC-SG-N2 scheme.
First, one focuses on the approximation of the phase equation \eqref{S_eqn}.
We insert the approximated solution
$$S(t,x,p,\bz)\approx \sum_{j=1}^{K}\tilde S_{j}(t,x,p)\tilde \psi_{j}(\bz)$$
into (\ref{S_eqn}) and conduct the Galerkin projection.
We denote by
\begin{align*}
\displaystyle \vv S(t,x,p)=(\tilde{S}_1(t,x,p), \cdots, \tilde{S}_K(t,x,p))^{T},
\end{align*}
the vector of unknown satisfying
\begin{equation}
\partial_t\vv S+p\partial_x\vv S+\partial_x U\partial_p\vv S
=2\int_{I_{\bz}}E(x,\bz)\tilde{\boldsymbol\psi}(\bz) \pi(\bz)d\bz, \qquad \vv S(0,x,p)=0,
\label{S_3D}
\end{equation}
with $\tilde{\boldsymbol\psi}(\bz)=(\tilde \psi_1(\bz), \cdots,\tilde  \psi_K(\bz))$ defined in $(\ref{psi})$.
A $4$-th order Runge-Kutta in time and pseudo-spectral method in space and velocity is used to solve $\vv S$.

Define the augmented functions
\begin{align*}&F^{+}(t,x,p,\tau,\bz)=W_1(S(t,x,p,\bz),x,p,\tau,\bz),\quad F^{-}(t,x,p,\tau,\bz)=W_2(S(t,x,p,\bz),x,p,\tau,\bz), \\[4pt]
& G(t,x,p,\tau,\bz)=W_3(S(t,x,p,\bz),x,p,\tau,\bz), \qquad H(t,x,p,\tau,\bz)=W_4(S(t,x,p,\bz),x,p,\tau,\bz)\,.
\end{align*}
Then $$\partial_t F^{+}=\partial_s W_1 \partial_t S, \quad \partial_x F^{+}=\partial_s W_1\partial_x S+\partial_x W_1, \quad
\partial_p F^{+}=\partial_s W_1 \partial_p S+\partial_p W_1\,.$$
Denote $${\mathscr E^{\pm}}(t,x,p,\bz)=2E(x,\bz)-\partial_x(2U(x)\pm E(x,\bz))\partial_p S(t,x,p,\bz)\,.$$
By (\ref{S_eqn}) and (\ref{SH new}), we have
\begin{align}
\label{W system}
\left\{
\begin{array}{l}
\displaystyle
\partial_s W_1 +\frac{p}{\tilde{\mathscr E^+}(s,x,p,\bz)}\partial_x W_1 -\frac{\partial_x (U+E)}{\tilde{\mathscr E^+}(s,x,p,\bz)}\partial_p W_1 = -\frac{1}{\varepsilon}\partial_{\tau}W_1 +
\frac{2b^i}{\tilde{\mathscr E^+}(s,x,p,\bz)}(W_3\cos\tau + W_4\sin\tau), \\[6pt]
\displaystyle \partial_s W_2 +\frac{p}{\tilde{\mathscr E^-}(s,x,p,\bz)}\partial_x W_2 -\frac{\partial_x (U-E)}{\tilde{\mathscr E^-}(s,x,p,\bz)}\partial_p W_2 = -\frac{1}{\varepsilon}\partial_{\tau}W_2 -\frac{2b^i}{\tilde{\mathscr E^-}(s,x,p,\bz)}(W_3\cos\tau + W_4\sin\tau), \\[6pt]
\displaystyle \partial_s W_3 + \frac{p}{2E(x,\bz)}\partial_x W_3 + \frac{\partial_x U}{2E(x,\bz)}\partial_p W_3 = -\frac{1}{\varepsilon}\partial_{\tau}W_3 + \frac{b^i}{2E(x,\bz)}(W_2-W_1)\cos\tau, \\[6pt]
 \displaystyle \partial_s W_4 + \frac{p}{2E(x,\bz)}\partial_x W_4 + \frac{\partial_x U}{2E(x,\bz)}\partial_p W_4 = -\frac{1}{\varepsilon}\partial_{\tau}W_4 +\frac{b^i}{2E(x,\bz)}(W_2-W_1)\sin\tau.
\end{array}\right.
\end{align}
How to obtain $\tilde{\mathscr E^{\pm}}(s,x,p,\bz)$ will be introduced in next subsection \ref{N2_3D}.

Now we focus on the approximation of the profiles $(W_1,W_2,W_3,W_4)$.
We insert the Galerkin approximation of each component  $$
(W_1,W_2,W_3,W_4)(s,x,p,\bz)=\sum_{j=1}^{K} ((\tilde W_1)_j, (\tilde W_2)_j, (\tilde W_3)_j, (\tilde W_4)_j)(s,x,p)\tilde \psi_{j}(\bz),
$$
into (\ref{W system}) and conduct the gPC Galerkin projection.
Denote the associated gPC coefficients vector of $(\vv W_1, \vv W_2, \vv W_3, \vv W_4)$ as previously, for instance for $\vv W_1$,
\begin{equation}\vv W_1=((\tilde W_1)_1, \cdots, (\tilde W_1)_K)^T. \label{gPC_vector} \end{equation}
Then, we get the gPC system
\begin{align}
\label{W_gPC1}
\left\{
\begin{array}{l}
\displaystyle
\partial_s \vv W_1 + p {\cal J}\partial_x \vv W_1 - {\cal C} \partial_p \vv W_1 = -\frac{1}{\varepsilon}\partial_{\tau}\vv W_1 +
2b^{i}{\cal J} (\vv W_3 \cos\tau +\vv W_4 \sin\tau), \\[4pt]
\displaystyle
\partial_s \vv W_2 + p {\cal L} \partial_x \vv W_2  + {\cal C}\partial_p \vv W_2 = -\frac{1}{\varepsilon}\partial_{\tau}\vv W_2 -
2b^{i}{\cal L} (\vv W_3 \cos\tau + \vv W_4 \sin\tau), \\[4pt]
\displaystyle
\partial_s \vv W_3 + p {\cal H}\partial_x \vv W_3 = -\frac{1}{\varepsilon}\partial_{\tau}\vv W_3 + b^{i}{\cal H}(\vv W_2 -
\vv W_1)\cos\tau, \\[4pt]
\displaystyle
\partial_s \vv W_4 + p {\cal H}\partial_x \vv W_4 = -\frac{1}{\varepsilon}\partial_{\tau}\vv W_4  + b^{i}{\cal H}(\vv W_2-
\vv W_1)\sin\tau,
\end{array}\right.
\end{align}
where matrices ${\cal J}$, ${\cal L}$, ${\cal C}$, ${\cal H}$ are given by
\begin{align*}
& {\cal J}_{mn}(s,x,p) = \int_{I_{\bz}}\frac{1}{\tilde{\mathscr E^+}(s,x,p,\bz)}\psi_m(\bz)\psi_n(\bz)\pi(\bz)d\bz, \\[2pt]
& {\cal L}_{mn}(s,x,p)= \int_{I_{\bz}}\frac{1}{\tilde{\mathscr E^-}(s,x,p,\bz)}\psi_m(\bz)\psi_n(\bz)\pi(\bz)d\bz, \\[2pt]
& {\cal C}_{mn}(s,x,p) =\int_{I_{\bz}}\frac{\partial_x E(x,\bz)}{\tilde{\mathscr E^+}(s,x,p,\bz)}\psi_m(\bz)\psi_n(\bz)\pi(\bz)d\bz, \\[2pt]
& {\cal H}_{mn}(x) =\int_{I_{\bz}}\frac{1}{2E(x,\bz)}\psi_m(\bz)\psi_n(\bz)\pi(\bz)d\bz.
\end{align*}

The initial conditions of $\vv W_1$, $\vv W_2$, $\vv W_3$, $\vv W_4$ are
\begin{align*}
& \vv W_1(0,x,p,\tau)=\vv F^{+}(0,x,p,\tau), \qquad
\vv W_2(0,x,p,\tau)=\vv F^{-}(0,x,p,\tau), \\[2pt]
& \vv W_3(0,x,p,\tau)=\vv G(0,x,p,\tau),  \qquad\quad
\vv W_4(0,x,p,\tau)=\vv H(0,x,p,\tau).
\end{align*}
According to (\ref{IC1})-(\ref{IC4}), by the gPC-SG method, the initial conditions for gPC coefficient vectors $\vv F^{+}$, $\vv F^{-}$, $\vv G$ and
$\vv H$ are given by
\begin{align}
\displaystyle
&\vv F^{+}(0,x,p,\tau)=f_{{in}}^{+}\vv Q_1-i\varepsilon\left(\bar b^{i}f_{{in}}^i (1-e^{-i\tau})-b^i\bar f_{{in}}^i
(1-e^{i\tau})\right)\vv Q_2, \\[4pt]
&\vv F^{-}(0,x,p,\tau)=f_{{in}}^{-}\vv Q_1+i\varepsilon\left(\bar b^{i}f_{{in}}^i (1-e^{-i\tau})-b^i\bar f_{{in}}^i
(1-e^{i\tau})\right)\vv Q_3, \\[4pt]
& \vv G(0,x,p,\tau)=\mbox{Re}(f_{{in}}^{i})\vv Q_1-\varepsilon b^{i}(f_{{in}}^{+}-f_{ in}^{-})\vv Q_4 \sin\tau,\\[4pt]
&\vv H(0,x,p,\tau)=\mbox{Im}(f_{{in}}^{i})\vv Q_1+\varepsilon b^{i}(f_{{in}}^{+}-f_{ in}^{-})\vv Q_4 (\cos\tau-1),
\end{align}
where the vectors $\vv Q_1$, $\vv Q_2$, $\vv Q_3$, $\vv Q_4$ are defined by
\begin{align}
\displaystyle
&\vv Q_1=\int_{I_{\bz}}\boldsymbol\psi(\bz)\pi(\bz)d\bz\, , &\vv Q_2(0,x,p)=\int_{I_{\bz}}\frac{\boldsymbol\psi(\bz)\pi(\bz)}{\mathscr E^{+}(0,x,p,\bz)}d\bz\,, \nonumber\\
&\vv Q_3(0,x,p)=\int_{I_{\bz}}\frac{\boldsymbol\psi(\bz)\pi(\bz)}{\mathscr E^{-}(0,x,p,\bz)}d\bz\, ,& \vv Q_4(x)=\int_{I_{\bz}}\frac{\boldsymbol\psi(\bz)\pi(\bz)}{2E(x,\bz)}d\bz \,.\nonumber
\end{align}

\subsubsection{The fully discrete scheme for gPC-SG-N2}
{\bf Monotonicity of $S$ in terms of $t$}

Assume $U=U(x)\in\mathbb R$, $E=E(x,\bz)>0$,
$$\partial_t S + p\partial_x S+ \partial_x U \partial_p S= 2E(x,\bz), \quad S_{{in}}(x,p,\bz)=0. $$
By the method of characteristics,
$$\frac{dx}{dt}=p, \qquad \frac{dp}{dt}=\partial_x U, \qquad x(0)=x_0, \qquad p(0)=p_0, $$
then $$x=X(t, x_0, p_0), \qquad p=P(t,x_0, p_0), $$
thus $$\frac{d}{dt}S(t,x,p,\bz)=2E(X(t,x_0,p_0)), $$
and the analytic solution is given by
\begin{align*} & S(t,x,p,\bz)=S_{\text{in}}(x_0, p_0, \bz)+2\int_0^t E(X(\mu,x_0, p_0, \bz))d\mu \\[4pt]
&\qquad\qquad\quad=2\int_0^t E(X(\mu,x_0, p_0, \bz))d\mu\,. \end{align*}
Therefore $S$ is an increasing function of $t$ for each $(x,p,\bz)$, since $E>0$. \\[2pt]

The time discretization of $W_1$, $W_2$, $W_3$, $W_4$ was defined in section \ref{N2_1D}.
Denote ${\mathscr E^{\pm}}(t,x,p,\bz)=\tilde{\mathscr E^{\pm}}(S(t,x,p,\bz),x,p,\bz)$. To find the values of $\tilde{\mathscr E^{\pm}}$ at $s_l$, namely
$$\tilde{\mathscr E^{\pm}}(s_l,x,p,\bz)={\mathscr E^{\pm}}(t^{\star},x,p,\bz), $$
where $t^{\star}=S^{-1}(s_l)$, for each $x_j$, $p_k$, and quadrature points $\bz^{(q)}$, $q=1, \cdots, N_s$, we apply an interpolation step shown here.

Search for the time interval $[t_{n_k}, t_{n_k+1}]$ such that $s_l$ falls between the interval $[S(t_{n_k}), S(t_{n_k+1})]$, thus $t^{\star}
\in [t_{n_k}, t_{n_k+1}]$, since $S$ is an increasing function of $t$. Linear interpolation is used to find $t^{\star}$. Denote $t_{n_k}=t^{(1)}$, $t_{n_k+1}=t^{(2)}$, and $S(t_{n_k})=S^{(1)}$, $S(t_{n_k+1})=S^{(2)}$. By the Lagrange interpolation,
$$s_l=S^{(1)}\frac{t^{\star}-t^{(2)}}{t^{(1)}-t^{(2)}}+S^{(2)}\frac{t^{\star}-t^{(1)}}{t^{(2)}-t^{(1)}}\,,$$
which gives $t_{\star}$. The values ${\mathscr E^{\pm}}$ at all $t_n$ have been obtained, one can use linear interpolation to approximate the value of
${\mathscr E^{\pm}}$ at $t^{\star}$, for each $x_j$, $p_k$, $\bz^{(q)}$.

Using the Gauss quadrature rule with quadrature points $\bz^{(q)}$ and the corresponding weights $\omega_q$, $q=1,\cdots, N_g$,
we update the matrices ${\cal J}$, ${\mathscr L}$, ${\cal C}$ at each $(s_l, x_j, p_k)$, and ${\cal H}$ at each $x_j$,
\begin{align*}
& {\cal J}_{mn}(s_l, x_j, p_k) = \sum_{q=1}^{N_g}\frac{1}{\tilde{\mathscr E^+}(s_l, x_j, p_k, \bz^{(q)})}\psi_m(\bz^{(q)})\psi_n(\bz^{(q)})\omega_q,  \\[4pt]
& {\cal L}_{mn}(s_l, x_j, p_k) = \sum_{q=1}^{N_g}\frac{1}{\tilde{\mathscr E^-}(s_l, x_j, p_k, \bz^{(q)})}\psi_m(\bz^{(q)})\psi_n(\bz^{(q)})\omega_q,    \\[4pt]
& {\cal C}_{mn}(s_l, x_j, p_k) =\sum_{q=1}^{N_g}\frac{\partial_x E(x_j, \bz^{(q)})}{\tilde{\mathscr E^+}(s_l, x_j, p_k, \bz^{(q)})}\psi_m(\bz^{(q)})\psi_n(\bz^{(q)})\omega_q, \\[4pt]
& {\cal H}_{mn}(x_j) =\sum_{q=1}^{N_g}\frac{1}{2E(x_j, \bz^{(q)})}\psi_m(\bz^{(q)})\psi_n(\bz^{(q)})\omega_q \,.
\end{align*}
\\[2pt]

\textbf{Step 1} \quad
We solve the transport part in $x$ for quantities $\vv W_1$, $\vv W_2$, $\vv W_3$, $\vv W_4$ in (\ref{W_gPC1}),
\begin{align*}
& \partial_s \vv W_1 + p{\cal J}\partial_x \vv W_1=0, \qquad \partial_s \vv W_2 + p{\cal J}\partial_x \vv W_2 =0, \\[4pt]
& \partial_s \vv W_3 + p{\cal H}\partial_x \vv W_3=0, \qquad \partial_s \vv W_4 + p{\cal H}\partial_x \vv W_4=0.
\end{align*}
This transport step is treated similarly as in section \ref{gPC_Dir_1D}, where we use a pseudo-spectral method in space and a three-stage Runge-Kutta method in time.
\\[2pt]

\textbf{Step 2} \quad
We solve the transport part in $p$ for $W_1$, $W_2$,
\begin{align*}
\partial_s \vv W_1-{\cal C}\partial_p \vv W_1=0, \qquad \partial_s \vv W_2+{\cal C}\partial_p\vv W_2=0\,.
\end{align*}
Similar procedure is taken as in the previous step.
\\[2pt]

\textbf{Step 3}\quad
We now solve the non-singular source part
\begin{align*}
&\partial_s \vv W_1 = 2b^{i} {\cal J}(\vv W_3 \cos\tau +\vv W_4\sin\tau), \\[4pt]
&\partial_s \vv W_2 = -2b^{i}{\cal L}(\vv W_3\cos\tau +\vv W_4\sin\tau), \\[4pt]
&\partial_s \vv W_3 = b^{i}{\cal H}(\vv W_2-\vv W_1)\cos\tau, \\[4pt]
&\partial_s \vv W_4 = b^{i}{\cal H}(\vv W_2-\vv W_1)\sin\tau\,.
\end{align*}

We use the forward Euler method,
\begin{align*}
&\vv W_1^{n+1}=\vv W_1^n + 2b^i \Delta s{\cal J}(\vv W_3^n \cos\tau +\vv W_4^n \sin\tau), \\[4pt]
&\vv W_2^{n+1}=\vv W_2^n - 2b^i \Delta s{\cal L}(\vv W_3^n\cos\tau +\vv W_4^n\sin\tau), \\[4pt]
&\vv W_3^{n+1}=\vv W_3^n + b^i \Delta s{\cal H}(\vv W_2^n-\vv W_1^n)\cos\tau, \\[4pt]
&\vv W_4^{n+1}=\vv W_4^n + b^i \Delta s{\cal H}(\vv W_2^n-\vv W_1^n)\sin\tau\,.
\end{align*}
Other Runge-Kutta methods can also be used to solve this ODE system.
\\[2pt]

\textbf{Step 4}\quad
Finally, we solve the highly oscillatory part
\begin{align*}
&\partial_s \vv W_1 = -\frac{1}{\varepsilon}\partial_{\tau}\vv W_1, \qquad \partial_s \vv W_2 = -\frac{1}{\varepsilon}\partial_{\tau}\vv W_2, \\[2pt]
&\partial_s \vv W_3 = -\frac{1}{\varepsilon}\partial_{\tau}\vv W_3, \qquad \partial_s \vv W_4 = -\frac{1}{\varepsilon}\partial_{\tau}\vv W_4.
\end{align*}
We use the Fourier transform $\hat W_1(\zeta)$ of $W_1(\tau)$ in the $\tau$ variable where $\zeta$ is the corresponding Fourier variable, then
$$ \partial_s \vv{\hat W_1} = -\frac{i\zeta}{\varepsilon}\vv{\hat W_1}, $$
which is solved by the backward Euler method in time,
$$\vv{\hat W_1}^{n+1}(x_j,p_k) = \left(I+\frac{i\zeta\Delta s}{\varepsilon}I\right)^{-1}\vv{\hat W_1}^n(x_j,p_k)\,.$$
The same method is used for calculating $\vv W_2$, $\vv W_3$, $\vv W_4$.

\section{Numerical examples for the surface hopping model}
\begin{itemize}
\item \textbf{Example $4.1$}
\end{itemize}
\noindent
Consider $x,p\in[-2\pi,2\pi]$. We use the similar data given in the numerical example in \cite{NGO},
with the following initial conditions,
\begin{align*}
\displaystyle
&f^{+}(t=0,x,p)=f^{-}(t=0,x,p)=(1+0.5\cos(x))\frac{e^{-p^2/2}}{\sqrt{2\pi}}, \\[4pt]
&f^{i}(t=0,x,p)=\left[1+0.5\sin(x))+i(1+0.5\cos(x))\right]\frac{e^{-p^2/2}}{\sqrt{2\pi}}.
\end{align*}
The expressions for $E$, $b_i$, and $b^{\pm}$ are given by
$$E(x,\bz)=(1-\cos(x/2)+\sqrt{\varepsilon})(1+0.5\bz), \quad b_i(x,p)=-\frac{1}{2}\sin(p+1), \quad b^{\pm}=0.$$
Periodic boundary conditions are considered in $x$ while the $p$ domain is chosen large enough so $f^{\pm, i}$ all vanish
outside the domain (thus a periodic boundary condition in $p$ can be used).  Without
loss of generality we set $U=0$. For all the following tests, we choose $N_{\tau}=8$, $N_p=32$ and $N_c=32$ quadrature points for gPC-SC methods (with a very small mesh size in $x$ so this solution is used as the reference solutions).

In Figure \ref{gPC_N2}, for gPC-SG-N2 scheme, even if the mesh
size is much larger than the wave length $\varepsilon=5\times 10^{-3}$, the solution $f^{\pm, i}$ (the mean and standard deviation), as well
as $\rho^{\pm, i}$, still agrees with
the reference solution at the grid points, despite a moderate gPC order $K=4$.

\begin{figure}[H]
 \begin{subfigure}{0.5\textwidth}
   \centering
 \includegraphics[trim={0 9cm 0 0},width=1.05\textwidth, height=0.8\textwidth]{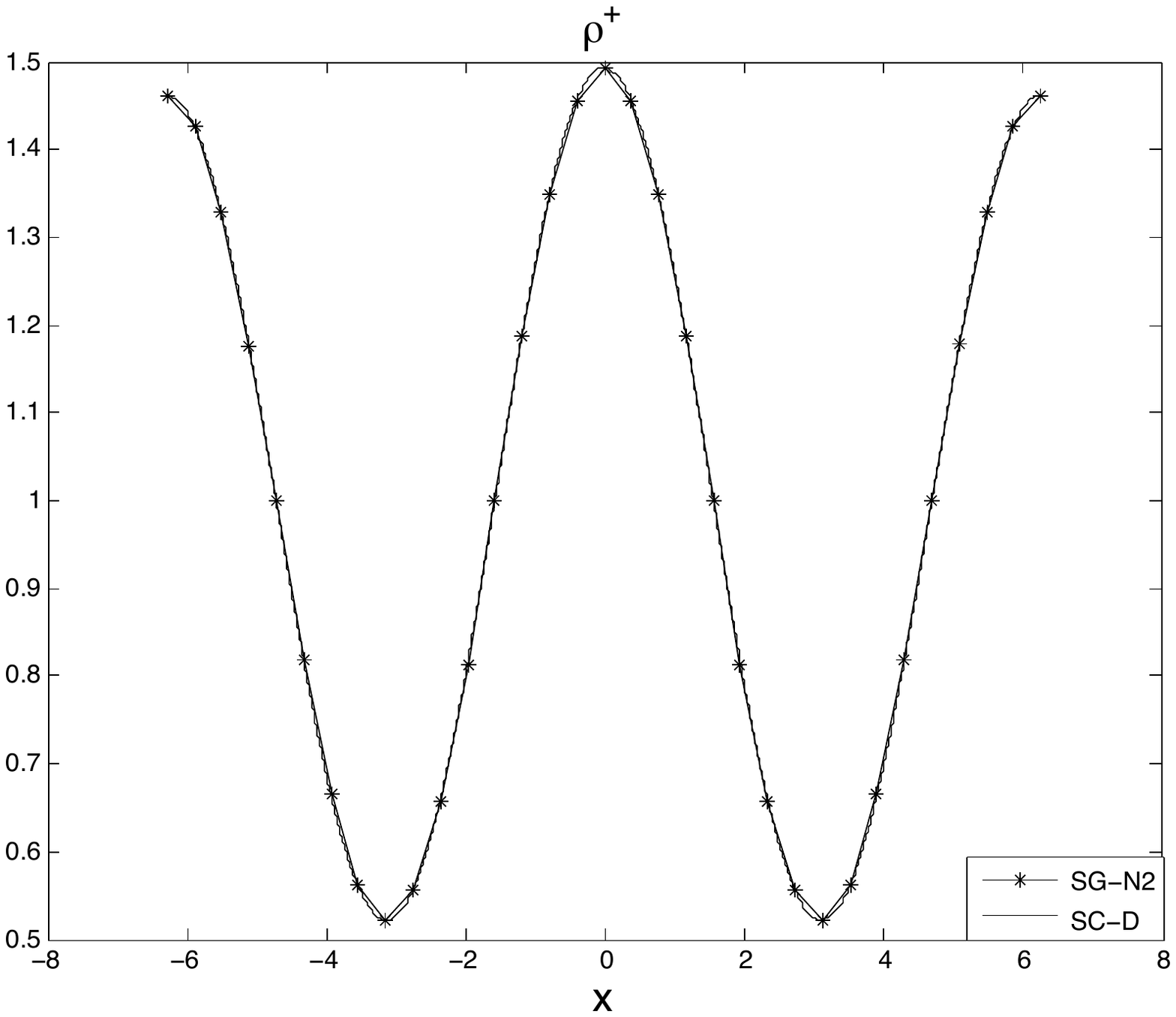}
  \end{subfigure}
   \begin{subfigure}{0.5\textwidth}
   \centering
 \includegraphics[trim={0 9cm 0 0},width=1.05\textwidth, height=0.8\textwidth]{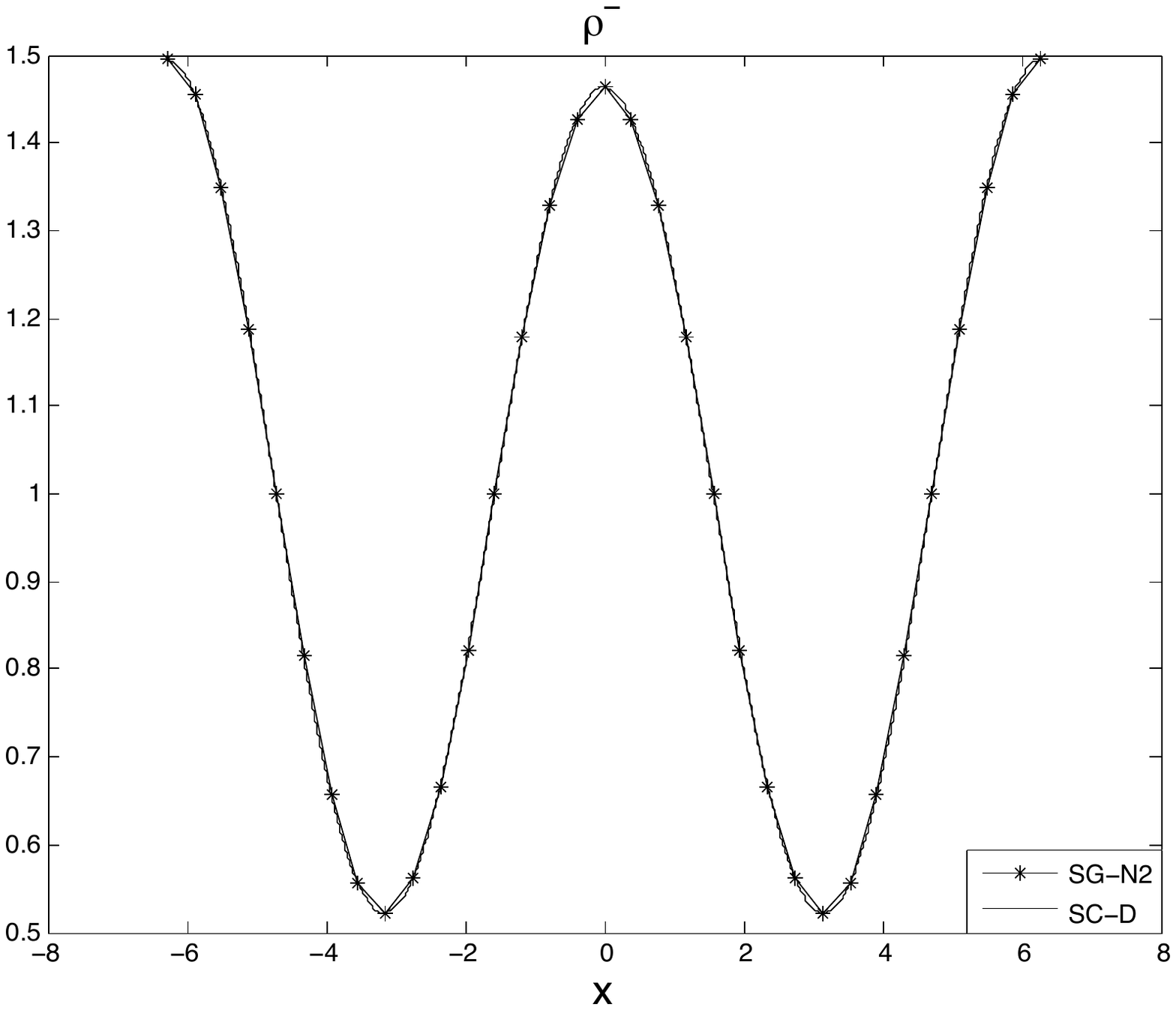}
  \end{subfigure}
\begin{subfigure}{0.5\textwidth}
   \centering
 \includegraphics[trim={0 9cm 0 0},width=1.05\textwidth, height=0.8\textwidth]{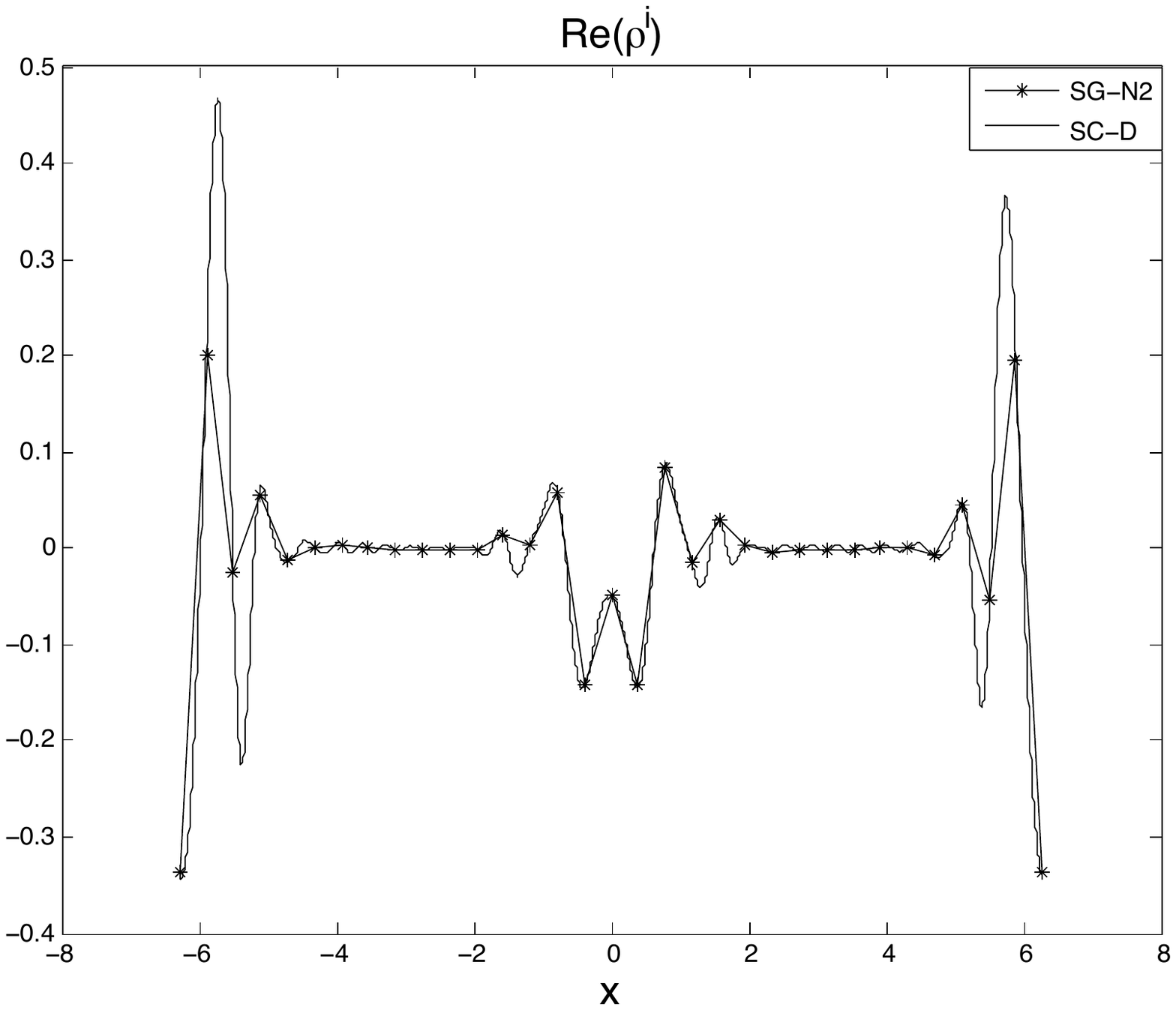}
  \end{subfigure}
 \begin{subfigure}{0.5\textwidth}
   \centering
 \includegraphics[trim={0 9cm 0 0},width=1.05\textwidth, height=0.8\textwidth]{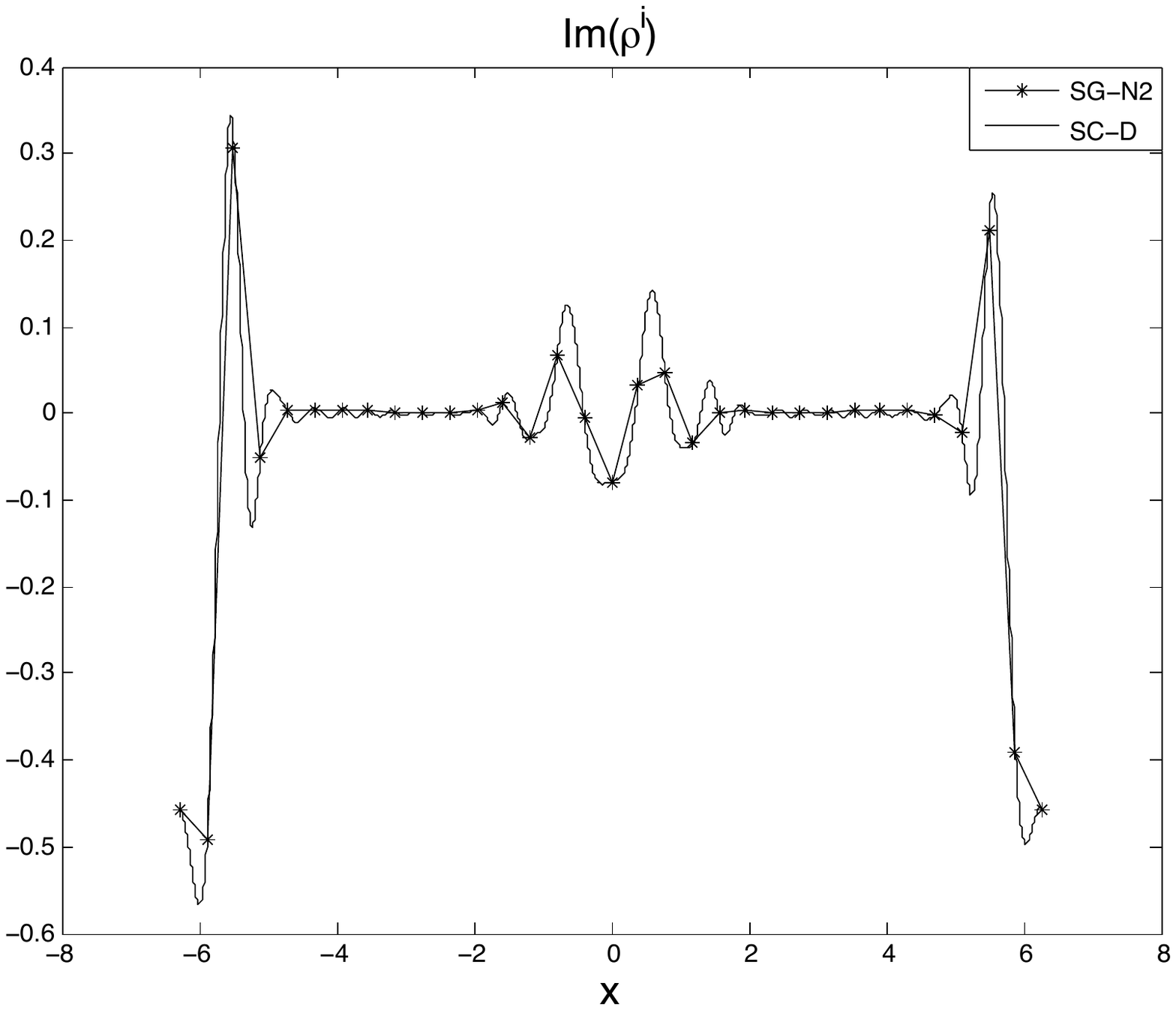}
  \end{subfigure}
  \begin{subfigure}{0.5\textwidth}
   \centering
 \includegraphics[trim={0 9cm 0 0},width=1.05\textwidth, height=0.8\textwidth]{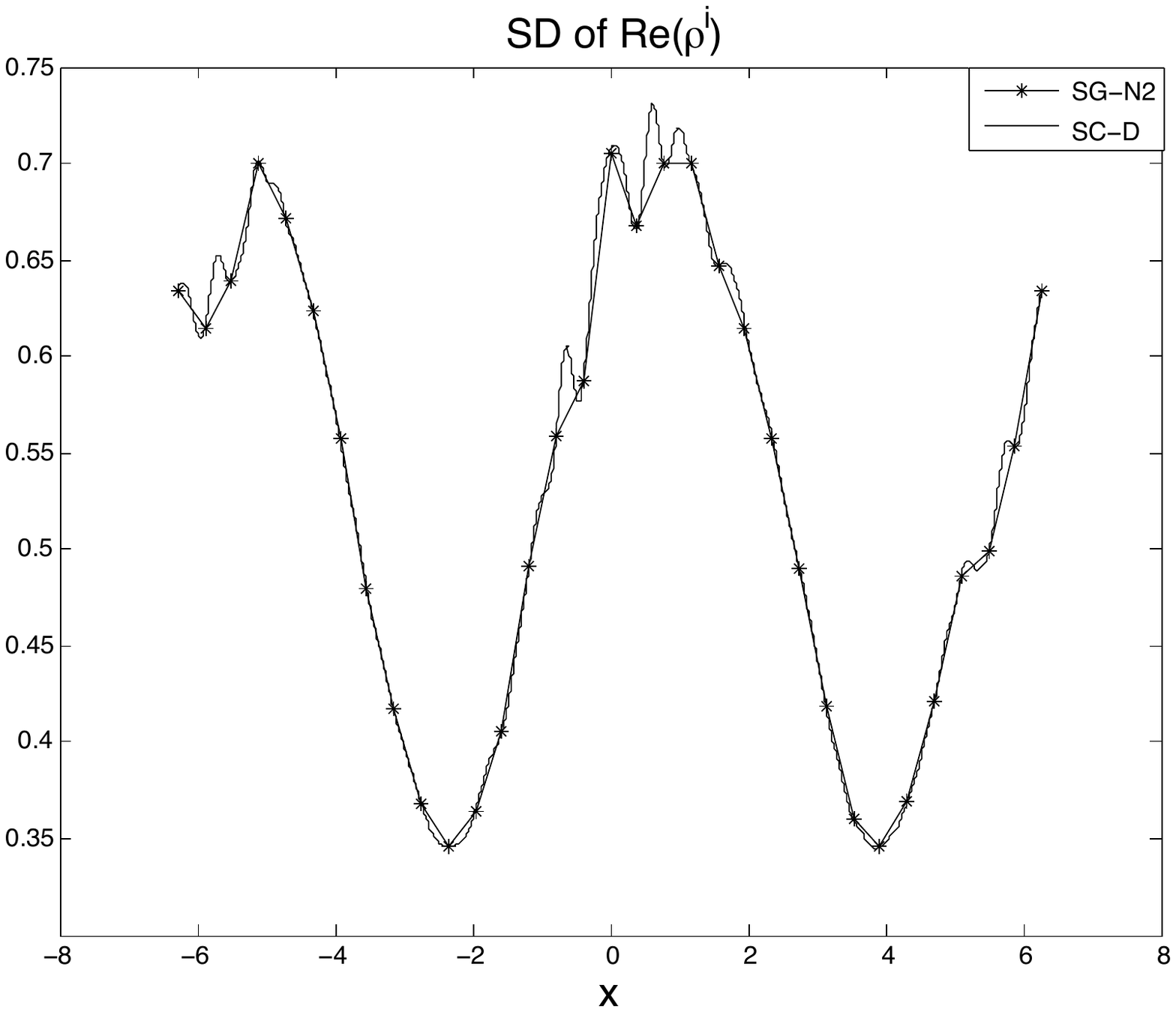}
  \end{subfigure}
 \begin{subfigure}{0.5\textwidth}
   \centering
 \includegraphics[trim={0 9cm 0 0},width=1.05\textwidth, height=0.8\textwidth]{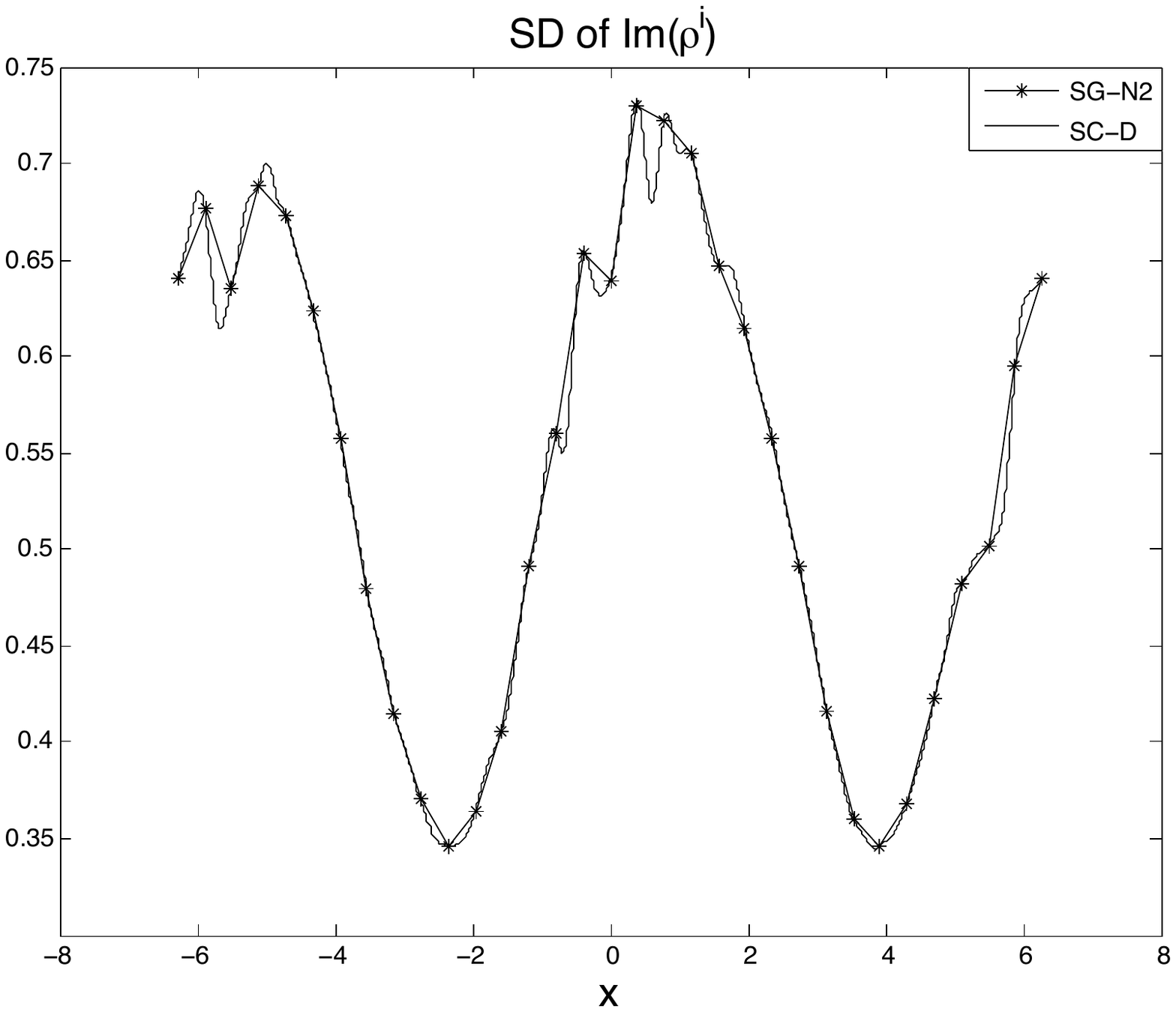}
  \end{subfigure}
\caption{Example $4.1$, with $\displaystyle E(x,\bz)=(10-\cos(x/2))(1+0.5\bz)$.
$\varepsilon=0.01$, $t=0.3$, $N_x=2000$, $\Delta t=5\times 10^{-3}$ (gPC-SC-D),
 and $N_x=32$, $\Delta t=0.02$ (gPC-SG-N2). Mean and standard deviation of the space dependence
of the densities $\rho^{\pm}$, $\text{Re}(\rho^{i})$ and $\text{Im}(\rho^{i})$.
Stars: gPC-SG-N2 with $K=4$. Solid lines: reference solution by gPC-SC-D.}
\label{SH_Ex2}
\end{figure}

Figure \ref{SH_Ex2} shows a similar result as in Figure \ref{gPC_N2}, using Example $4.1$ except that here $E(x,\bz)={\cal O}(1)$.
As discussed in Remark \ref{rmk},  small $\Delta t
=O(\sqrt{\varepsilon}\Delta x)$ needs to be chosen in Figure \ref{gPC_N2} due to the CFL condition in the
transport steps for $E=O(\sqrt{\varepsilon})$. However, larger $\Delta t=O(\Delta x)$, which is {\it independent of} $\varepsilon$, can be used in Figure \ref{SH_Ex2} since now $E(x,\bz)={\cal O}(1)$.

\begin{figure} [H]
\begin{subfigure}{0.5\textwidth}
   \centering
 \includegraphics[trim={0 9cm 0 0},width=1.0\textwidth, height=0.75\textwidth]{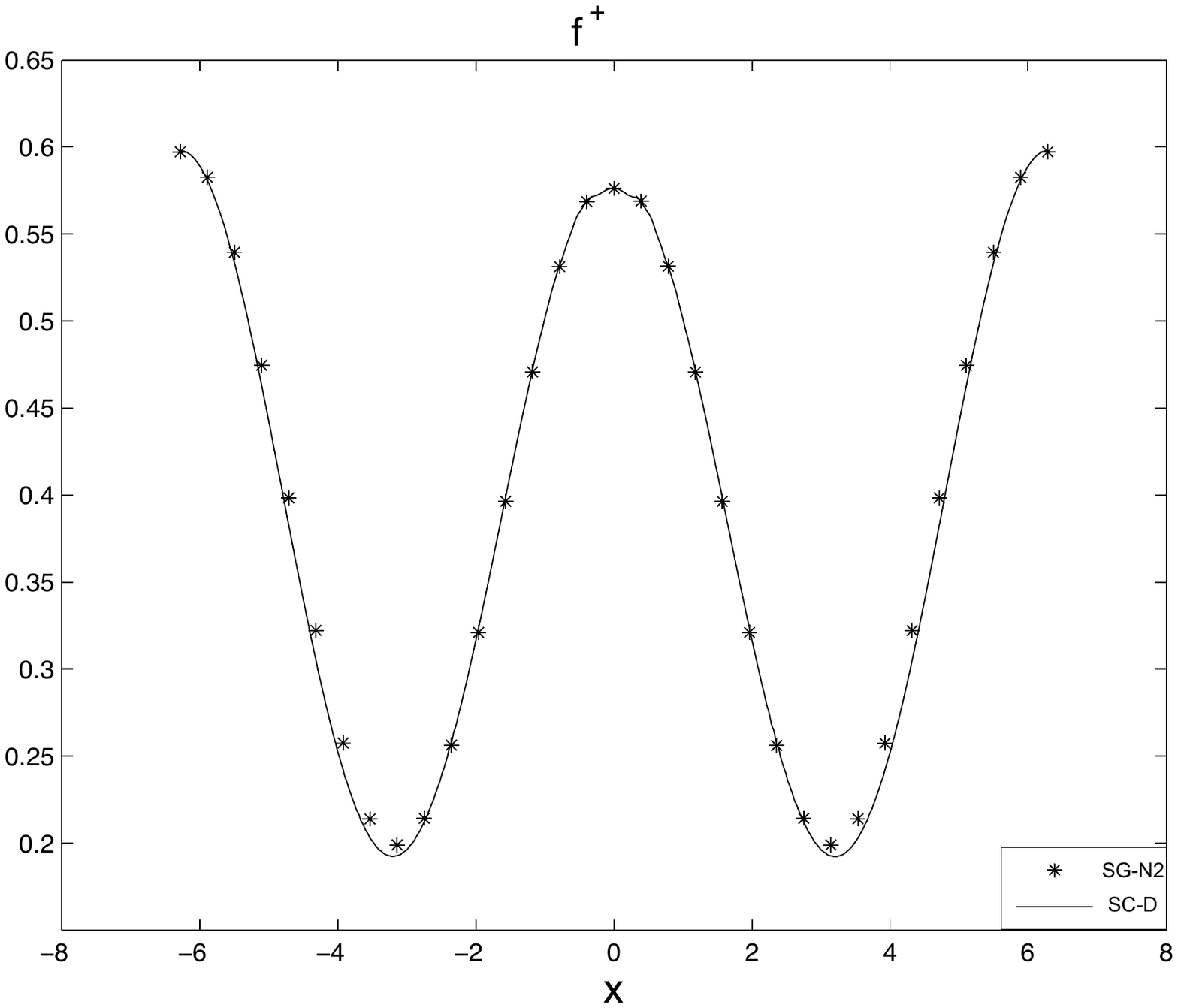}
  \end{subfigure}
  \begin{subfigure}{0.5\textwidth}
 \centering
 \includegraphics[trim={0 9cm 0 0},width=1.0\textwidth, height=0.75\textwidth]{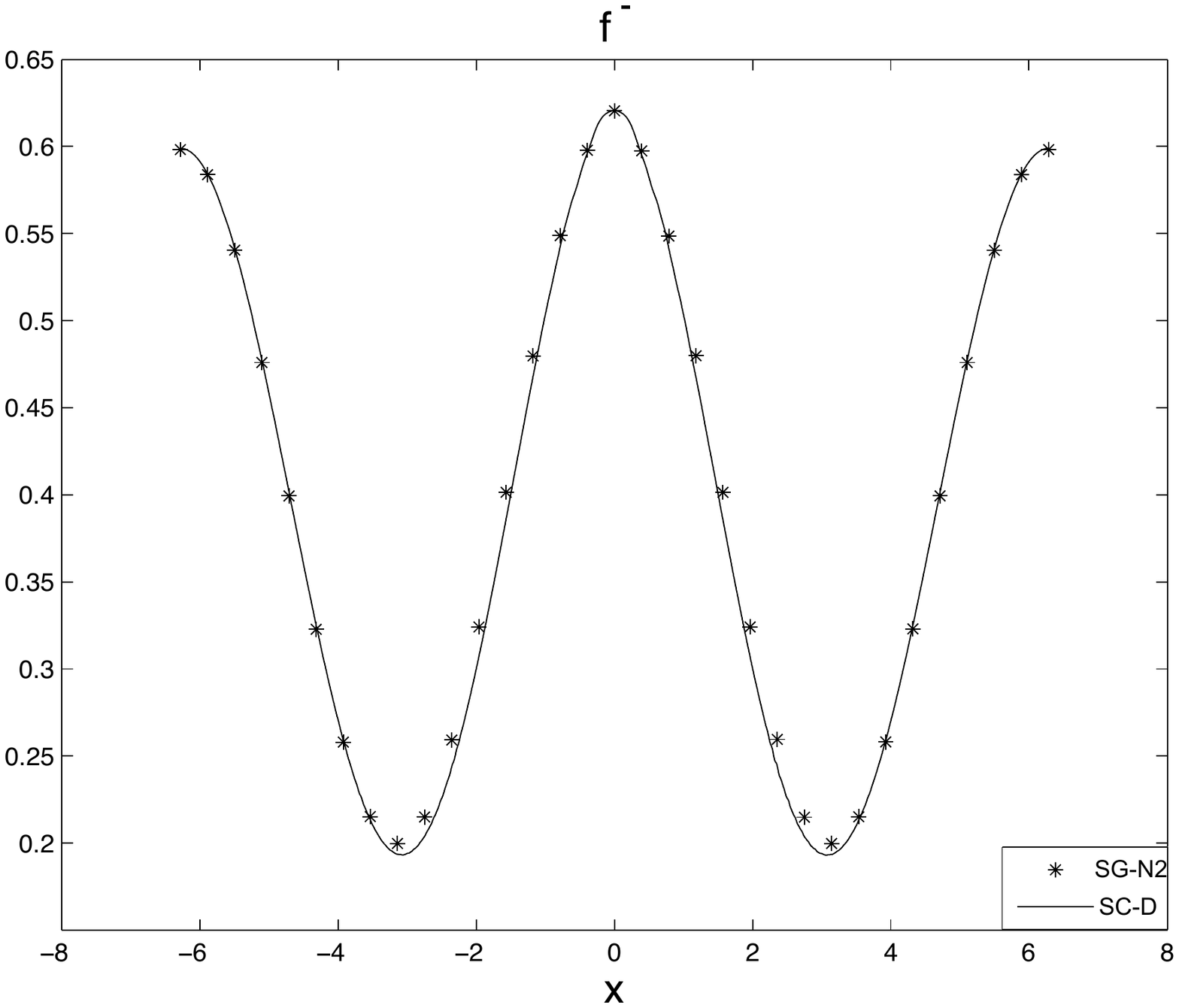}
  \end{subfigure}
  \begin{subfigure}{0.5\textwidth}
 \centering
 \includegraphics[trim={0 9cm 0 0},width=1.0\textwidth, height=0.8\textwidth]{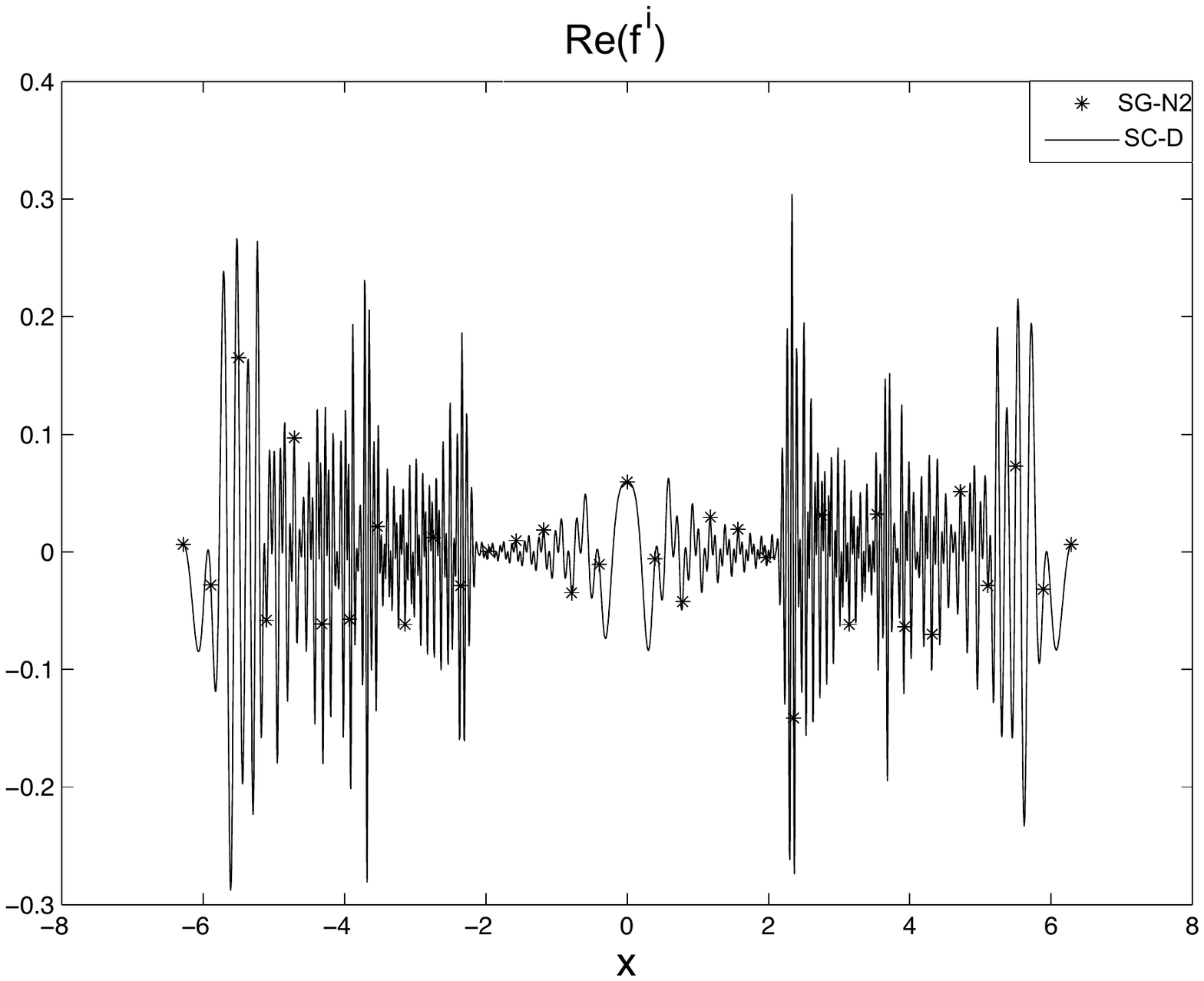}
  \end{subfigure}
  \begin{subfigure}{0.5\textwidth}
   \centering
 \includegraphics[trim={0 9cm 0 0},width=1.0\textwidth, height=0.8\textwidth]{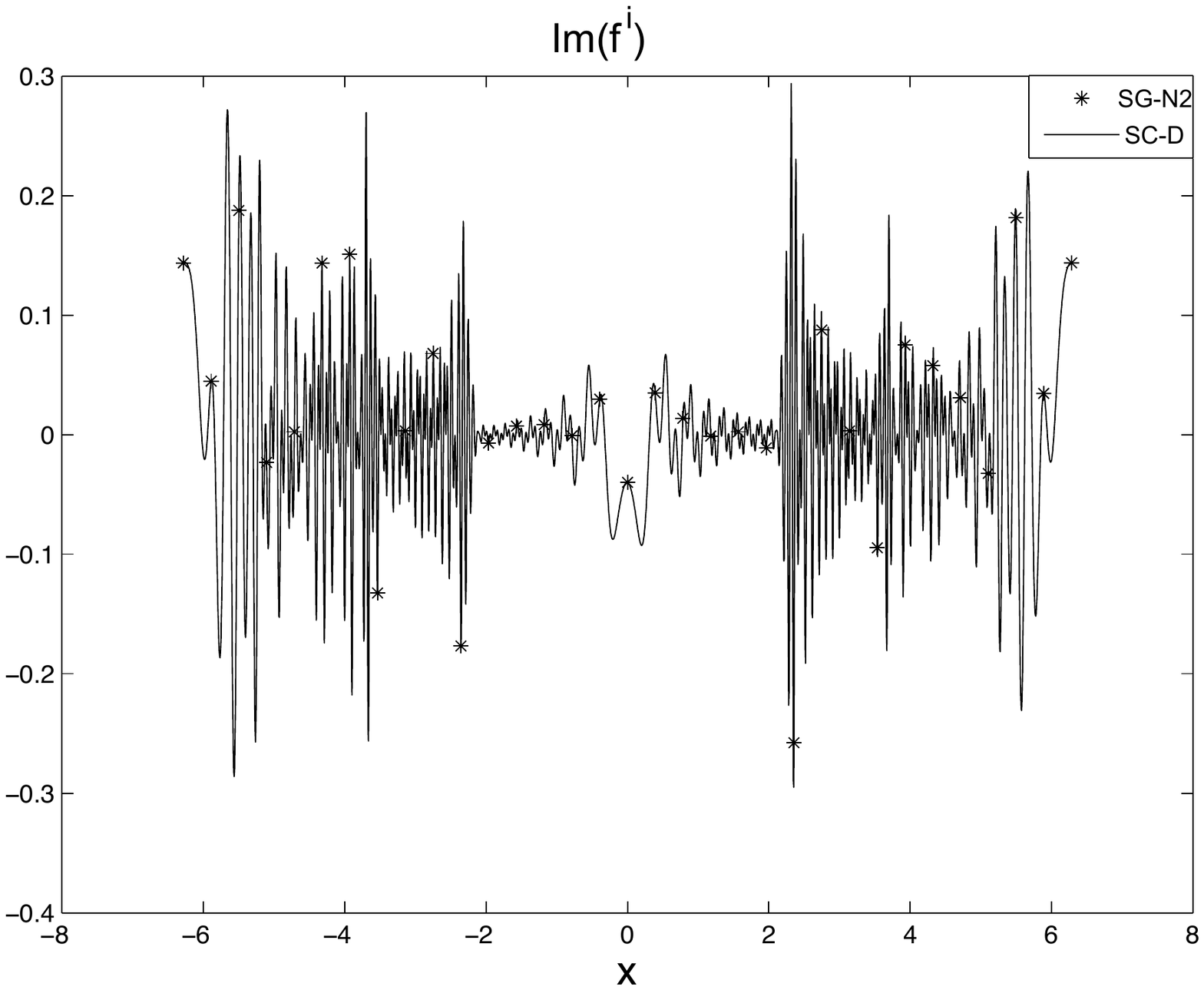}
  \end{subfigure}
   \begin{subfigure}{0.5\textwidth}
 \centering
 \includegraphics[trim={0 9cm 0 0},width=1.0\textwidth, height=0.75\textwidth]{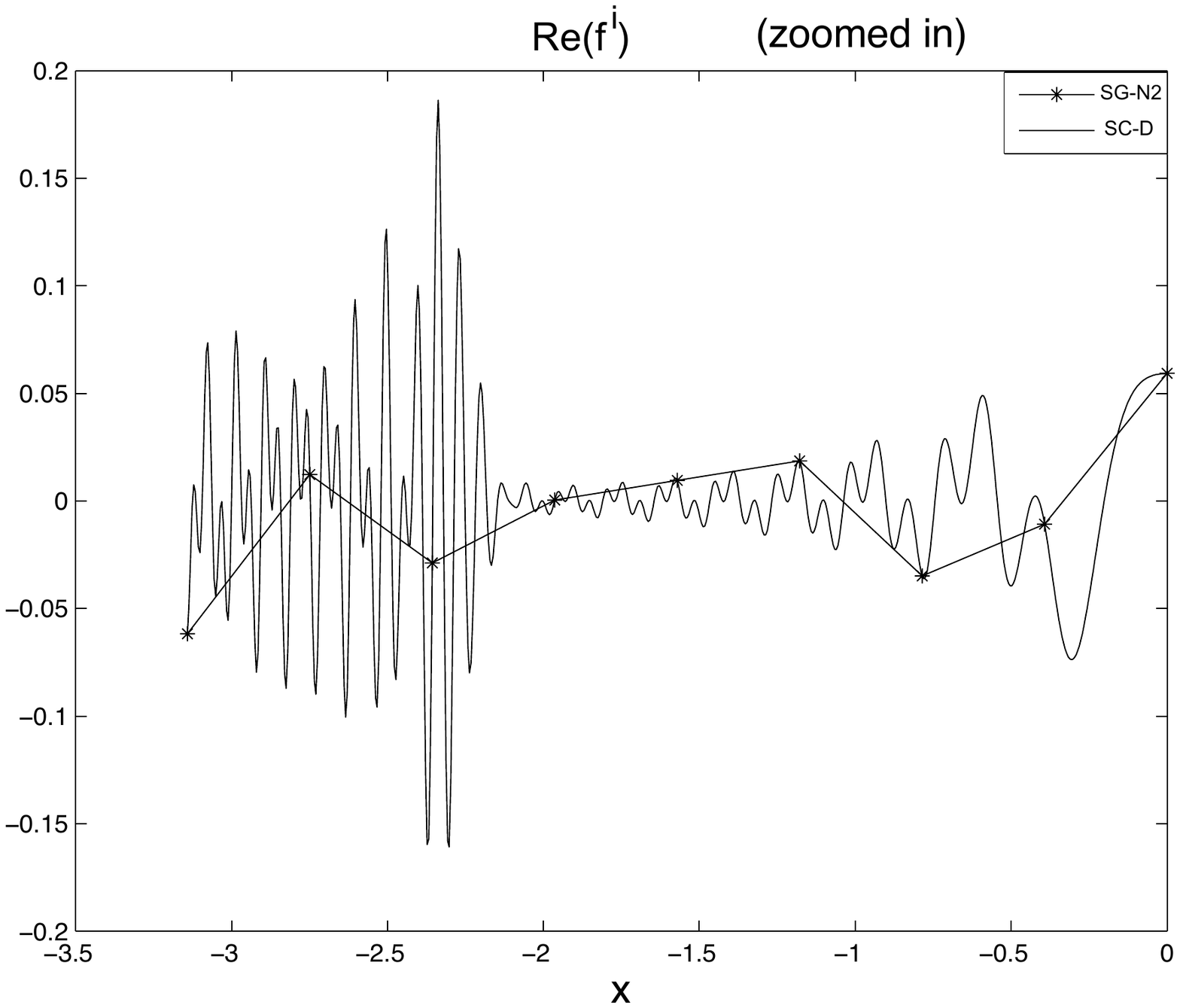}
  \end{subfigure}
  \begin{subfigure}{0.5\textwidth}
   \centering
 \includegraphics[trim={0 9cm 0 0},width=1.0\textwidth, height=0.75\textwidth]{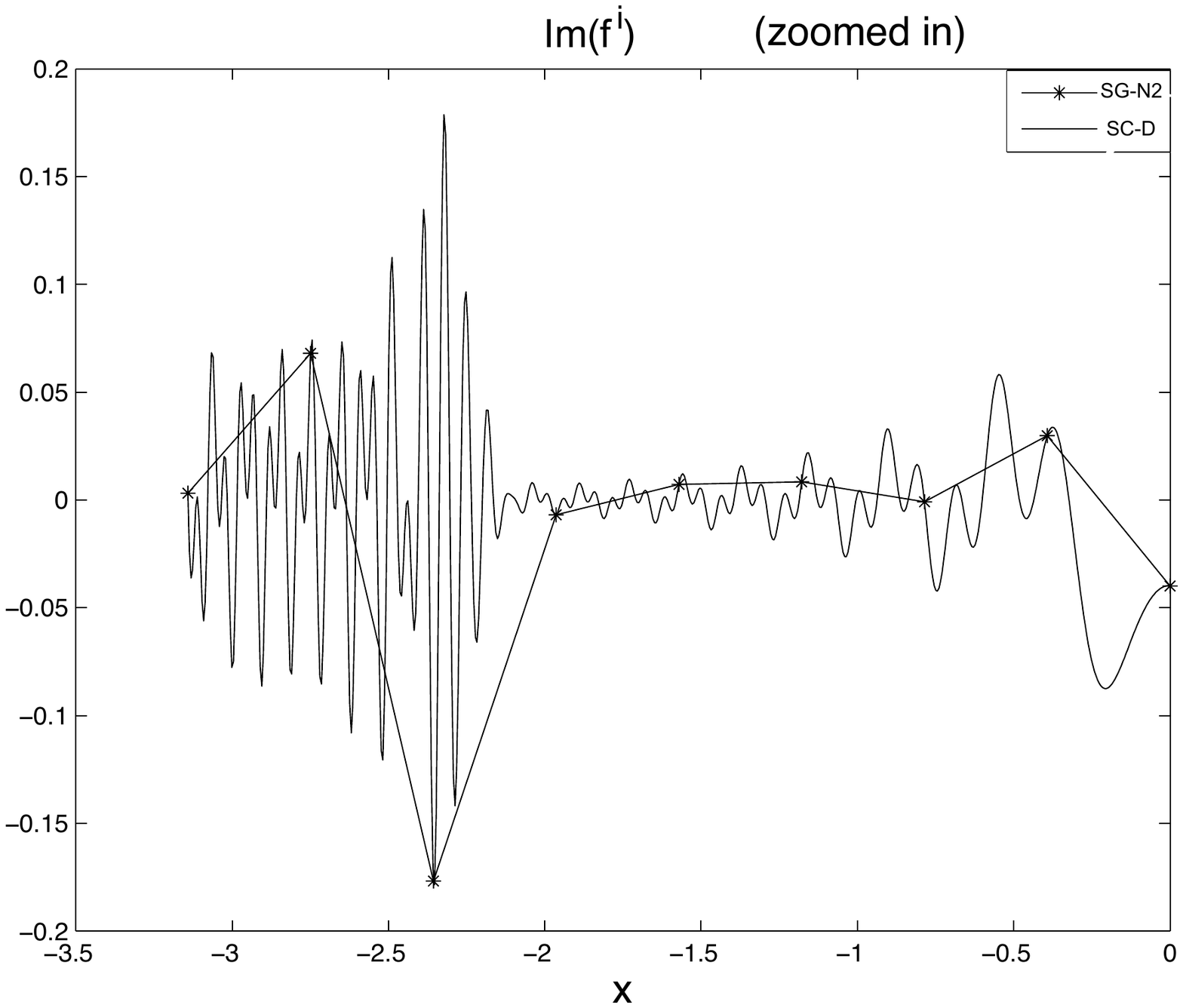}
  \end{subfigure}
  \makeatletter
   \renewcommand{\fnum@figure}{\figurename~\thefigure (a)}
  \makeatother
  \caption{Example $4.1$. $\varepsilon=5\times 10^{-3}$, $t=0.5$, $N_x=2500$, $\Delta t=5\times 10^{-3}$ (gPC-SC-D),
 and $N_x=32$, $\Delta t=5\times 10^{-5}$ (gPC-SG-N2). Mean of the space dependence
 of $f^{\pm}$, $\text{Re}(f^{i})$ and $\text{Im}(f^{i})$ (and their zoomed in solutions) at $p=0$.
 Stars: gPC-SG-N2 with $K=4$. Solid lines: reference solution by gPC-SC-D.}
  \end{figure}
 \begin{figure}[H]\ContinuedFloat
  \begin{subfigure}{0.5\textwidth}
   \centering
 \includegraphics[trim={0 9cm 0 0},width=1.0\textwidth, height=0.7\textwidth]{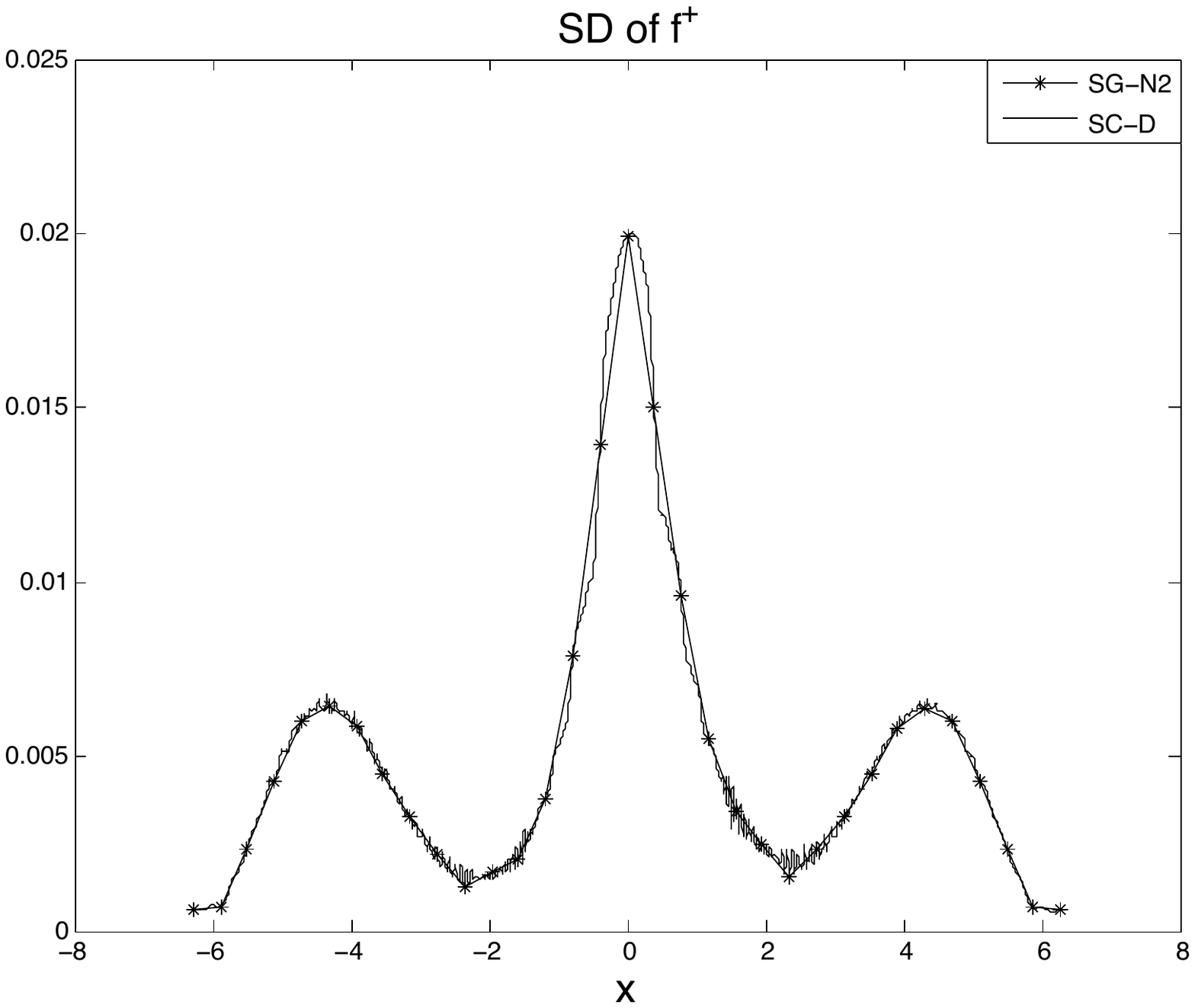}
  \end{subfigure}
 \begin{subfigure}{0.5\textwidth}
   \centering
 \includegraphics[trim={0 9cm 0 0},width=1.0\textwidth, height=0.7\textwidth]{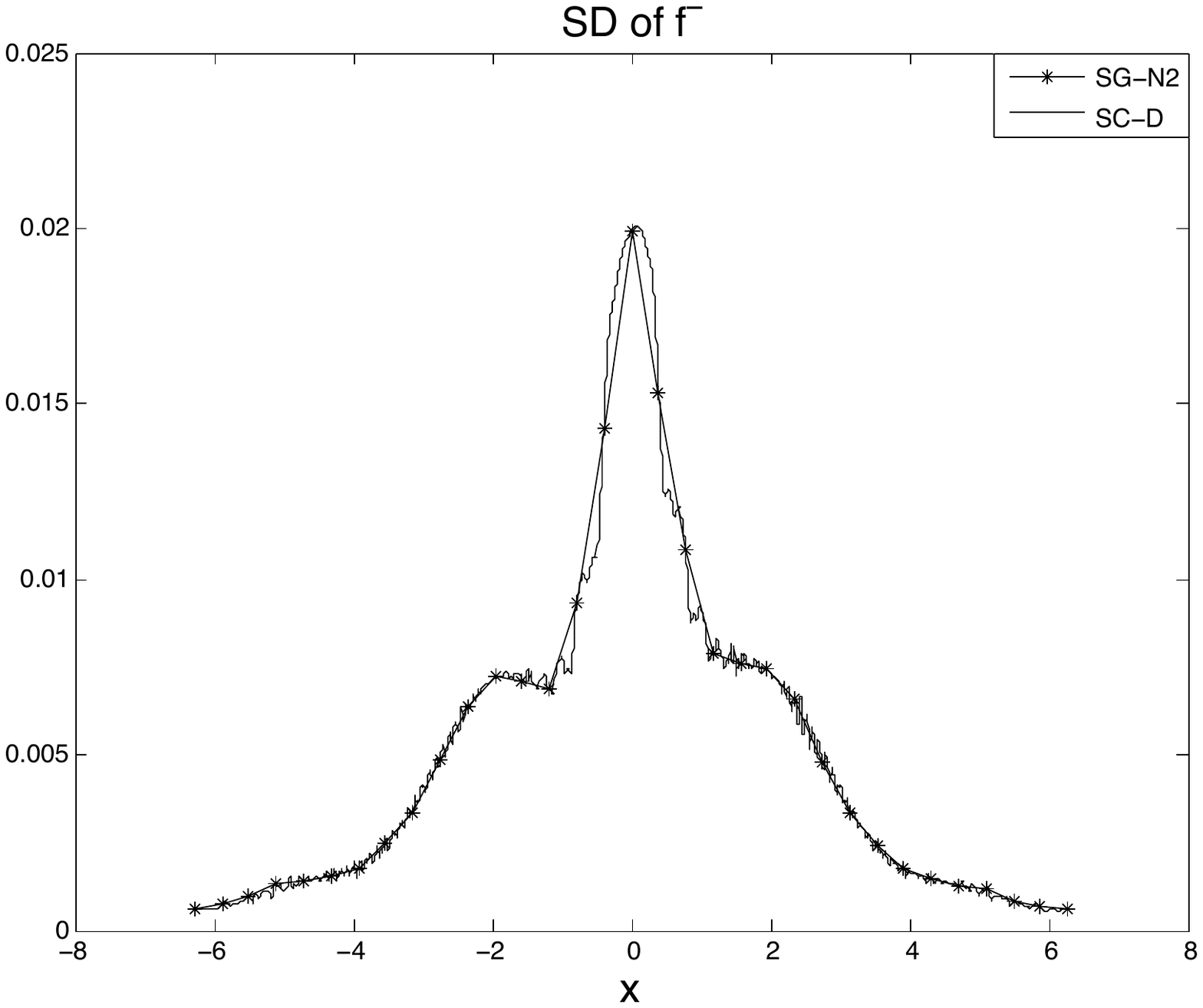}
  \end{subfigure}
 \begin{subfigure}{0.5\textwidth}
   \centering
 \includegraphics[trim={0 9cm 0 0},width=1.0\textwidth, height=0.7\textwidth]{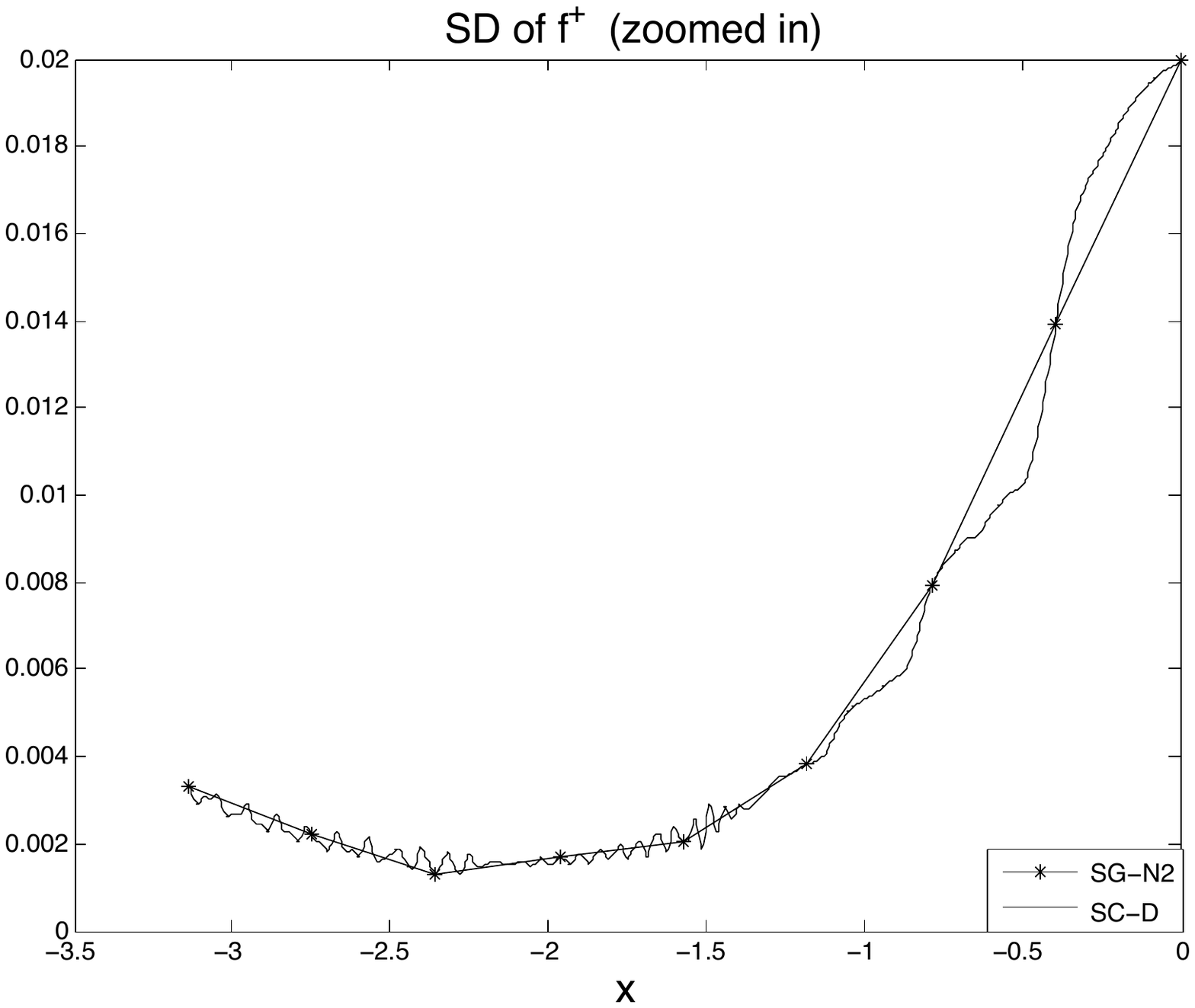}
  \end{subfigure}
  \begin{subfigure}{0.5\textwidth}
   \centering
 \includegraphics[trim={0 9cm 0 0},width=1.0\textwidth, height=0.7\textwidth]{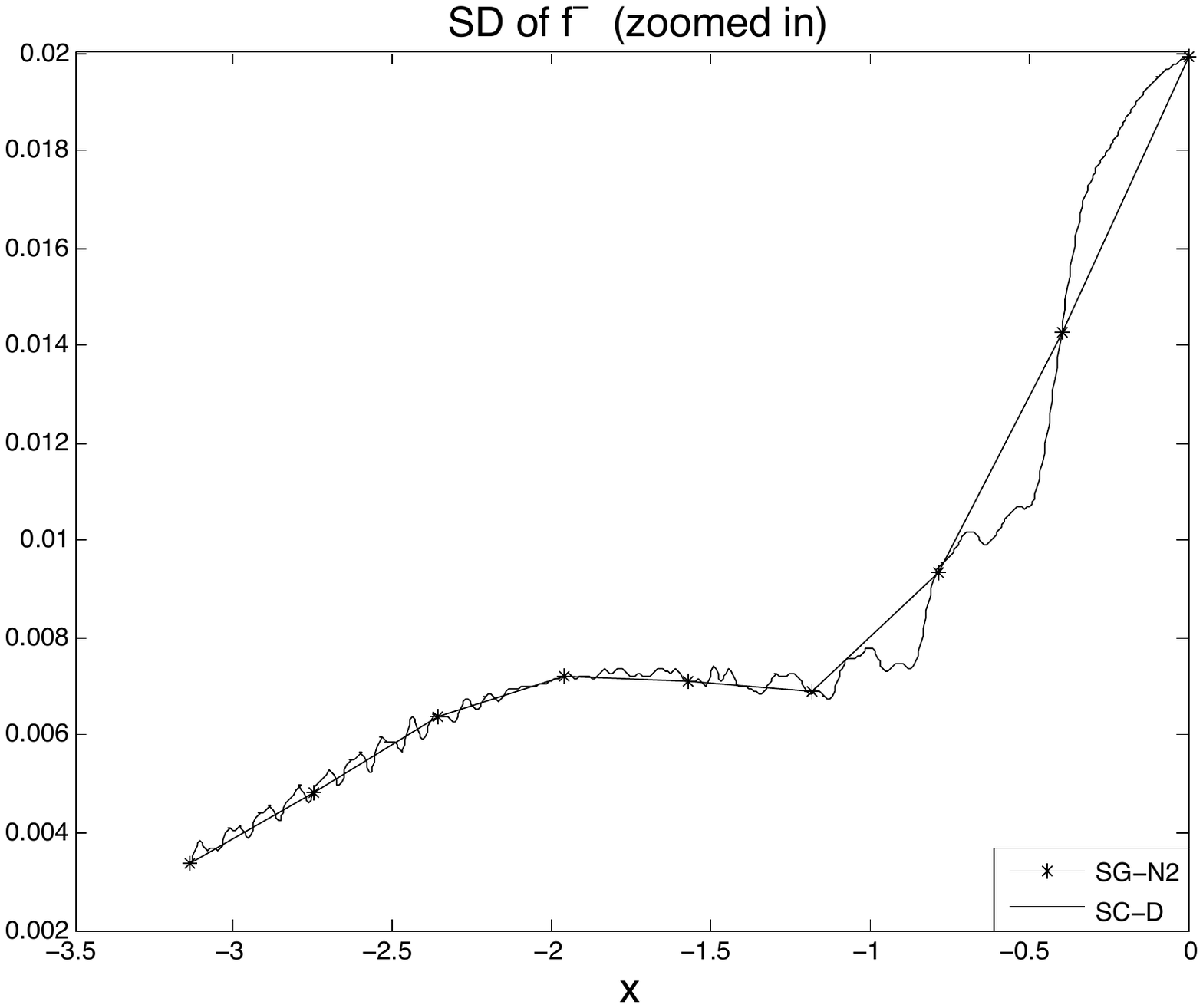}
  \end{subfigure}
 \begin{subfigure}{0.5\textwidth}
   \centering
 \includegraphics[trim={0 9cm 0 0},width=1.0\textwidth, height=0.7\textwidth]{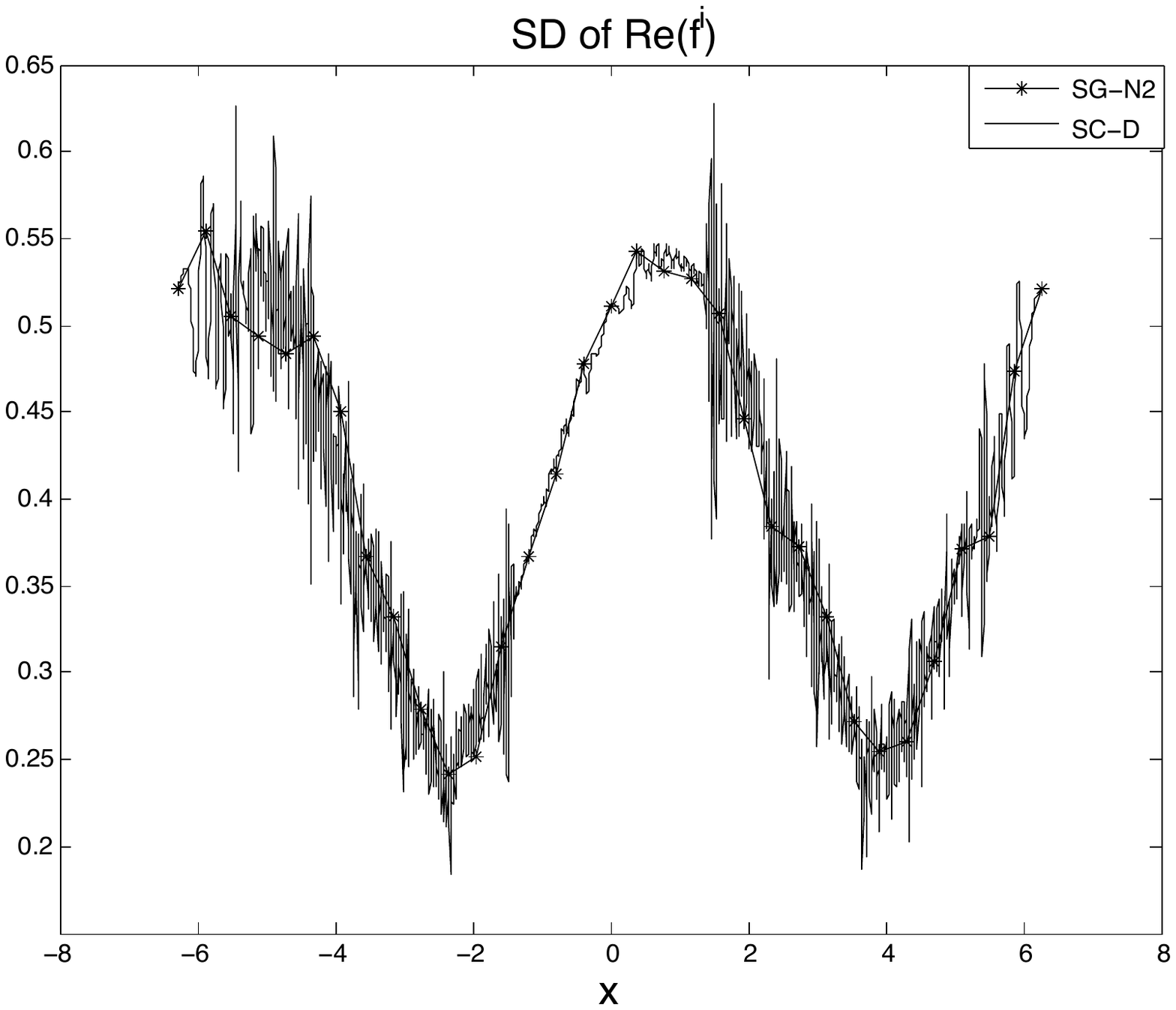}
  \end{subfigure}
 \begin{subfigure}{0.5\textwidth}
   \centering
 \includegraphics[trim={0 9cm 0 0},width=1.0\textwidth, height=0.7\textwidth]{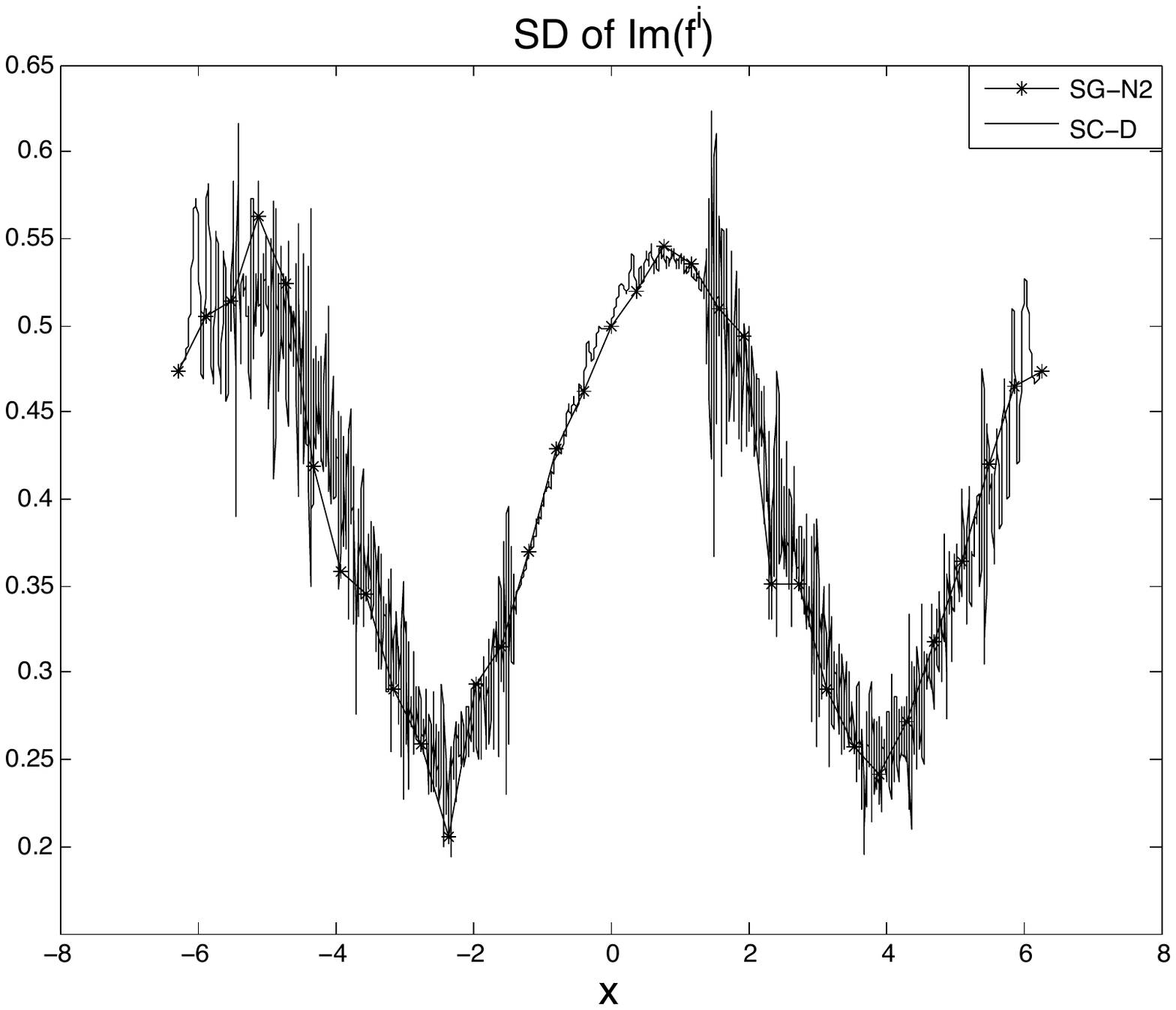}
  \end{subfigure}
  \begin{subfigure}{0.5\textwidth}
   \centering
 \includegraphics[trim={0 9cm 0 0},width=1.0\textwidth, height=0.7\textwidth]{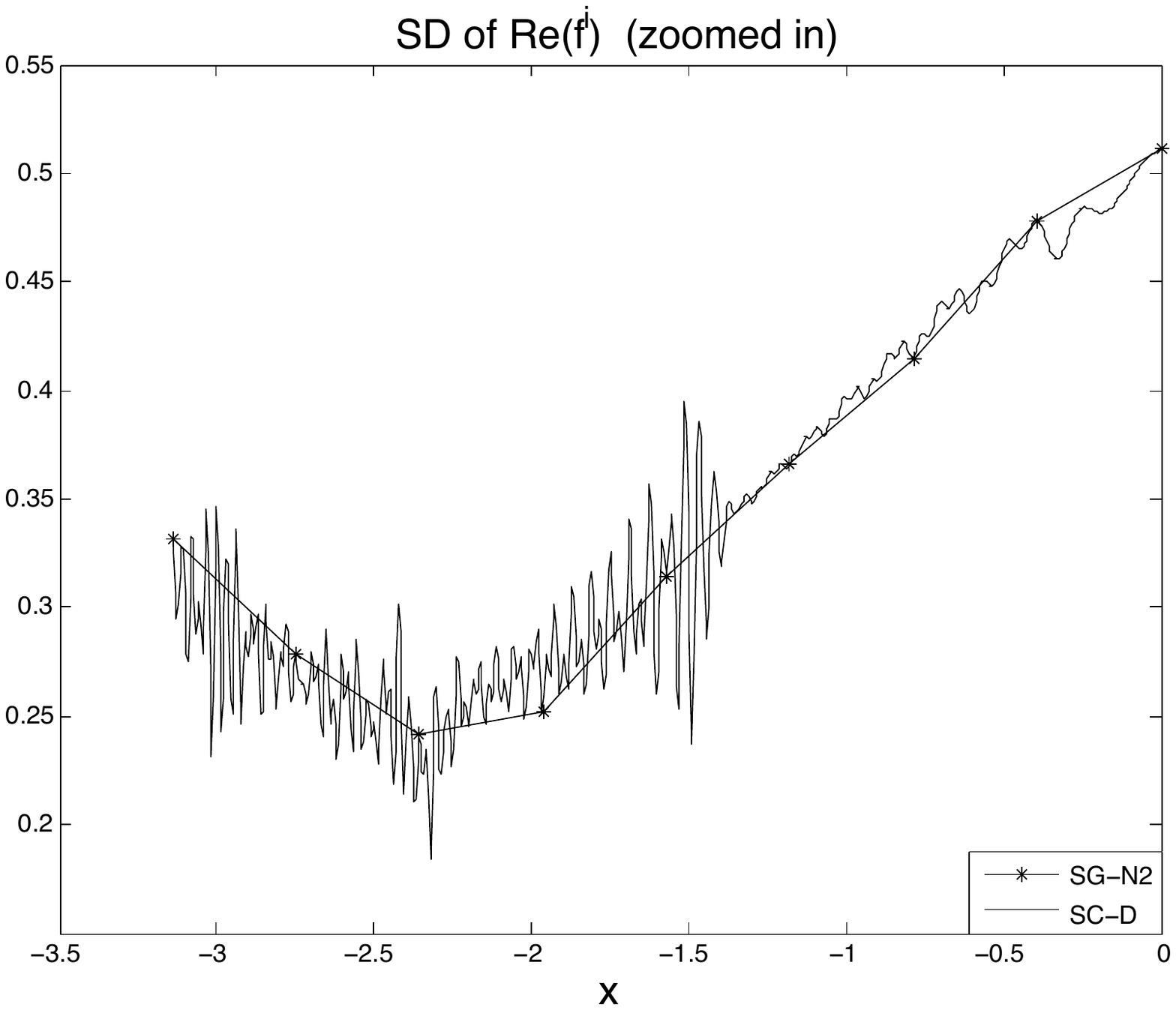}
  \end{subfigure}
 \begin{subfigure}{0.5\textwidth}
   \centering
 \includegraphics[trim={0 9cm 0 0},width=1.0\textwidth, height=0.7\textwidth]{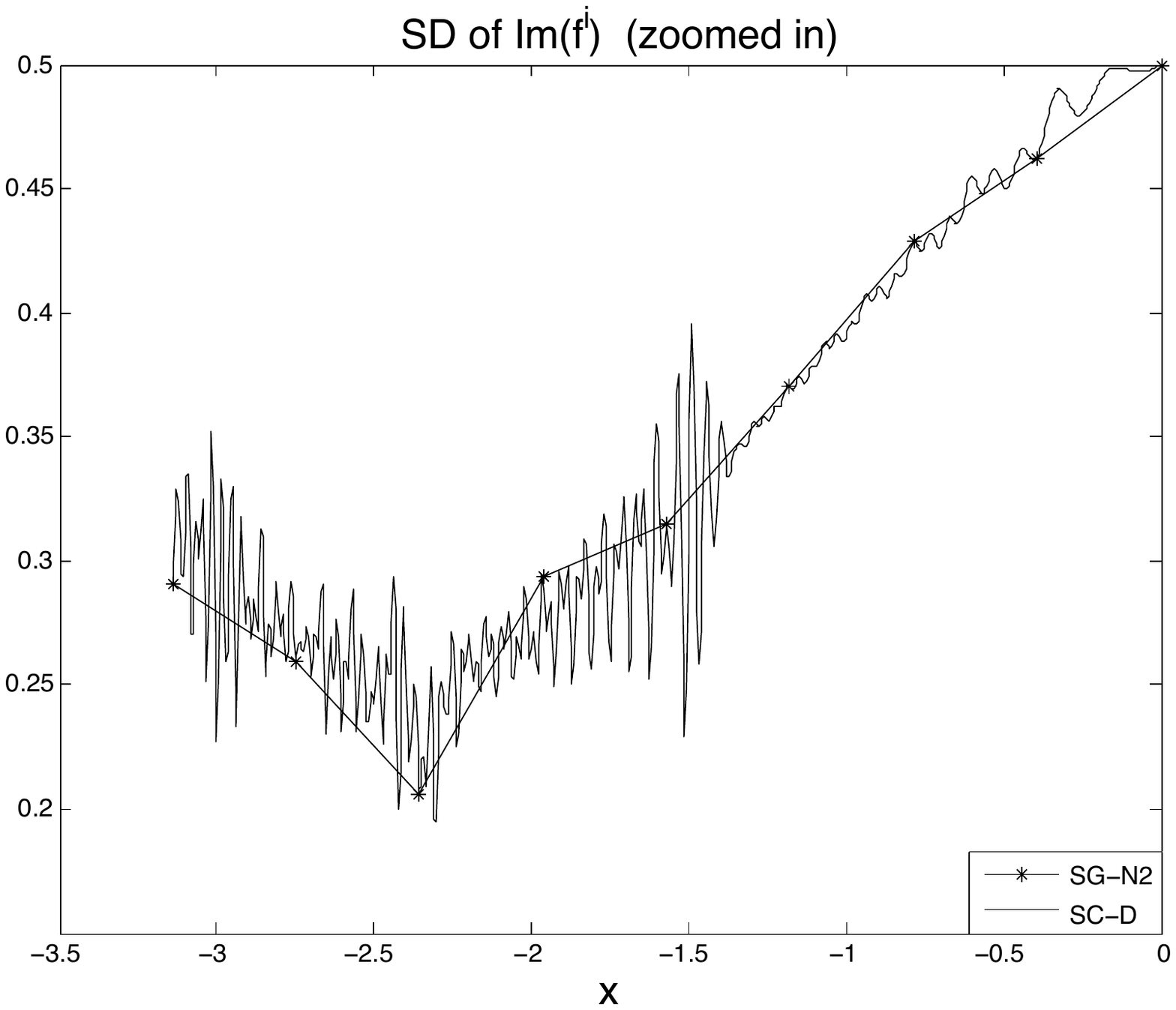}
  \end{subfigure}
  \makeatletter
   \renewcommand{\fnum@figure}{\figurename~\thefigure (b)}
  \makeatother
 \caption{Standard deviation of the space dependence
 of $f^{\pm}$, $\text{Re}(f^{i})$ and $\text{Im}(f^{i})$ (and their zoomed in solutions) at $p=0$.}
\end{figure}
 \begin{figure}[H]\ContinuedFloat
 \begin{subfigure}{0.5\textwidth}
   \centering
 \includegraphics[trim={0 9cm 0 0},width=1.05\textwidth, height=0.8\textwidth]{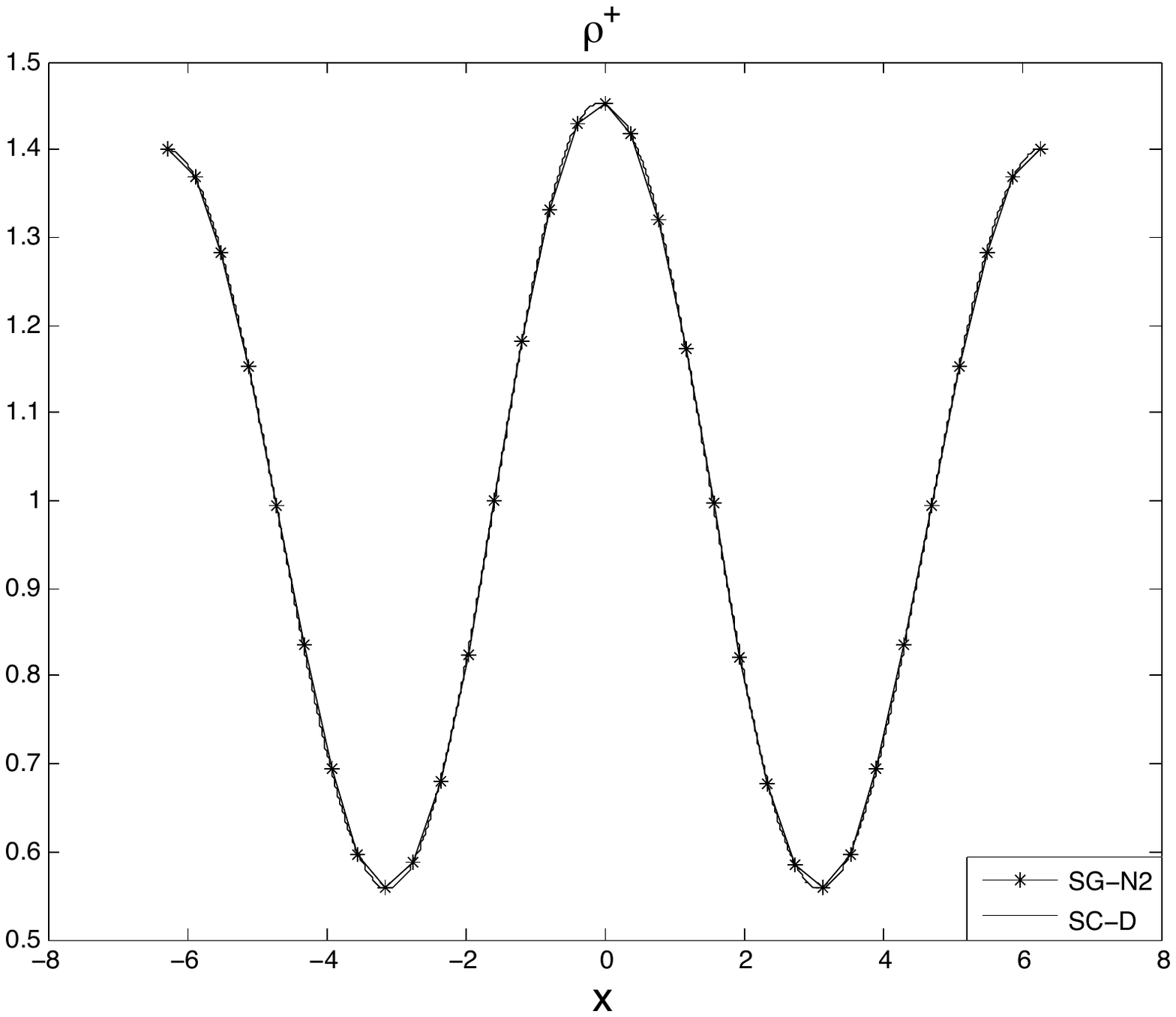}
  \end{subfigure}
   \begin{subfigure}{0.5\textwidth}
   \centering
 \includegraphics[trim={0 9cm 0 0},width=1.05\textwidth, height=0.8\textwidth]{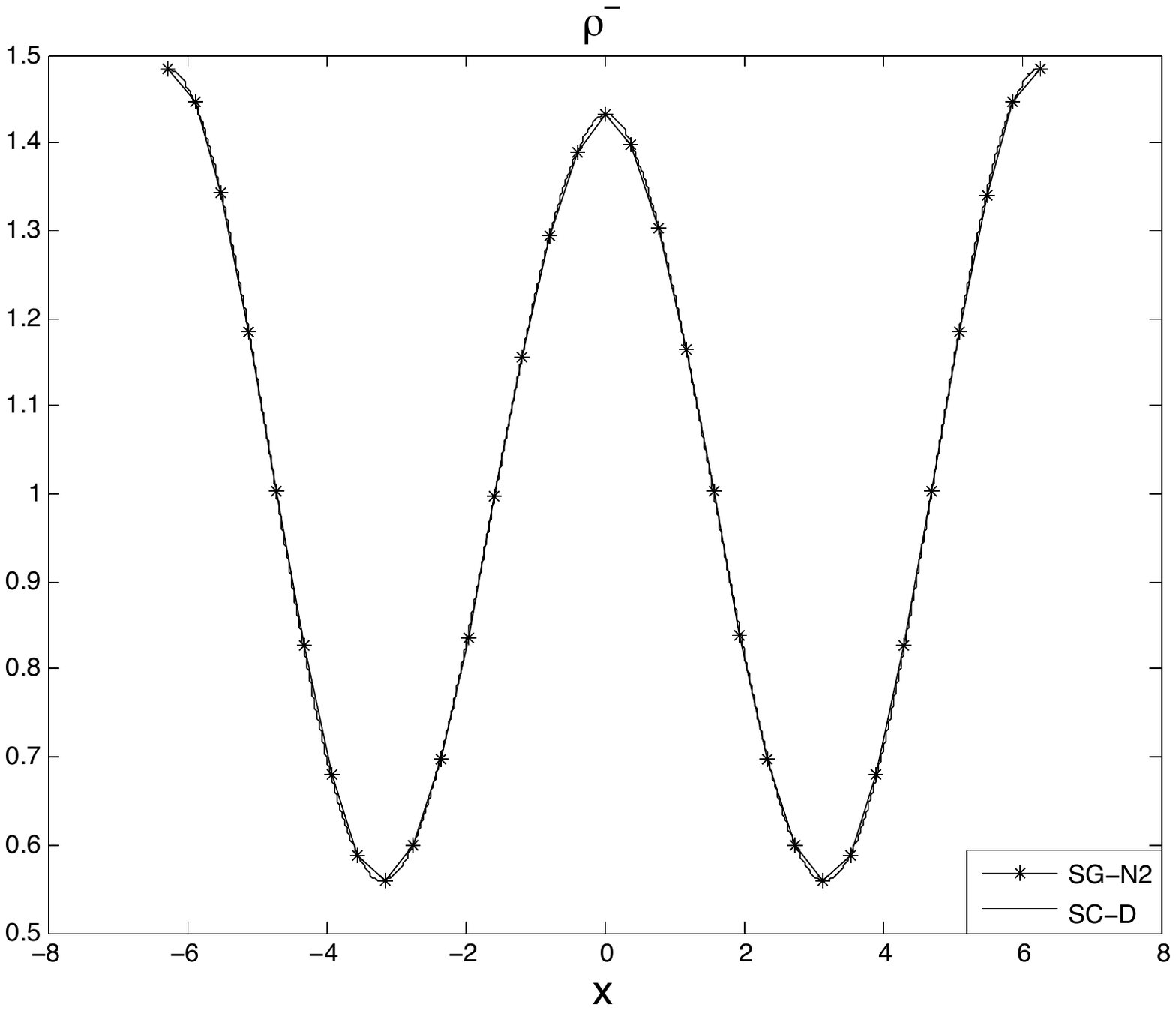}
  \end{subfigure}
\begin{subfigure}{0.5\textwidth}
   \centering
 \includegraphics[trim={0 9cm 0 0},width=1.05\textwidth, height=0.8\textwidth]{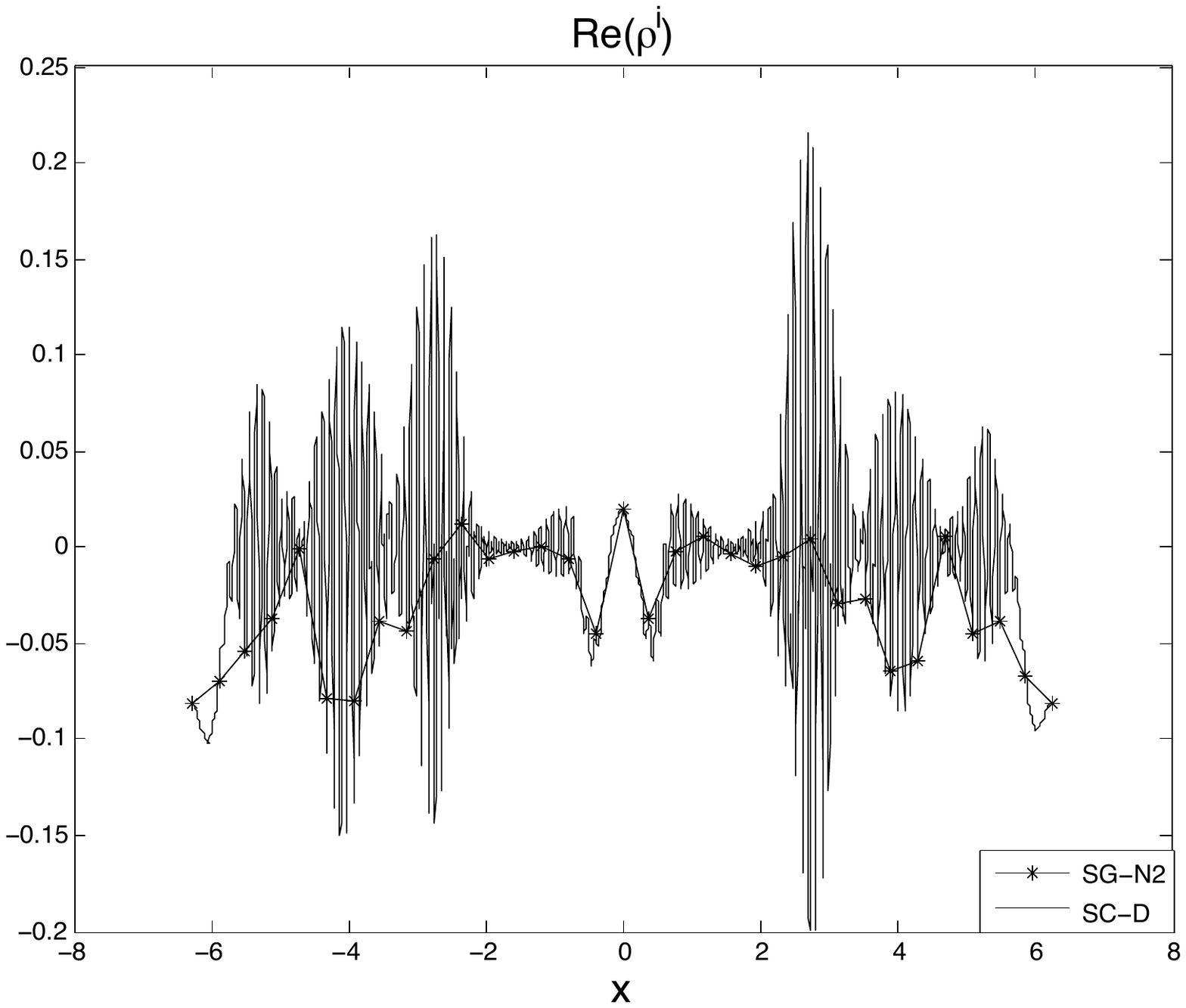}
  \end{subfigure}
 \begin{subfigure}{0.5\textwidth}
   \centering
 \includegraphics[trim={0 9cm 0 0},width=1.05\textwidth, height=0.8\textwidth]{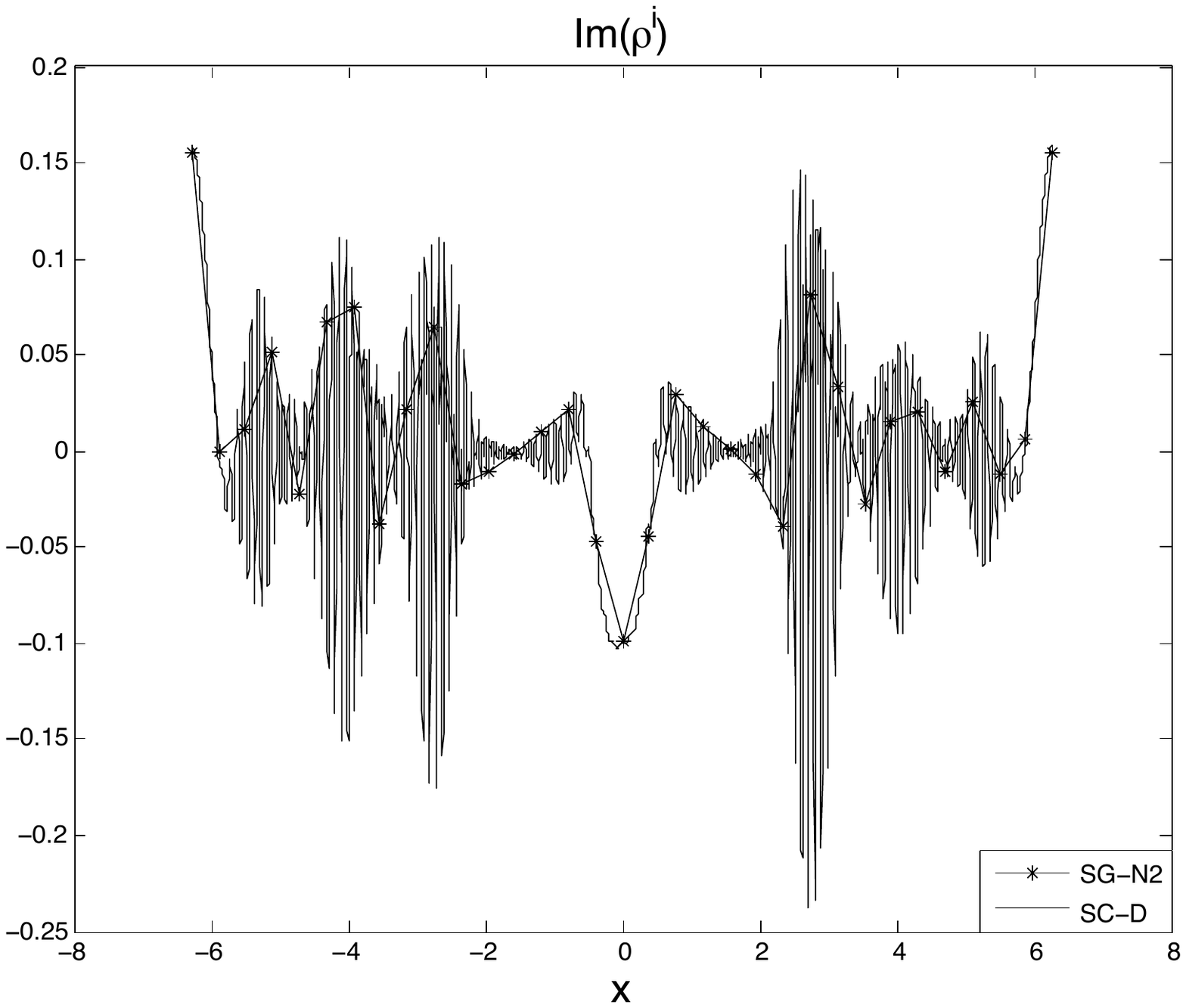}
  \end{subfigure}
  \begin{subfigure}{0.5\textwidth}
   \centering
 \includegraphics[trim={0 9cm 0 0},width=1.05\textwidth, height=0.8\textwidth]{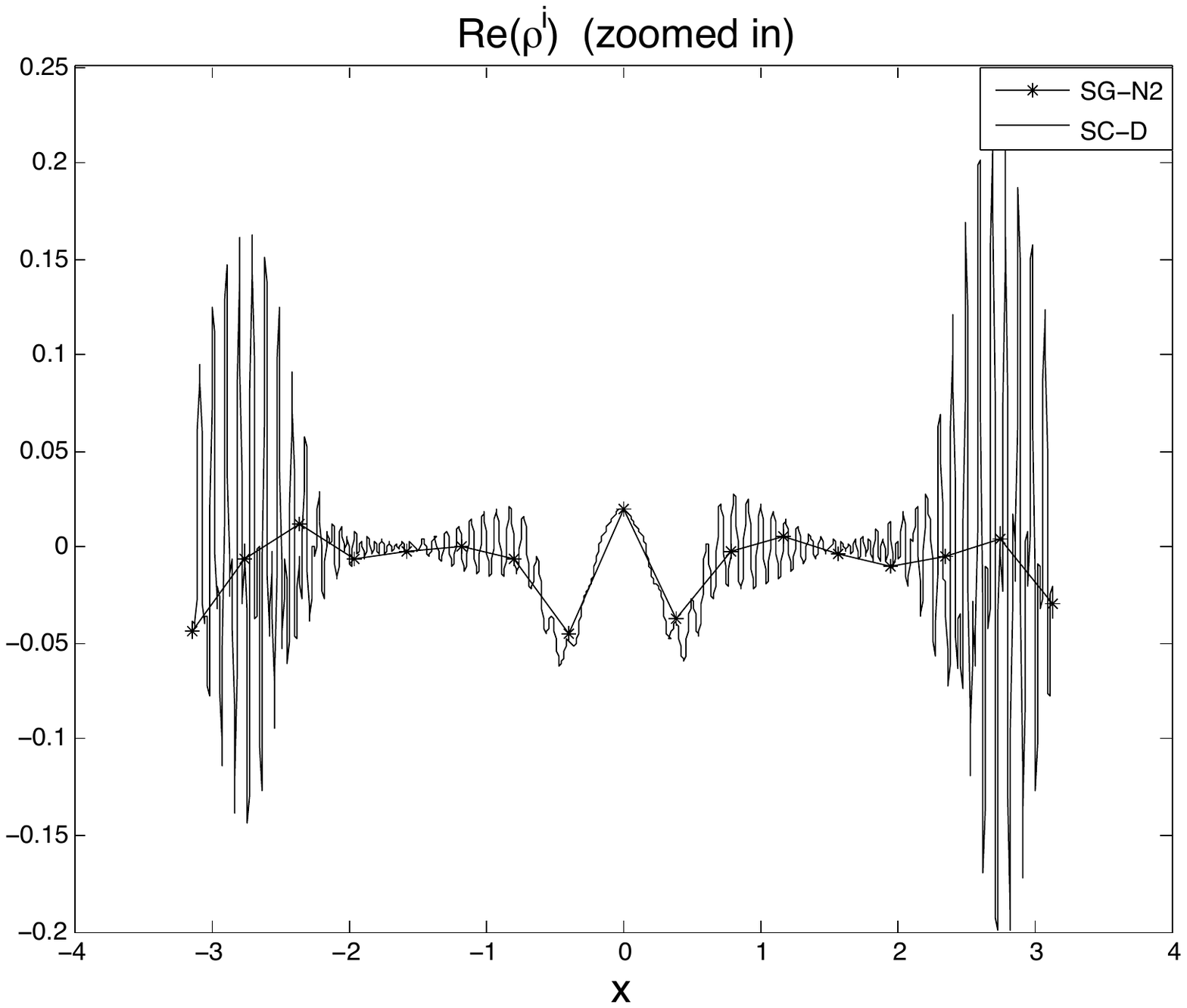}
  \end{subfigure}
 \begin{subfigure}{0.5\textwidth}
   \centering
 \includegraphics[trim={0 9cm 0 0},width=1.05\textwidth, height=0.8\textwidth]{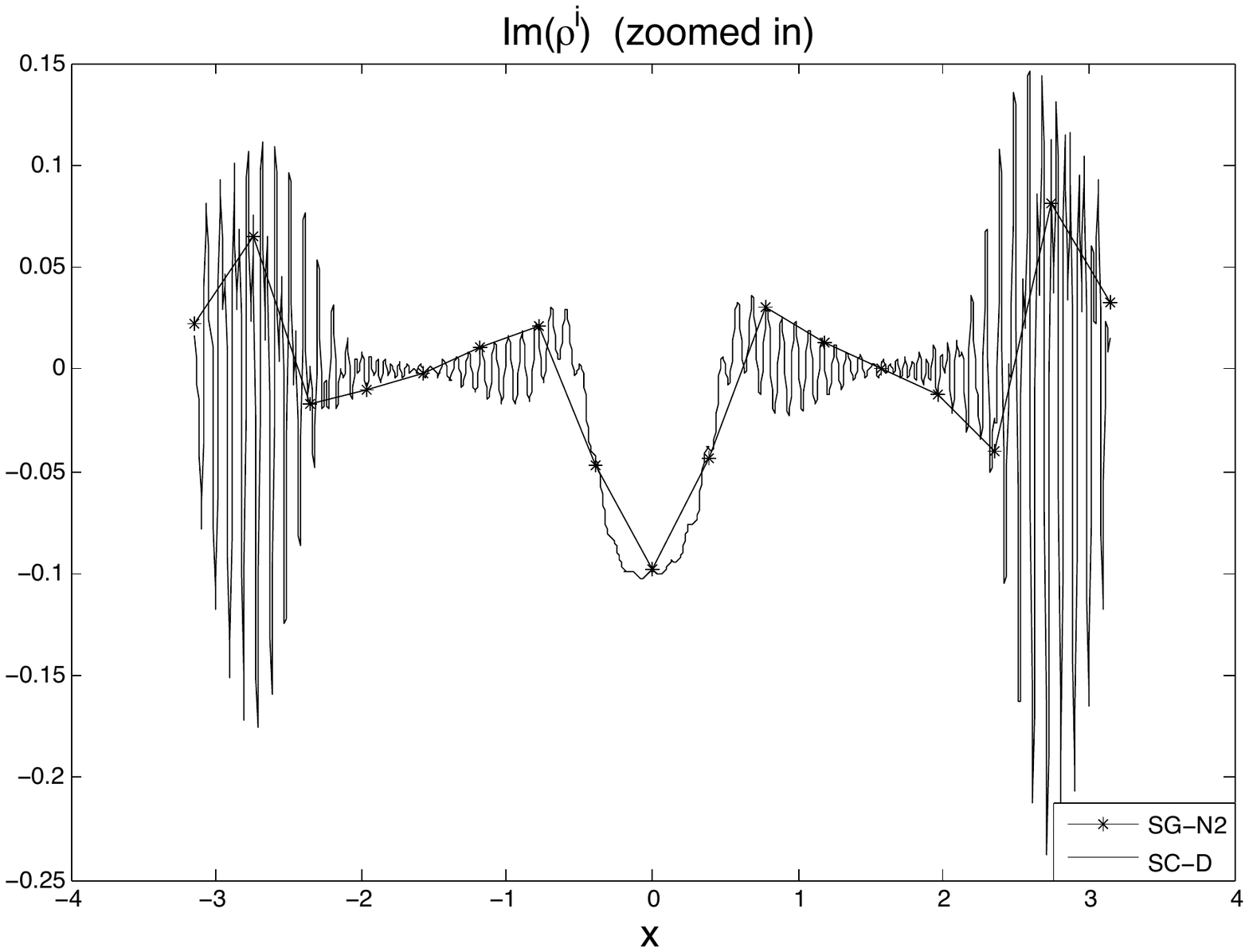}
  \end{subfigure}
  \makeatletter
   \renewcommand{\fnum@figure}{\figurename~\thefigure (c)}
  \makeatother
\caption{Mean of the space dependence of the densities $\rho^{\pm}$, $\text{Re}(\rho^{i})$ and $\text{Im}(\rho^{i})$ (and their zoomed in solutions). }
 \label{gPC_N2}
 \end{figure}

\section {Conclusions}
In this paper, based on a nonlinear geometric optics (NGO) based numerical
method developed in \cite{NGO}, with a new ``time'' variable defined from the phase,
we obtain a stochastic Galerkin (SG) method
for highly oscillatory transport equations that arise in semiclassical models
of non-adiabatic quantum dynamics, in which the potential energy surfaces are assumed to be random, due to uncertainties in modeling or measurement errors.
We prove that the generalized polynomial chaos (gPC) based SG method applied directly to
the models will require the order of gPC to depend on the possibly very small
wave length, while the new method does not have such a requirement.
This important property allows us to use this method to solve these
highly oscillatory problems with uncertain coefficients with all numerical parameters independent of the wave length, yet still capture the solution statistics pointwisely.

There are several projects along this direction. First it will be desirable to develop an implicit scheme for the transport steps in order to
obtain an improved time step constraint when the bad gap becomes very small. Second, methods for higher dimensional--in space, velocity as well as
the random variables--remain to be developed.

\appendix\newpage\markboth{Appendix}{Appendix}
\renewcommand{\thesection}{\Alph{section}}
\numberwithin{equation}{section}
\section{Appendix: The gPC-SG-N1 and gPC-SG-D for the surface hopping model with random inputs}

\subsubsection*{A.1. The gPC Approximation}
In this part, we detail the gPC-SG-N1 strategy for model \eqref{SH old}, based on the method developed in \cite{NGO} for the deterministic problem.
Solving for $S$ follows the same strategy presented in section \ref{N2_3D} on the gPC-SG-N2 method.

Now we focus on the approximation of the profiles $(F^+, F^-, G, H)$. One inserts the Galerkin approximation of  each component  $$
(F^+, F^-, G, H)(t,x,p,\bz)=\sum_{j=1}^{K} (F^+_j, F^-_j, G_j, H_j)(t,x,p)\tilde \psi_{j}(\bz),
$$
into (\ref{SH new}) and conduct the gPC Galerkin projection.
Denote the associated gPC coefficients vector of $(\vv F^+, \vv F^-, \vv G, \vv H)$ as previously, for instance for $\vv F^+$,
\begin{equation}\vv F^{+}=(F^{+}_1, \cdots, F^{+}_K)^T. \label{gPC_vector_app} \end{equation}
Then, we get the gPC system
\begin{align}
\label{SH gPC}
\left\{
\begin{array}{l}
\displaystyle
\partial_t \vv F^{+}+p\partial_{x}\vv F^{+}-{\cal C}\partial_{p}\vv F^{+}=\frac{1}{\varepsilon}(-2{\cal C}+{\cal H})\partial_{\tau}\vv F^{+}
+2b^i( \vv G\cos\tau+\vv H\sin\tau), \\[8pt]
\displaystyle
\partial_t \vv F^{-}+p\partial_{x}\vv F^{-}+{\cal C}\partial_{p}\vv F^{-}=\frac{1}{\varepsilon}(-2{\cal C}-{\cal H})\partial_{\tau}\vv F^{-}
-2b^i(\vv G\cos\tau+ \vv H\sin\tau), \\[8pt]
\displaystyle
\partial_t \vv G+p\partial_{x}\vv G=-\frac{2{\cal C}}{\varepsilon}\partial_{\tau}\vv G-b^i\vv F^{+}\cos\tau+b^i\vv F^{-}\cos\tau, \\[8pt]
\displaystyle
\partial_t \vv H+p\partial_{x}\vv H=-\frac{2{\cal C}}{\varepsilon}\partial_{\tau}\vv H-b^i\vv F^{+}\sin\tau+b^i\vv F^{-}\sin\tau,
\end{array}\right.
\end{align}
where matrices ${\cal C}\in\mathbb R_{K\times K}$ and ${\cal H}\in\mathbb R_{K\times K}$ are defined by
\begin{align*}
& {\cal C}_{mn}(x)=\int_{I_{\bz}} E(x,\bz)\tilde \psi_{m}(\bz)\tilde \psi_{n}(\bz)d\bz, \qquad 1\leq m, n\leq K. \\[4pt]
& {\cal H}_{k m}(t,x,p)=\sum_{l=1}^K\partial_p\tilde S_{l}(t,x,p)D_{l m k}(x),  \qquad 1\leq k, m\leq K,
\end{align*}
with $\partial_p\tilde S_{l}(t,x,p)$ ($l=1, \cdots, K$) the gPC coefficients of $\partial_p S(t,x,p,\bz)$, i.e.,
\begin{align*}\partial_p\vv S(t,x,p)=(\partial_p\tilde{S}_1(t,x,p), \cdots, \partial_p\tilde{S}_K(t,x,p))^{T} \,,
\end{align*}
and the tensor $D(x)$ is given by
$$D_{l m k}(x)=\int_{I_{\bz}}\partial_x E(x,\bz)\tilde \psi_{l}(\bz)\tilde \psi_{m}(\bz)\tilde \psi_{k}(\bz)\pi(\bz)d\bz.$$
Note that the matrix ${\cal H}$ is symmetric.

The initial conditions have been discussed in subsection \ref{N2_3D}.

\subsubsection*{A.2. The fully discrete scheme gPC-SG-N1}
\label{full_new1}
To solve \eqref{SH gPC}, we use a time splitting procedure. We split the equation into four steps: two transport steps
(in $x$ and in $p$), a highly oscillatory part and a non-singular source part. We detail how we solve each
step in the sequel.

\textbf{Step 1} \quad
We solve the transport part in $x$ for each quantity $F^+, F^-, G, H$. For example,
to solve the transport in $F^+$,
$$\partial_t \vv F^{+}+p\partial_x \vv F^+=0, $$
we use a spectral method in space and an exact integration in time in the
Fourier space, that is
\begin{equation}\vv{\hat F}^{+, n+1}(\xi)=e^{-i p\xi\Delta t}\vv{\hat F}^{+, n}(\xi),\end{equation}
where $\xi$ denotes the Fourier space variable and $\hat F^+(\xi)$
the corresponding (discrete) Fourier transform. Finally,
$\vv F^{+, n+1}$ is obtained by the inverse Fourier transform. \\[2pt]

\textbf{Step 2}\quad
We solve the transport part in $p$ for $F^\pm$,
$$\partial_t \vv F^\pm \mp {\cal C}\partial_p \vv F^\pm=0. $$
We also use a spectral method in $p$ and an exact integration in time in the Fourier
space, that is
\begin{equation}\vv{\hat F}^{\pm,n+1}(\eta)=e^{\pm i\eta {\cal C}(x)\Delta t}\vv{\hat F}^{\pm,n}(\eta), \end{equation}
where $\eta$ is the Fourier velocity variable. Finally
$\vv F^{\pm, n+1}$ are obtained by the inverse Fourier transform.  \\[2pt]


\textbf{Step 3}\quad
We now solve the non-singular source part
\begin{align}
&\label{3a} \partial_t\vv F^+ = 2b^i (\vv G \cos\tau + \vv H \sin\tau),  \\[4pt]
&\label{3b} \partial_t\vv F^- = -2b^i (\vv G \cos\tau+ \vv H \sin\tau), \\[4pt]
&\label{3c} \partial_t\vv G = -b^i \vv F^+ \cos\tau+ b^i \vv F^-\cos\tau, \\[4pt]
&\label{3d} \partial_t \vv H = -b^i \vv F^+ \sin\tau+ b^i \vv F^-\sin\tau\,.
\end{align}
This linear system with time independent coefficients can be solved exactly.    \\[2pt]

\textbf{Step 4}\quad
Finally, we solve the highly oscillatory part
\begin{align*}
&\displaystyle\partial_t\vv F^+ = \frac{1}{\varepsilon}(-2{\cal C}+{\cal H})\partial_{\tau}\vv F^+, \qquad
\partial_t\vv F^- = \frac{1}{\varepsilon}(-2{\cal C}-{\cal H})\partial_{\tau}\vv F^-, \\[4pt]
&\displaystyle\partial_t\vv G = -\frac{2{\cal C}}{\varepsilon}\partial_{\tau}\vv G, \qquad\qquad\quad
\partial_t\vv H = -\frac{2{\cal C}}{\varepsilon}\partial_{\tau}\vv H.
\end{align*}
We use the Fourier transform $\hat F(\zeta)$ of $F(\tau)$ in the $\tau$ variable where $\zeta$ is the corresponding Fourier variable. Then, we have
\begin{align*}
\displaystyle
&\partial_t\vv{\hat F}^+=\frac{i\zeta}{\varepsilon}(-2{\cal C}+{\cal H})\vv{\hat F}^+, \qquad \partial_t\vv{\hat F}^-=\frac{i\zeta}{\varepsilon}(-2{\cal C}-{\cal H})\vv{\hat F}^-, \\[4pt]
&\partial_t\vv{\hat G}=\frac{i\zeta}{\varepsilon}(-2{\cal C})\vv{\hat G}, \qquad\qquad\partial_t\vv{\hat H}=\frac{i\zeta}{\varepsilon}(-2{\cal C})\vv{\hat H},
\end{align*}
which is solved by the backward Euler method in time,  \begin{align}
&\displaystyle\vv{\hat F}^{+,n+1}=\left[I-\frac{i\zeta}{\varepsilon}\Delta t(-2{\cal C}+{\cal H})\right]^{-1}\vv{\hat F}^{+,n}, \label{step4_eqn1}\\[2pt]
&\displaystyle\vv{\hat F}^{-,n+1}=\left[I-\frac{i\zeta}{\varepsilon}\Delta t(-2{\cal C}-{\cal H})\right]^{-1}\vv{\hat F}^{-,n},\label{step4_eqn2}\\[2pt]
&\displaystyle\vv{\hat G}^{n+1}=\left[I+2\frac{i\zeta}{\varepsilon}\Delta t {\cal C}\right]^{-1}\vv{\hat G}^n,\label{step4_eqn3}\\[2pt]
&\displaystyle\vv{\hat H}^{n+1}=\left[I+2\frac{i\zeta}{\varepsilon}\Delta t {\cal C}\right]^{-1}\vv{\hat H}^n. \label{step4_eqn4}
\end{align}
Notice that we have already solved the gPC coefficients vector $\vv S(t,x,p)$ from (\ref{S_3D}).
This allows to compute $\partial_p\vv S(t,x,p)$ spectrally (using Fourier in $p$) which is used in the expression of matrix
$\displaystyle\left(I-\frac{i\zeta}{\varepsilon}\Delta t(-2{\cal C}\pm{\cal H})\right)$. Note that the matrices in (\ref{step4_eqn1})-(\ref{step4_eqn4})
are invertible since they are symmetric and have real eigenvalues.

\subsubsection*{A.3 The gPC-SG-D scheme}
We briefly introduce the gPC-SG-D for the surface hopping model with random inputs. One first inserts the gPC expansions
$$(f^+, f^-, g, h)(t_n,x,p,\bz)\approx \sum_{j=1}^K (\tilde{f}^{+,n}_j, \tilde{f}^{-,n}_j, \tilde{g}^n_j, \tilde{h}^n_j)(x,p)\tilde{\psi}_j(\bz),   $$
into equation \eqref{F_3}. Then, denoting as previously $\vv g^n = (\tilde{g}^n_1, \tilde{g}^n_2, \cdots, \tilde{g}^n_K)$
(and the same notations for $\vv f^{+,n}, \vv f^{-,n}, \vv h^n$), we conduct the Galerkin projection to get
\begin{align*}
\displaystyle
&\vv g^{n+1}=\left[I+ (b^i)^2 (\Delta t)^2 I+\frac{\Delta t^2}{\varepsilon^2}\mathcal P\right]^{-1}\\[4pt]
&\qquad\quad\quad \left[\vv g^n +\Delta t\left(-b^i \vv f^{+,n}+b^i \vv f^{-,n} +\frac{2}{\varepsilon}\mathcal W \vv h^n-(b^i)^2\Delta t \vv g^n -\frac{\Delta t}{\varepsilon^2}\mathcal P
\vv g^n\right)\right] \,,
\end{align*}
where the matrices $\mathcal W_{K\times K}$ and $\mathcal P_{K\times K}$  are defined by
\begin{align*}
\displaystyle
& \mathcal W(x)_{ml }=\int_{I_{\bz}}E(x,\bz)\psi_m(\bz)\psi_l(\bz)\pi(\bz)d\bz, \quad 1\leq m, l\leq K, \\[4pt]
& \mathcal P(x)_{m l}=\int_{I_{\bz}}E(x,\bz)^2 \psi_m(\bz)\psi_l(\bz)\pi(\bz)d\bz, \quad 1\leq m, l\leq K.
\end{align*}
Since $\mathcal P$ is symmetric and positive definite, with positive eigenvalues, so does the matrix $\displaystyle\left(I+ (b^i)^2 (\Delta t)^2 I+\frac{\Delta t^2}{\varepsilon^2}\mathcal P\right)$, thus it is invertible.
Since $\vv g^{n+1}$ is known, we can compute $\vv f^{+,n+1}$, $\vv f^{-,n+1}$, $\vv h^{n+1}$ by gPC approximation of (\ref{eqn1}),
(\ref{eqn2}) and \eqref{eqn4} as
\begin{align*}
\displaystyle
& \vv f^{+,n+1}=\vv f^{+,n}+b^i\Delta t(\vv g^n+\vv g^{n+1}), \\[6pt]
& \vv f^{-,n+1}=\vv f^{-,n}-b^i\Delta t(\vv g^n+\vv g^{n+1}), \\[6pt]
& \vv h^{n+1}=\vv h^n-\frac{\Delta t}{\varepsilon}\mathcal W(\vv g^n+\vv g^{n+1}).
\end{align*}


\bibliographystyle{siam}
\bibliography{gPC_bibtex1}

\begin{thebibliography}{10}

\bibitem{Geim}
{\sc A.~Castro~Neto, F.~Guinea, N.~Peres, K.~Novoselov, and A.~Geim}, {\em The
  electronic properties of graphene}, Phys. Mod. Phys., 81 (2009),
  pp.~109--162.

\bibitem{CJL}
{\sc L.~Chai, S.~Jin, and Q.~Li}, {\em Semi-classical models for the
  {S}chr\"odinger equation with periodic potentials and band crossings}, Kinet.
  Relat. Models, 6 (2013), pp.~505--532.

\bibitem{CQJ}
{\sc L.~Chai, S.~Jin, Q.~Li, and O.~Morandi}, {\em A multiband semiclassical
  model for surface hopping quantum dynamics}, Multiscale Model. Simul., 13
  (2015), pp.~205--230.

\bibitem{PCL}
{\sc P.~Chartier, N.~Crouseilles, M.~Lemou, and F.~M\'{e}hats}, {\em Uniformly
  accurate numerical schemes for highly-oscillatory {K}lein-{G}ordon and
  nonlinear {S}chr\"{o}dinger equations}, Numer. Math., 129 (2015),
  pp.~211--250.

\bibitem{Choi-Liu}
{\sc H.~Choi and J.~Liu}, {\em The reconstruction of upwind fluxes for
  conservation laws: Its behavior in dynamic and steady state calculations}, J.
  Comput. Phys., 144 (1998), pp.~237--256.

\bibitem{NGO}
{\sc N.~Crouseilles, S.~Jin, and M.~Lemou}, {\em Nonlinear geometric optics
  method based multi-scale numerical schemes for highly-oscillatory transport
  equations}, arXiv:1605.09676,  (2016).

\bibitem{FCL}
{\sc N.~Crouseilles, M.~Lemou, and F.~M\'{e}hats}, {\em Asymptotic preserving
  schemes for highly-oscillatory {V}lasov-{P}oisson equations}, J. Comput.
  Phys., 248 (2013), pp.~287--308.

\bibitem{ER2}
{\sc B.~Engquist and O.~Runborg}, {\em Computational high frequency wave
  propagation}, Acta Numer., 12 (2003), pp.~181--266.

\bibitem{FW}
{\sc C.~L. Fefferman and M.~I. Weinstein}, {\em Honeycomb lattice potentials
  and {D}irac points}, J. Amer. Math. Soc., 25 (2012), pp.~1169--1220.

\bibitem{GS}
{\sc R.~G. Ghanem and P.~D. Spanos}, {\em Stochastic finite elements: A
  spectral approach}, Springer-Verlag, New York,  (1991).

\bibitem{GX}
{\sc D.~Gottlieb and D.~Xiu}, {\em Galerkin method for wave equations with
  uncertain coefficients}, Commun. Comput. Phys., 3 (2008), pp.~505--518.

\bibitem{GWZ}
{\sc M.~D. Gunzburger, C.~G. Webster, and G.~Zhang}, {\em Stochastic finite
  element methods for partial differential equations with random input data},
  Acta Numer., 23 (2014), pp.~521--650.

\bibitem{LMK}
{\sc O.~P. Le~Ma{\^{\i}}tre and O.~M. Knio}, {\em Spectral methods for
  uncertainty quantification, with applications to computational fluid
  dynamics}, Scientific Computation, Springer, New York,  (2010).

\bibitem{Morandi2}
{\sc O.~Morandi}, {\em Effective classical {L}iouville-like evolution equation
  for the quantum phase-space dynamics}, J. Phys. A, 43 (2010).

\bibitem{Morandi}
{\sc O.~Morandi and F.~Sch{\"u}rrer}, {\em Wigner model for {K}lein tunneling
  in graphene}, J. Phys. A, 44 (2011).

\bibitem{STW}
{\sc J.~Shen, T.~Tang, and L.-L. Wang}, {\em Spectral methods: Algorithms,
  analysis and applications}, Springer, Heidelberg, 41 (2011).

\bibitem{Tully2}
{\sc J.~Tully}, {\em Molecular dynamics with electonic transitions}, J. Chem.
  Phys., 93 (1990), pp.~1061--1071.

\bibitem{Tully1}
{\sc J.~Tully and R.~Preston}, {\em Trajectory surface hopping approach to
  nonadiabatic molecular collisions: the reaction of $h^+$ with $d_2$}, J.
  Chem. Phys., 55 (1971), pp.~562--572.

\bibitem{WuNiu}
{\sc B.~Wu and Q.~Niu}, {\em Nonlinear {L}andau-{Z}ener tunneling}, Phys. Rev.
  A., 61 (2000).

\bibitem{Xiubook}
{\sc D.~Xiu}, {\em Numerical Methods for Stochastic Computations}, Princeton
  University Press, Princeton, New Jersey, 2010.

\bibitem{XH}
{\sc D.~Xiu and J.~S. Hesthaven}, {\em High-order collocation methods for
  differential equations with random inputs}, SIAM J. Sci. Comput., 27 (2005),
  pp.~1118--1139.

\bibitem{XK}
{\sc D.~Xiu and G.~E. Karniadakis}, {\em The {W}iener-{A}skey polynomial chaos
  for stochastic differential equations}, SIAM J. Sci. Comput., 24 (2002),
  pp.~619--644.

\end{thebibliography}
\end{document}